\documentclass[10pt,a4paper]{article}

\usepackage{amsmath}
\usepackage{graphicx,latexsym,euscript,makeidx,color}
\usepackage{amsmath,amsfonts,amssymb,amsthm,mathrsfs}
\usepackage{enumerate}

\usepackage{graphicx,colordvi}

\usepackage{amssymb}
\usepackage{amsmath}
\usepackage{amsfonts}
\usepackage{graphicx}
\usepackage{graphicx}          
\usepackage{color}

\usepackage[colorlinks,linkcolor=blue,anchorcolor=green,citecolor=red]{hyperref}

\usepackage{geometry}
\def\dbE{{\mathbb{E}}}

\font\tenbb=msbm10 \font\sevenbb=msbm7 \font\fivebb=msbm5
\newfam\bbfam
\scriptscriptfont\bbfam=\fivebb \textfont\bbfam=\tenbb
\scriptfont\bbfam=\sevenbb

\newtheorem{theorem}{\indent Theorem}[section]
\newtheorem{definition}[theorem]{\indent Definition}

\newtheorem{corollary}[theorem]{\indent Corollary}
\newtheorem{lemma}[theorem]{\indent Lemma}
\newtheorem{remark}[theorem]{\indent Remark}

\makeatletter
   
   \@addtoreset{equation}{section}
\makeatother

\allowdisplaybreaks[4]

\begin{document}

\title{\bf Mixed Equilibrium Solution of Time-Inconsistent Stochastic LQ Problem
\thanks{This work is supported in part by the National Natural Science Foundation of China (61227902, 11471242, 61773222), the National Key Basic Research Program (973 Program) of China (2014CB845301), Hong Kong RGC grants 15224215 and 15255416.}
}
\author{Yuan-Hua Ni\thanks{College of Artificial Intelligence, Nankai University, Tianjin 300350, P.R. China. Email: {\tt yhni@nankai.edu.cn}.}~~~~~~Xun Li\thanks{Department of Applied Mathematics, The Hong Kong
Polytechnic University, Kowloon, Hong Kong, P.R. China. Email: {\tt malixun@polyu.edu.hk}.}~~~~~~Ji-Feng Zhang\thanks{Key Laboratory of Systems and Control,
Institute of Systems Science, Academy of Mathematics and Systems Science,
Chinese Academy of Sciences, Beijing 100190; School of
Mathematical Sciences, University of Chinese Academy of Sciences, Beijing
100049, P.R. China. Email: {\tt jif@iss.ac.cn}.}~~~~~~Miroslav Krstic\thanks{Department of Mechanical and Aerospace Engineering, University of California, San Diego, USA. Email: {\tt krstic@ucsd.edu}.}}
\maketitle

{\bf Abstract:} In this paper, we propose a novel equilibrium solution notion for the time-inconsistent stochastic linear-quadratic optimal control problem. {This notion is called the mixed equilibrium solution, which} consists of two parts: a pure-feedback-strategy part and an open-loop-control part.
When the pure-feedback-strategy part is zero or the open-loop-control part does not depend on the initial state, the mixed equilibrium solution reduces to the open-loop equilibrium control and the feedback equilibrium strategy, respectively.
Using a maximum-principle-like methodology with forward-backward stochastic difference equations,
a necessary and sufficient condition is established to characterize the existence of a mixed equilibrium solution. Then, by decoupling the forward-backward stochastic difference equations, three sets of difference equations, which together portray the existence of a mixed equilibrium solution, are obtained.
Moreover, the case with a fixed time-state initial pair and the case with all the initial pairs are separately investigated.
Furthermore,
an example is constructed to show that the mixed
equilibrium solution exists for all the initial pairs, although neither the open-loop equilibrium control nor the feedback
equilibrium strategy exists for some initial pairs.


\section{Introduction}

In this paper, we consider a class of mean-field stochastic linear-quadratic (LQ, for short) control problems. The system dynamics are described by the following discrete-time stochastic difference equation (S$\Delta$E, for short)
\begin{eqnarray}\label{system-1}
\left\{\begin{array}{l}
X^t_{k+1}=\big{(}A_{t,k}X^t_k+\bar{A}_{t,k}\mathbb{E}_tX^t_k+B_{t,k}u_k+\bar{B}_{t,k}
\mathbb{E}_tu_k+f_{t,k}\big{)}\\[1mm]
\hphantom{X^t_{k+1}=}+\sum_{i=1}^p\big{(}C^i_{t,k}X^t_k+\bar{C}^i_{t,k}
\mathbb{E}_tX^t_k+D^i_{t,k}u_k+\bar{D}^i_{t,k}\mathbb{E}_tu_k
+d^i_{t,k}\big{)}w^i_k, \\[1mm]
X^t_t=x,~~k\in  \mathbb{T}_t,~~t\in \mathbb{T},
\end{array}
\right.
\end{eqnarray}
where
$\mathbb{T}=\{0,\dots,N-1\}$, $\mathbb{T}_t=\{t,\cdots,N-1\}$ and
$A_{t,k},\bar{A}_{t,k},C^i_{t,k},\bar{C}^i_{t,k}\in \mathbb{R}^{n\times
n}$, $B_{t,k},\bar{B}_{t,k},D^i_{t,k},\bar{D}^i_{t,k}\in
\mathbb{R}^{n\times m}, f_{t,k},
d^i_{t,k}\in \mathbb{R}^n$ are deterministic matrices, and $\{X^t_{k}, k\in
\widetilde{{\mathbb{T}}}_t\}\triangleq X^t$ and $\{u_k, k\in
\mathbb{T}_t\}\triangleq u$ with $\widetilde{\mathbb{T}}_t=\{t,1,...,N\}$ are the state process and control
process, respectively. The noise $\{w_k, k\in
\mathbb{T}\}$ is assumed to be a vector-valued martingale difference sequence
defined on a probability space $(\Omega, \mathcal{F}, P)$ with
\begin{eqnarray}\label{w-moment}
\mathbb{E}_{k}[w_{k}]=0,~~\mathbb{E}_{k}[(w_{k})^2]=\Delta_k,~~k\in \mathbb{T},
\end{eqnarray}
where $\Delta_{k}=(\delta_{k}^{ij})_{p\times p}$, $k\in \mathbb{T}$, are assumed to be deterministic.
$\mathbb{E}_t[\,\cdot\,]$ in (\ref{system-1}) denotes the conditional mathematical expectation
$\mathbb{E}[\,\cdot\,|\mathcal{F}_t]$, where $\mathcal{F}_t$ is defined as $\sigma\{w_l, l=0, 1,\cdots,t-1\}$ and $\mathcal{F}_{0}$ is understood as
$\{\emptyset, \Omega\}$, and $\mathbb{E}_{k}[\,\cdot\,]$ in (\ref{w-moment}) is similarly defined.
In (\ref{system-1}),
$x$ belongs to $l^2_\mathcal{F}(t; \mathbb{R}^n)$, which is defined as
$\big{\{}\zeta \in \mathbb{R}^n\,\big{|}$ $\zeta$ is $\mathcal{F}_t$-measurable, $\mathbb{E}|\zeta|^2<\infty \big{\}}$.
%
%
%
We introduce the cost function
\begin{eqnarray}\label{cost-1}
&&\hspace{-3.5em}J(t,x;u)
 =\sum_{k=t}^{N-1}\mathbb{E}_t\Big{\{}(X_k^t)^TQ_{t,k}X^t_k +(\mathbb{E}_tX_k^t)^T
 \bar{Q}_{t,k}\mathbb{E}_tX^t_k+ u_k^TR_{t,k}u_k+(\mathbb{E}_tu_k)^T \bar{R}_{t,k}
 \mathbb{E}_tu_k
\nonumber\\
&&\hspace{-3.5em}\hphantom{J(t,x;u)=}  +2q_{t,k}^T X^t_k+ 2\rho_{t,k}^T u_k\Big{\}}+\mathbb{E}_t\big{[}(X_N^t)^TG_tX^t_N\big{]}
+(\mathbb{E}_tX^t_N)^T \bar{G}_t\mathbb{E}_tX^t_{N}+2(F_tx+g_t)^T \mathbb{E}_tX^t_N, 
\end{eqnarray}
where $Q_{t,k}, \bar{Q}_{t,k},R_{t,k},\bar{R}_{t,k}, k\in
\mathbb{T}_t$, $G_t, \bar{G}_t$ are deterministic symmetric matrices
of appropriate dimensions, and $q_{t,k}, \rho_{t,k},  k\in
\mathbb{T}_t, g_t$ are deterministic vectors.
Let
$$ 
l^2_\mathcal{F}(\mathbb{T}_{t}; \mathbb{R}^m)=
\big\{\nu=\{\nu_k, k\in \mathbb{T}_{t}\}\, \big{|}\,
\nu_k\mbox{ is }\mathcal{F}_k\mbox{-measurable},
\mathbb{E}|\nu_k|^2 < \infty,  k\in\mathbb{T}_{t}\big\}.$$
%
Then, we pose the
following optimal control problem.

\textbf{Problem (LQ).}  \emph{For the time-state initial pair $(t,x)$, find
${u}^*\in l^2_\mathcal{F}(\mathbb{T}_t; \mathbb{R}^m)$ such that
\begin{eqnarray*}\label{Problem-LQ}
J(t,x;{u}^*) = \inf_{u\in l^2_\mathcal{F}(\mathbb{T}_t; \mathbb{R}^m)}J(t,x;u).
\end{eqnarray*}
}

Compared with the standard stochastic LQ problems, Problem (LQ) has three unconventional features. First, the cost weighting and system matrices depend explicitly on the initial time $t$. Second, the term $2(F_tx+g_t)^T\mathbb{E}_tX_N^t$ makes $J(t,x;u)$ a state-dependent (or rank-dependent) utility. Third, $J(t,x;u)$ contains nonlinear terms of the conditional expectation of state and control. These three features are deeply rooted in the fields of economics and finance. The first feature is an abstraction of the general discounting functions; see \cite{Bjork,Krusell}
for examples of hyperbolic discounting and quasi-geometric discounting.
The second feature is of rank-dependent utility \cite{Bjork-2}, and a notable example of the third feature is the mean-variance utility \cite{Basak,Bjork,{Cui-2017},{Li-duan-2013},Li-Duan,Markowitz-1}.
It is known that any of the three features will ruin the time consistency of the optimal
control, namely, Bellman's principle of optimality will no longer work for Problem (LQ).

Problems with the nonlinear term of the conditional expectation (in the cost functional) are classified as mean-field stochastic optimal control problems \cite{Yong-2013}.
Realizing the time inconsistency (called nonseparability there), Li and Ng \cite{Li-Duan} used an embedding scheme to derive the optimal policy for the multi-period mean-variance portfolio selection.
Note that the optimal policy of \cite{Li-Duan} is with respect to the initial pair, that is, it is an optimal policy only when viewed at the initial time.
This derivation is now called the pre-committed optimal solution. 
%
%
However, we find that a pre-committed optimal control (with respect to an initial pair) will no longer  serve as an optimal control for an intertemporal initial pair.
Although the pre-committed optimal solution is of some practical and theoretical value, it neglects and does not
fully address the time inconsistency.
Another approach is to handle the time inconsistency in a
dynamic manner, by seeking time-consistent equilibrium solutions instead of a pre-committed optimal control; this has mainly been motivated by
practical applications in economics and finance, and has recently
attracted considerable research interest.

The qualitative analysis of time inconsistency can be traced back to the ideas of the father of free market economics and moral philosopher Adam Smith \cite{Smith}. In 1955, Strotz \cite{Strotz} gave the first quantitative formulation of time inconsistency and studied the general discounting problem. His approach successfully tackled time inconsistency using a lead-follower game with a hierarchical structure.
Inspired by Strotz,
hundreds of works have sought to tackle practical problems in economics and finance by focusing on the time inconsistency of dynamic
systems described by ordinary difference or differential equations;
see, for example,
\cite{Ekeland,Ekland-2,Goldman,Krusell,Laibson,Palacios} and
the references therein. Unfortunately, as Ekeland
\cite{Ekeland,Ekland-2} pointed out, it is hard to prove the existence of Strotz's equilibrium policy. Therefore, it is necessary and of great
importance to develop a general theory of time-inconsistent
control. 
In recent years, this topic
has attracted considerable attention from the theoretical control
community; see, for example,
\cite{Bjork,Hu-jin-Zhou,Hu-jin-zhou-2,Wang-haiyang,Yong-1,Yong-0,Yong-2013}
and the references therein.

With respect to the time-inconsistent LQ problems, two kinds of
time-consistent equilibrium solutions have been investigated, namely, the
open-loop equilibrium control and the closed-loop equilibrium
strategy \cite{Hu-jin-Zhou,Hu-jin-zhou-2,Yong-1,Yong-0,Yong-2013}. The two formulations are investigated
separately because in dynamic game theory, open-loop control
differs significantly from the closed-loop strategy
\cite{Basar,Yong-game-2015}. To compare, the aim of open-loop formulation is to find an open-loop equilibrium
``control," while the ``strategy" is the object of closed-loop
formulation. 
%
%
Yong further developed Strotz's equilibrium solution \cite{Strotz}, which is essentially a
closed-loop equilibrium strategy, into the LQ optimal control \cite{Yong-1,Yong-2013}
and the nonlinear optimal control \cite{Yong-0,Yong-2,Wei-Yu-Yong2017}. The open-loop equilibrium control has been extensively studied by
Hu-Jin-Zhou \cite{Hu-jin-Zhou,Hu-jin-zhou-2}, Yong
\cite{Yong-2013}, Ni-Zhang-Krstic \cite{Ni-Zhang-Krstic2017}, and Qi-Zhang\cite{Zhang-huanshui}. In particularly, the closed-loop
formulation can be viewed as an extension of Bellman's dynamic
programming, and the corresponding equilibrium strategy ({if it
exists}) is constructed by a backward procedure
\cite{Yong-1,Yong-2,Yong-0,Yong-2013}. Differently, the open-loop
equilibrium control is characterized via a maximum-principle-like
methodology \cite{Hu-jin-Zhou,Hu-jin-zhou-2,Ni-Zhang-Krstic2017}.

It is well known that the aim of portfolio selection is to seek the best allocation of wealth among a basket of securities.
The (single-period) mean-variance formulation initiated by Markowitz \cite{Markowitz-1} is the cornerstone of modern portfolio theory and is widely used in both academic studies and the financial industry.
The multi-period mean-variance portfolio selection is the natural extension of \cite{Markowitz-1}, which has been extensively studied.
Li-Ng \cite{Li-Duan} and Zhou-Li \cite{zhou-xunyu-2000} were the first to report the analytical pre-commitment optimal policies for the discrete-time case and the continuous-time case, respectively.
In fact, the multi-period mean-variance portfolio selection problem, which is a particular example of time-inconsistent problem, stimulated the recent developments in time-inconsistent problems and  the revisits to multi-period mean-variance portfolio selection \cite{Basak,{Bjork-2},Cui-2017,Cui-2012,Hu-jin-Zhou,Hu-jin-zhou-2}.

In this paper, we examine the aforementioned Problem (LQ). In Section \ref{Section-Mixed}, we introduce the mixed equilibrium solution to Problem (LQ). The solution contains two different parts: a pure-feedback-strategy part and an open-loop-control part. By letting the open-loop-control part be independent of the initial state or the pure-feedback-strategy part be zero, the corresponding mixed equilibrium solution is reduced to a linear feedback equilibrium strategy and open-loop equilibrium control, respectively.
Section \ref{Section-characterization} characterizes the mixed equilibrium solution using a maximum-principle-like
methodology with convexity, stationarity, and forward-backward stochastic difference equations (FBS$\Delta$Es). It is shown that the convexity and stationarity conditions can be equivalently characterized via solutions to three sets of difference equations. Based on the results for the mixed equilibrium solution, we then obtain the results for the open-loop equilibrium control and linear feedback equilibrium strategy (with respect to a fixed initial pair).
{For the case with all the initial pairs, conditions in terms of solvability of three sets of difference equations are given to ensure the existence of mixed equilibrium solution}. These conditions are necessary and sufficient to determine the open-loop equilibrium control and linear feedback equilibrium strategy.
Interestingly, for all of the initial pairs, the existence of general feedback equilibrium strategy is shown to be equivalent to the existence of linear feedback equilibrium strategy, which can be obtained by a backward procedure. Furthermore, the backward procedure works only when the feedback equilibrium strategy exists for all of the initial pairs, and cannot be applied to the case where we know only of the existence of a feedback equilibrium strategy for a fixed initial pair.
Section \ref{Section-example} gives an example to illustrate the developed theory.
Finally, in Section \ref{Section-Summary}, we discuss future topics that are worth investigating.
%
%

This paper makes the following novelties.
\begin{itemize}
\item[$\bullet$]Most of the existing results for time-inconsistent LQ problems are for the continuous-time case \cite{Hu-jin-Zhou,Hu-jin-zhou-2,Wei-Yu-Yong2017,Yong-1,Yong-0,Yong-2013}. The discrete-time multi-period mean-variance portfolio selection problem is a notable example of Problem (LQ), and its investigation calls for the development of general theory of discrete-time time-inconsistent LQ optimal control. Furthermore, the model and methodology developed in this paper are more general than those in \cite{Ni-Zhang-Krstic2017}.

\item[$\bullet$]The notion of mixed equilibrium solution is introduced, and it seems that no similar notion has been reported for time-inconsistent optimal control. Necessary and sufficient conditions are established
to characterize a pair of pure-feedback strategy and open-loop control as a mixed equilibrium solution
(for a time-state initial pair). Using the notion of mixed equilibrium solution, the conditions to equivalently ensure the existence of an open-loop equilibrium control and a linear feedback equilibrium strategy can be simultaneously obtained. In other words, we can investigate the two equilibrium solutions in a unified way.

\vspace{0.15em}

Importantly, the mixed equilibrium solution is not a hollow concept. In Section \ref{Section-example}, it is shown that neither the open-loop equilibrium control nor the feedback equilibrium strategy exists for the initial pair $(t,x)$ with $t=0,1$ and $x\in l^2_{\mathcal{F}}(t; \mathbb{R}^2)$, although we are able to construct 10 mixed equilibrium solutions. Therefore, it is necessary to study the mixed equilibrium solution, which gives us more flexibility to deal with the time-inconsistent optimal control.

\end{itemize}

The work of \cite{Ni-Li-Zhang-Krstic2018} serves as a companion to this paper in terms of testing our developed theory and pursuing the solvability of the multi-period mean-variance portfolio selection problem. The non-degenerate assumption was removed in \cite{Ni-Li-Zhang-Krstic2018}, which is popular in the literature on multi-period mean-variance portfolio selection. Neat conditions have been obtained in \cite{Ni-Li-Zhang-Krstic2018} to characterize the existence of the equilibrium solutions.
To emphasize the dependence on the initial pair, Problem (LQ) for the initial pair $(t,x)$ is denoted as Problem (LQ)$_{tx}$ throughout this paper.
Furthermore, for notational simplicity, we denote in  this paper
\begin{eqnarray*}
&&\mathcal{A}_{t,k}=A_{t,k}+\bar{A}_{t,k},~~\mathcal{B}_{t,k}=B_{t,k}+\bar{B}_{t,k}, ~~\mathcal{C}^i_{t,k}=C^i_{t,k}+\bar{C}^i_{t,k},~~\mathcal{D}^i_{t,k}=D^i_{t,k}+\bar{D}^i_{t,k},\\[1mm]
&&\mathcal{Q}_{t,k}=Q_{t,k}+\bar{Q}_{t,k},~~\mathcal{R}_{t,k}=R_{t,k}+\bar{R}_{t,k},~~\mathcal{G}_{t}=G_{t}+\bar{G}_{t},~~t\in \mathbb{T},~~  k\in \mathbb{T}_t.
\end{eqnarray*}

\section{Mixed equilibrium solution}\label{Section-Mixed}

Before introducing the mixed equilibrium solution, we give the definition of feedback equilibrium strategy.
By a strategy, we mean a decision rule that a
controller uses to select a control action based on the available
information set. Mathematically, a strategy is a mapping or an operator defined
on the information set. Substituting the available information
into a strategy, we obtain the open-loop value or realization of
this strategy.
%
%
%

\begin{definition}
i). At stage $k\in \mathbb{T}_t$, a function $f_k(\cdot)$ is called an admissible feedback strategy (or simply a feedback strategy) at $k$ if for $\zeta\in l^2_{\mathcal{F}}(k; \mathbb{R}^n)$, $f_k(\zeta)\in l^2_{\mathcal{F}}(k; \mathbb{R}^m)$. The set of this type of feedback strategies is denoted by $\mathbb{F}_k$, and $\mathbb{F}_t\times \cdots\times \mathbb{F}_{N-1}$ is denoted by $\mathbb{F}_{\mathbb{T}_t}$.

ii). Let $f=(f_t,...,f_{N-1})\in \mathbb{F}_{\mathbb{T}_t}$. For $k\in \mathbb{T}_t$ and $\zeta\in l^2_{\mathcal{F}}(k; \mathbb{R}^n)$, $f_k(\zeta)$ can be divided into two parts, namely, $f_k(\zeta)=f^c_k+f^p_k(\zeta)$, where $f^c_k=f_k(0)$ is the inhomogeneous part and the remainder $f_k^p(\cdot)$ is the pure-feedback-strategy part of $f_k$. Furthermore, $(f_t^p,...,f_{N-1}^p)$ is called a pure-feedback strategy.

\end{definition}
%
%

\begin{definition}\label{definition-closed-loop-2}
i). A strategy $\psi\in \mathbb{F}_{\mathbb{T}_t}$ is called a feedback equilibrium strategy of Problem (LQ)$_{t,x}$, if the following two points hold:

\begin{itemize}

\item[a)] $\psi$ does not depend on $x$;

\item[b)] For any $k\in \mathbb{T}_t$ and any $u_k\in l^2_\mathcal{F}(k; \mathbb{R}^m)$, it holds that
\begin{eqnarray}\label{defi-closed-loop}
J\big{(}k, X_k^{t,x,*}; (\psi\cdot X^{k,\psi})|_{\mathbb{T}_k}\big{)}\leq J\big{(}k, X_k^{t,x,*}; (u_k,(\psi\cdot X^{k,u_k,\psi})|_{\mathbb{T}_{k+1}})\big{)}.
\end{eqnarray}
In (\ref{defi-closed-loop}), $(\psi\cdot X^{k,\psi})|_{\mathbb{T}_k}$ and $(\psi\cdot X^{k,u_k,\psi})|_{\mathbb{T}_{k+1}}$ (with $\mathbb{T}_k=\{k,...,N-1\}, \mathbb{T}_{k+1}=\{k+1,...,N-1\}$) are given by
\begin{eqnarray*}
&&(\psi\cdot X^{k,\psi})|_{\mathbb{T}_k}=\big{(}\psi_k(X^{k,\psi}_k),...,\psi_{N-1}(X^{k,\psi}_{N-1})\big{)},\\
&&(\psi\cdot X^{k,u_k,\psi})|_{\mathbb{T}_{k+1}}=\big{(}\psi_{k+1}(X^{k,u_k,\psi}_{k+1}),...,\psi_{N-1}(X^{k,u_k,\psi}_{N-1})\big{)},
\end{eqnarray*}
where
$X^{k,\psi}, X^{k,u_k,\psi}$ are as follows
\begin{eqnarray}
%
%
&&\hspace{-3em}\label{Sec3:Defini-closed-X---2}
\left\{\begin{array}{l}
X^{k,\psi}_{\ell+1} =A_{k,\ell}X^{k,\psi}_\ell+B_{k,\ell}\psi_\ell(X^{k,\psi}_\ell)+\bar{A}_{k,\ell}\mathbb{E}_kX^{k,\psi}_\ell+\bar{B}_{k,\ell}\mathbb{E}_k\psi_\ell(X^{k,\psi}_\ell)+f_{k,\ell}\\[1mm]
\hphantom{X^{k,\psi}_{\ell+1} =}+\sum_{i=1}^p\Big{[}C^i_{k,\ell}X^{k,\psi}_\ell+D^i_{k,\ell}\psi_\ell(X^{k,\psi}_\ell)+\bar{C}^i_{k,\ell}\mathbb{E}_kX^{k,\psi}_\ell+\bar{D}^i_{k,\ell}\mathbb{E}_k\psi_\ell(X^{k,\psi}_\ell)+d^i_{k,\ell}
\Big{]} w^i_\ell,\\[1mm]
X^{k,\psi}_{k} = X^{t,x,*}_k,~~\ell\in \mathbb{T}_k,
\end{array}
\right.\\
&&\label{sec3:X-u---2}
\hspace{-3em}\left\{\begin{array}{l}
X^{k,u_k,\psi}_{\ell+1} =A_{k,\ell}X^{k,u_k,\psi}_\ell+B_{k,\ell}\psi_\ell(X^{k,u_k,\psi}_\ell)+\bar{A}_{k,\ell}\mathbb{E}_kX^{k,u_k,\psi}_\ell +\bar{B}_{k,\ell}\mathbb{E}_k\psi_\ell(X^{k,u_k,\psi}_\ell)+f_{k,\ell}\\[1mm]
\hphantom{X^{k,u_k,\psi}_{\ell+1} =}+\sum_{i=1}^p\Big{[}C^i_{k,\ell}X^{k,u_k,\psi}_\ell+D^i_{k,\ell}\psi_\ell(X^{k,u_k,\psi}_\ell)+\bar{C}^i_{k,\ell}\mathbb{E}_kX^{k,u_k,\psi}_\ell \\[1mm]
\hphantom{X^{k,u_k,\psi}_{\ell+1} =} +\bar{D}^i_{k,\ell}\mathbb{E}_k\psi_\ell(X^{k,u_k,\psi}_\ell)+d^i_{k,\ell}
\Big{]} w^i_\ell,\\[1mm]
X^{k,u_k,\psi}_{k+1} =\big{[}\mathcal{A}_{k,k}X^{k,u_k,\psi}_k+\mathcal{B}_{k,k}u_k+f_{k,k}\big{]}+\sum_{i=1}^p\big{[}\mathcal{C}^i_{k,k}X^{k,u_k,\psi}_k+\mathcal{D}^i_{k,k}u_k+d^i_{k,k}\big{]}w^i_k,\\[1mm]
X^{k,u_k,\psi}_{k} = X^{t,x,*}_k,~~\ell\in \mathbb{T}_{k+1}.
\end{array}\right.
\end{eqnarray}
Furthermore, in (\ref{defi-closed-loop}), (\ref{Sec3:Defini-closed-X---2}), and (\ref{sec3:X-u---2}), $X_k^{t,x,*}$ is computed via
\begin{eqnarray*}
\left\{\begin{array}{l}
X^{t,x,*}_{k+1} =\big{[}\mathcal{A}_{k,k}X^{t,x,*}_{k}+\mathcal{B}_{k,k}\psi_k(X^{t,x,*}_k)+f_{k,k}\big{]}+\sum_{i=1}^p\big{[}\mathcal{C}^i_{k,k}X^{t,x,*}_{k}+\mathcal{D}^i_{k,k}\psi_k(X^{t,x,*}_k)+d^i_{k,k}\big{]} w^i_k,\\[1mm]
X^{t,x,*}_{t} = x,~~k\in \mathbb{T}_t.
\end{array}
\right.
\end{eqnarray*}

\end{itemize}

ii). Let $(\Psi, \gamma)\in l^2(\mathbb{T}_t; \mathbb{R}^{m\times n})\times  l^2_{\mathcal{F}}(\mathbb{T}_t;\mathbb{R}^m)$ with
$$ 
l^2(\mathbb{T}_t; \mathbb{R}^{m\times n})=\Big{\{} \nu=\{\nu_k, k\in \mathbb{T}_t\}\Big{|}v_k\in \mathbb{R}^{m\times n} ~\mbox{is deterministic}, |\nu_k|^2<\infty, k\in \mathbb{T}_t
\Big{\}}.
$$ 
If $\Psi$ and $\gamma$ do not depend on $x$, and $\psi$ of i) is equal to $(\Psi, \gamma)$, namely,
$\psi_k(\xi)=\Psi_k\xi+\gamma_k$,
$k\in \mathbb{T}_t$, $\xi\in l^2_{\mathcal{F}}(k;\mathbb{R}^n)$,
then $(\Psi, \gamma)$ is called a linear feedback equilibrium strategy of Problem (LQ)$_{tx}$.

\end{definition}

\begin{definition}
A control $u^{t,x}\in l^2_{\mathcal{F}}(\mathbb{T}_t;\mathbb{R}^m)$ is called an open-loop equilibrium control of Problem (LQ)$_{tx}$, if
\begin{eqnarray}
J\big{(}k, X_k^{t,x,*}; u^{t,x}|_{\mathbb{T}_k}\big{)}\leq J\big{(}k, X_k^{t,x,*}; (u_k,u^{t,x}|_{\mathbb{T}_{k+1}})\big{)}
\end{eqnarray}
holds for any $k\in \mathbb{T}_t$ and any $u_k\in
l^2_\mathcal{F}(k; \mathbb{R}^m)$.
Here, $u^{t,x}|_{\mathbb{T}_k}$ and
$u^{t,x}|_{\mathbb{T}_{k+1}}$  are the restrictions of $u^{t,x}$
on $\mathbb{T}_k$ and $\mathbb{T}_{k+1}$, respectively; and $X_k^{t,x,*}$ is computed via
\begin{eqnarray*}
\left\{
\begin{array}{l}
X^{t,x,*}_{k+1} = \big{[}\mathcal{A}_{k,k}X^{t,x,*}_{k}+\mathcal{B}_{k,k}u^{t,x}_k+f_{k,k}\big{]}+\sum_{i=1}^p\big{[}\mathcal{C}^i_{k,k}X^{t,x,*}_{k}+\mathcal{D}^i_{k,k}u^{t,x}_k+d^i_{k,k}\big{]}w^i_k,\\[1mm]
X^{t,x,*}_{t} = x,~~ k \in  \mathbb{T}_t.
\end{array}
\right.
\end{eqnarray*}

\end{definition}

\begin{definition}\label{definition-linear feedback}
i). A pair $(\Phi, v^{t,x}) \in l^2(\mathbb{T}_t; \mathbb{R}^{m\times n})\times  l^2_{\mathcal{F}}(\mathbb{T}_t;\mathbb{R}^m)$ is called a mixed equilibrium solution of Problem (LQ)$_{tx}$, if the following two points hold:

\begin{itemize}

\item[a)] $\Phi$ does not depend on $x$, and $v^{t,x}$ depends on $x$;

\item[b)] For any $k\in \mathbb{T}_t$ and any $u_k\in l^2_\mathcal{F}(k; \mathbb{R}^m)$, it holds that
\begin{eqnarray}\label{defi-linear feedback}
J\big{(}k, X_k^{t,x,*}; (\Phi\cdot X^{k,\Phi}+v^{t,x})|_{\mathbb{T}_k}\big{)}\leq J\big{(}k, X_k^{t,x,*}; (u_k,(\Phi \cdot X^{k,u_k,\Phi}+v^{t,x})|_{\mathbb{T}_{k+1}})\big{)}.
\end{eqnarray}
In (\ref{defi-linear feedback}), $(\Phi\cdot X^{k,\Phi}+v^{t,x})|_{\mathbb{T}_k}$ and $(\Phi\cdot X^{k,u_k,\Phi}+v^{t,x})|_{\mathbb{T}_{k+1}}$ are given, respectively, by
\begin{eqnarray*}
&&\big{(}\Phi_k X_k^{k,\Phi}+v^{t,x}_k,\cdots, \Phi_{N-1} X_{N-1}^{k,\Phi}+v^{t,x}_{N-1}\big{)},\\
&&\big{(}\Phi_{k+1} X_{k+1}^{k,u_k,\Phi}+v^{t,x}_{k+1},\cdots, \Phi_{N-1} X_{N-1}^{k,u_k,\Phi}+v^{t,x}_{N-1}\big{)},
\end{eqnarray*}
where $X^{k,\Phi}$, $X^{k,u_k,\Phi}$ are defined by
\begin{eqnarray}
%
%
&&\hspace{-2em}\label{Sec2:Defini-closed-X-1}
\left\{\begin{array}{l}
X^{k,\Phi}_{\ell+1} =\Big{[}\big{(}A_{k,\ell}+B_{k,\ell}\Phi_\ell\big{)} X^{k,\Phi}_{\ell}+\big{(}\bar{A}_{k,\ell}+\bar{B}_{k,\ell}\Phi_\ell\big{)} \mathbb{E}_kX^{k,\Phi}_{\ell}+B_{k,\ell}v^{t,x}_\ell+\bar{B}_{k,\ell}\mathbb{E}_kv^{t,x}_\ell+f_{k,\ell}\Big{]}\\[1mm]
\hphantom{X^{k,\Phi}_{\ell+1} =}+\sum_{i=1}^p\Big{[}\big{(}C^i_{k,\ell}+D^i_{k,\ell}\Phi_\ell\big{)} X^{k,\Phi}_{\ell}+\big{(}\bar{C}^i_{k,\ell}+\bar{D}^i_{k,\ell}\Phi_\ell\big{)} \mathbb{E}_kX^{k,\Phi}_{\ell}\\[1mm]
\hphantom{X^{k,\Phi}_{\ell+1} =}+D^i_{k,\ell}v^{t,x}_\ell+\bar{D}^i_{k,\ell}\mathbb{E}_kv^{t,x}_\ell+d^i_{k,\ell}\Big{]} w^i_\ell,\\[1mm]
X^{k,\Phi}_{k} = X^{t,x,*}_k,~~k\in \mathbb{T}_t,~~\ell\in \mathbb{T}_k,
\end{array}
\right.\\
&&\label{X-u}
\hspace{-2em}\left\{\begin{array}{l}
X^{k,u_k,\Phi}_{\ell+1} =\Big{[}\big{(}A_{k,\ell}+B_{k,\ell}\Phi_\ell\big{)} X^{k,u_k,\Phi}_{\ell}+\big{(}\bar{A}_{k,\ell}+\bar{B}_{k,\ell}\Phi_\ell\big{)} \mathbb{E}_kX^{k,u_k,\Phi}_{\ell}\\[1mm]
\hphantom{X^{k,u_k,\Phi}_{\ell+1} =}+B_{k,\ell}v^{t,x}_\ell+\bar{B}_{k,\ell}\mathbb{E}_kv^{t,x}_\ell+f_{k,\ell}\Big{]}\\[1mm]
\hphantom{X^{k,u_k,\Phi}_{\ell+1} =}+\sum_{i=1}^p\Big{[}\big{(}C^i_{k,\ell}+D^i_{k,\ell}\Phi_\ell\big{)} X^{k,u_k,\Phi}_{\ell}+\big{(}\bar{C}^i_{k,\ell}+\bar{D}^i_{k,\ell}\Phi_\ell\big{)} \mathbb{E}_kX^{k,u_k,\Phi}_{\ell}\\[1mm]
\hphantom{X^{k,u_k,\Phi}_{\ell+1} =}+D^i_{k,\ell}v^{t,x}_\ell+\bar{D}^i_{k,\ell}\mathbb{E}_kv^{t,x}_\ell+d^i_{k,\ell}\Big{]} w^i_\ell,\\[1mm]
X^{k,u_k,\Phi}_{k+1} =\big{[}\mathcal{A}_{k,k}X^{k,u_k,\Phi}_k+\mathcal{B}_{k,k}u_k+f_{k,k}\big{]}+\sum_{i=1}^p\big{[}\mathcal{C}^i_{k,k}X^{k,u_k,\Phi}_k+\mathcal{D}^i_{k,k}u_k+d^i_{k,k}\big{]}w^i_k,\\[1mm]
X^{k,u_k,\Phi}_{k} = X^{t,x,*}_k,~~\ell\in \mathbb{T}_{k+1}.
\end{array}\right.
\end{eqnarray}
The state $X^{t,x,*}_k$ in (\ref{defi-linear feedback}), (\ref{Sec2:Defini-closed-X-1}), and (\ref{X-u}) is computed via
\begin{eqnarray}\label{equili-state-closed}
\left\{\begin{array}{l}
X^{t,x,*}_{k+1} =\big{[}\big{(}\mathcal{A}_{k,k}+\mathcal{B}_{k,k}\Phi_k\big{)} X^{t,x,*}_{k}+\mathcal{B}_{k,k}v^{t,x}_k+f_{k,k}\big{]}\\[1mm]
\hphantom{X^{t,x,*}_{k+1} =}+\sum_{i=1}^p\big{[}\big{(}\mathcal{C}^i_{k,k}+\mathcal{D}^i_{k,k}\Phi_k\big{)} X^{t,x,*}_{k}+\mathcal{D}^i_{k,k}v^{t,x}_k+d^i_{k,k}\big{]} w^i_k,\\[1mm]
X^{t,x,*}_{t} = x,~~k\in \mathbb{T}_t.
\end{array}
\right.
\end{eqnarray}

\end{itemize}

ii). $\Phi$ and $v^{t,x}$ in i) are, respectively, the pure-feedback-strategy part and the open-loop-control part of the mixed equilibrium solution $(\Phi, v^{t,x})$.

iii). Letting $\Phi=0$ in i), the corresponding $v^{t,x}$ satisfying (\ref{defi-linear feedback}), denoted as $\widehat{v}^{t,x}$, is then an open-loop equilibrium control of Problem (LQ)$_{tx}$.

iv). If $v^{t,x}$ of i) {happens not to depend} on $x$ and denote such $v^{t,x}$ as $v$, then the corresponding $(\Phi, v)$ is a linear feedback equilibrium strategy of Problem (LQ)$_{tx}$. Here, $ 
l^2(\mathbb{T}_{t}; \mathbb{R}^m)=\big\{\nu=\{\nu_k, k\in \mathbb{T}_{t}\}\, \big{|}\,
|\nu_k|^2 < \infty,  k\in\mathbb{T}_{t}\big\}$.

\end{definition}

\begin{remark}\label{Remark-linear feedback}

By the definition, a mixed equilibrium solution $(\Phi,v^{t,x})$ is time consistent along $X^{t,x,*}$, namely, for any $k\in \mathbb{T}_t$, $(\Phi, v^{t,x})|_{\mathbb{T}_k}$ is a mixed equilibrium solution for the initial pair $(k, X^{t,x,*}_k)$.
Noting $(\Phi\cdot X^{k,\Phi}+v^{t,x})|_{\mathbb{T}_k}=\big{(}\Phi_k X_k^{k,\Phi}+v^{t,x}_k, (\Phi\cdot X^{k,\Phi}+v^{t,x})|_{\mathbb{T}_{k+1}}\big{)}$,
$(u_k,(\Phi\cdot X^{k,u_t,\Phi}+v^{t,x})|_{\mathbb{T}_{k+1}})\big{)}$ is obtained by replacing $\Phi_k X_k^{k,\Phi}+v^{t,x}_k$ and $X^{k,\Phi}$ of $(\Phi\cdot X^{k,\Phi}+v^{t,x})|_{\mathbb{T}_k}$ with $u_k$ and $X^{k,u_k, \Phi}$, respectively.
Furthermore, note that the $v^{t,x}$'s on both sides of (\ref{defi-linear feedback}) are the same. This is why we call $\Phi$ the pure-feedback-strategy part and $v^{t,x}$ the open-loop-control part.

\end{remark}

%
%
%

%

%

\section{Characterization of the mixed equilibrium solution}\label{Section-characterization}

\subsection{The case with the fixed time-state initial pair $(t,x)$}

The following lemma describes the cost difference formula under control perturbation.

\begin{lemma}\label{Lemma-difference}
Let $\bar{u}_k\in l^2_{\mathcal{F}} (k;\mathbb{R}^m)$ and  $\lambda \in
\mathbb{R}$. Then, we have
\begin{eqnarray}\label{appendix-A-J-0}
&&\hspace{-1.5em}J\big{(}k, X_k^{t,x,*}; (\Phi_k \bar{X}^{k,\bar{u}_k,\lambda}_k+v^{t,x}_k +\lambda \bar{u}_k,(\Phi \cdot \bar{X}^{k,\bar{u}_k,\lambda}+v^{t,x})|_{\mathbb{T}_{k+1}})\big{)}-J\big{(}k, X_k^{t,x,*}; (\Phi\cdot X^{k,\Phi}+v^{t,x})|_{\mathbb{T}_k}\big{)}\nonumber\\
&&\hspace{-1.5em}=2\lambda \Big{[}\mathcal{R}_{k,k}(\Phi_kX^{k,\Phi}_k+v^{t,x}_k)+\mathcal{B}^T_{k,k}\mathbb{E}_kY_{k+1}^{k,\Phi} +\sum_{i=1}^p(\mathcal{D}^i_{k,k})^T\mathbb{E}_k({Y}_{k+1}^{k,\Phi}w^i_{k})+\rho_{k,k}
\Big{]}^T\bar{u}_k+ \lambda^2\widetilde{J}(k,0;\bar{u}_k),
\end{eqnarray}
where
\begin{eqnarray}\label{Sec3:lemma-J}
&&\hspace{-1.5em}\widetilde{J}(k,0;\bar{u}_k)=\sum_{\ell=k}^{N-1}\mathbb{E}_k\Big{[}(\alpha^{k,\bar{u}_k}_\ell)^T \big{(}Q_{k,\ell}+\Phi_\ell^TR_{k,\ell}\Phi_\ell \big{)}\alpha^{k,\bar{u}_k}_\ell+(\mathbb{E}_k\alpha^{k,\bar{u}_k}_\ell)^T\big{(}\bar{Q}_{k,\ell}+\Phi_\ell^T\bar{R}_{k,\ell}\Phi_\ell\big{)} \mathbb{E}_k\alpha^{k,\bar{u}_k}_\ell  \Big{]}\nonumber\\
&&\hspace{-1.5em}\hphantom{\widetilde{J}(k,0;\bar{u}_k)=}+\mathbb{E}_k\big{[}\bar{u}_k^T\mathcal{R}_{k,k}\bar{u}_k
\big{]}+\mathbb{E}_k\big{[}(\alpha^{k,\bar{u}_k}_N)^TG_k \alpha_N^{k,\bar{u}_k} \big{]}+(\mathbb{E}_k\alpha^{k,\bar{u}_k}_N)^T\bar{G}_k \mathbb{E}_k\alpha_N^{k,\bar{u}_k},
\end{eqnarray}
and $\bar{X}^{k,\bar{u}_k,\lambda}$, $\alpha^{k,\bar{u}_k}$, ${Y}^{k,\Phi}$
are given, respectively, by the S$\Delta$Es
\begin{eqnarray}
&&\label{Sec3:barX}
\left\{\begin{array}{l}
\bar{X}^{k,\bar{u}_k,\lambda}_{\ell+1} =\big{[}\big{(}A_{k,\ell}+B_{k,\ell}\Phi_\ell\big{)} \bar{X}^{k,\bar{u}_k,\lambda}_{\ell}+\big{(}\bar{A}_{k,\ell}+\bar{B}_{k,\ell}\Phi_\ell\big{)} \mathbb{E}_k\bar{X}^{k,\bar{u}_k,\lambda}_{\ell}\\[1mm]
\hphantom{\bar{X}^{k,\bar{u}_k,\lambda}_{\ell+1} =}+B_{k,\ell} v^{t,x}_\ell+\bar{B}_{k,\ell}\mathbb{E}_k v^{t,x}_\ell+f_{k,\ell}\big{]}\\[1mm]
\hphantom{\bar{X}^{k,\bar{u}_k,\lambda}_{\ell+1} =}+\sum_{i=1}^p\big{[}\big{(}C^i_{k,\ell}+D^i_{k,\ell}\Phi_\ell\big{)} \bar{X}^{k,\bar{u}_k,\lambda}_{\ell}+\big{(}\bar{C}^i_{k,\ell}+\bar{D}^i_{k,\ell}\Phi_\ell\big{)} \mathbb{E}_k\bar{X}^{k,\bar{u}_k,\lambda}_{\ell}\\[1mm]
\hphantom{\bar{X}^{k,\bar{u}_k,\lambda}_{\ell+1} =}+D^i_{k,\ell} v^{t,x}_\ell+\bar{D}^i_{k,\ell}\mathbb{E}_k v^{t,x}_\ell+d^i_{k,\ell}\big{]} w^i_\ell,\\[1mm]
\bar{X}^{k,\bar{u}_k,\lambda}_{k+1} =\big{[}\big{(}\mathcal{A}_{k,k}+\mathcal{B}_{k,k}\Phi_k\big{)}\bar{X}^{k,\bar{u}_k,\lambda}_k+\mathcal{B}_{k,k} v^{t,x}_k+\lambda\mathcal{B}_{k,k}\bar{u}_k+f_{k,k}\big{]}\\[1mm]
\hphantom{\bar{X}^{k,\bar{u}_k,\lambda}_{k+1}= }+\sum_{i=1}^p\big{[}\big{(}\mathcal{C}^i_{k,k}+\mathcal{D}^i_{k,k}\Phi_k\big{)}\bar{X}^{k,\bar{u}_k,\lambda}_k+\mathcal{D}^i_{k,k} v^{t,x}_k+\lambda\mathcal{D}^i_{k,k}\bar{u}_k+d^i_{k,k}\big{]}w^i_k,\\[1mm]
\bar{X}^{k,\bar{u}_k,\lambda}_{k} = X^{t,x,*}_k,~~\ell\in \mathbb{T}_{k+1},
\end{array}\right.\\
&&\label{system-y-k}
\left\{
\begin{array}{l}
{\alpha}^{k,\bar{u}_k}_{\ell+1}=\big{(}A_{k,\ell}+B_{k,\ell}\Phi_\ell\big{)}\alpha^{k,\bar{u}_k}_\ell +\big{(}\bar{A}_{k,\ell}+\bar{B}_{k,\ell}\Phi_\ell\big{)}\mathbb{E}_k\alpha^{k,\bar{u}_k}_\ell\\[1mm]
\hphantom{{Y}^{k,\bar{u}_k}_{\ell+1}=}+\sum_{i=1}^p\big{[}\big{(}C^i_{k,\ell}+D^i_{k,\ell}\Phi_\ell\big{)} \alpha^{k,\bar{u}_k}_\ell+\big{(}\bar{C}^i_{k,\ell} +\bar{D}^i_{k,\ell}\Phi_\ell\big{)}\mathbb{E}_k\alpha^{k,\bar{u}_k}_\ell\big{]}w^i_\ell,\\[1mm]
\alpha^{k,\bar{u}_k}_{k+1}=\mathcal{B}_{k,k}\bar{u}_k+\sum_{i=1}^p{\mathcal{D}}^i_{k,k}\bar{u}_k w^i_k,\\[1mm]
{\alpha}^{k,\bar{u}_k}_k=0,~~~\ell\in \mathbb{T}_{k+1},
\end{array}
\right.
\end{eqnarray}
and the backward stochastic difference equation (BS$\Delta$E, for short)
\begin{eqnarray}\label{Sec3:bsde}
\left\{
\begin{array}{l}
{Y}^{k,\Phi}_\ell=\big{[}Q_{k,\ell}+\Phi_\ell^TR_{k,\ell}\Phi_\ell\big{]} X^{k,\Phi}_\ell+\big{[}\bar{Q}_{k,\ell}+\Phi_\ell^T\bar{R}_{k,\ell}\Phi_\ell\big{]} \mathbb{E}_kX^{k,\Phi}_\ell\\[1mm]
\hphantom{{Z}^{k,\Phi}_\ell=}+(A_{k,\ell}+B_{k,\ell}\Phi_\ell)^T\mathbb{E}_\ell {Y}^{k,\Phi}_{\ell+1}+(\bar{A}_{k,\ell}+\bar{B}_{k,\ell}\Phi_\ell)^T\mathbb{E}_k {Y}_{\ell+1}^{k,\Phi}\\[1mm]
\hphantom{{Z}^{k,\Phi}_\ell=}
+\sum_{i=1}^p\big{[}(C^i_{k,\ell}+D^i_{k,\ell}\Phi_{\ell})^T\mathbb{E}_\ell({Y}_{\ell+1}^{k,\Phi}w^i_\ell)+(\bar{C}^i_{k,\ell} +\bar{D}^i_{k,\ell}\Phi_\ell)\mathbb{E}_{k}({Y}^{k,\Phi}_{\ell+1}w^i_\ell)\big{]}\\[1mm]
\hphantom{{Z}^{k,\Phi}_\ell=}+\Phi_\ell^TR_{k,\ell}v^{t,x}_\ell+\Phi_\ell^T\bar{R}_{k,\ell}\mathbb{E}_kv^{t,x}_\ell+\Phi_\ell^T\rho_{k,\ell}+q_{k,\ell}, \\[1mm]
{Y}^{k,\Phi}_N=G_k{X}_N^{k,\Phi}+\bar{G}_{k}\mathbb{E}_k{X}_N^{k,\Phi}+F_kX^{t,x,*}_k+g_k,~~~~~\ell\in \mathbb{T}_k.
\end{array}
\right.
\end{eqnarray}

\end{lemma}

\emph{Proof.}  See Appendix \ref{appendix-0}. \hfill $\square$


\begin{theorem}\label{Theorem-Equivalentce-linear feedback}
The following statements are equivalent:

\begin{itemize}
\item[i)] Problem (LQ)$_{tx}$ admits a mixed equilibrium solution.

\item[ii)] There exists a pair $(\Phi, v^{t,x})\in l^2(\mathbb{T}_t; \mathbb{R}^{m\times n})\times l^2_\mathcal{F}(\mathbb{T}_t; \mathbb{R}^{m})$ such that the stationary condition
\begin{eqnarray}\label{stationary-condition-closed}
0=\mathcal{R}_{k,k} (\Phi_kX^{k,\Phi}_k+v^{t,x}_k)+\mathcal{B}^T_{k,k}\mathbb{E}_kY_{k+1}^{k,\Phi} +\sum_{i=1}^p(\mathcal{D}^i_{k,k})^T\mathbb{E}_k\big{(}Y_{k+1}^{k,\Phi}w^i_k\big{)} +\rho_{k,k},~~~k\in \mathbb{T}_t
\end{eqnarray}
and the convexity condition
\begin{eqnarray}\label{convex-closed}
&&\inf_{\bar{u}_k\in l^2_\mathcal{F}(k; \mathbb{R}^m)}\widetilde{J}(k,0;\bar{u}_k)\geq 0,~~~k\in \mathbb{T}_t
\end{eqnarray}
hold. Here, $Y^{k,\Phi}_{k+1}$ is computed via the following FBS$\Delta$E
\begin{eqnarray}\label{system-adjoint-closed}
\left\{\begin{array}{l}
X^{k,\Phi}_{\ell+1} =\big{[}\big{(}A_{k,\ell}+B_{k,\ell}\Phi_\ell\big{)} X^{k,\Phi}_{\ell}+\big{(}\bar{A}_{k,\ell}+\bar{B}_{k,\ell}\Phi_\ell\big{)} \mathbb{E}_kX^{k,\Phi}_{\ell}+B_{k,\ell}v^{t,x}_\ell+\bar{B}_{k,\ell}\mathbb{E}_kv^{k,x}_\ell +f_{k,\ell}\big{]}\\[1mm]
\hphantom{X^{k,\Phi}_{\ell+1} =}+\sum_{i=1}^p\big{[}\big{(}C^i_{k,\ell}+D^i_{k,\ell}\Phi_\ell\big{)} X^{k,\Phi}_{\ell}+\big{(}\bar{C}^i_{k,\ell}+\bar{D}^i_{k,\ell}\Phi_\ell\big{)} \mathbb{E}_kX^{k,\Phi}_{\ell}\\[1mm]
\hphantom{X^{k,\Phi}_{\ell+1} =}+D^i_{k,\ell}v^{t,x}_\ell+\bar{D}^i_{k,\ell}\mathbb{E}_kv^{k,x}_\ell+d^i_{k,\ell}\big{]} w^i_\ell,\\[1mm]
Y_{\ell}^{k,\Phi}=Q_{k,\ell}X^{k,\Phi}_\ell+\bar{Q}_{k,\ell}\mathbb{E}_kX^{k,\Phi}_\ell+\Phi_\ell^TR_{k,\ell}\Phi_\ell X^{k,\Phi}_\ell+\Phi_\ell^T\bar{R}_{k,\ell}\Phi_\ell \mathbb{E}_kX^{k,\Phi}_\ell+\Phi_\ell^TR_{k,\ell}v^{t,x}_\ell\\[1mm]
\hphantom{Z_{\ell}^{k,\Phi}=}+\Phi_\ell^T{\bar{R}}_{k,\ell}\mathbb{E}_kv^{t,x}_\ell+\big{(}A_{k,\ell}+B_{k,\ell}\Phi_\ell\big{)}^T \mathbb{E}_\ell Y_{\ell+1}^{k,\Phi}+\big{(}\bar{A}_{k,\ell}+\bar{B}_{k,\ell}\Phi_\ell\big{)}^T\mathbb{E}_kY_{\ell+1}^{k,\Phi} \\[1mm]
\hphantom{Z_{\ell}^{k,\Phi}=} +\sum_{i=1}^p\big{[}\big{(}C^i_{k,\ell}+D^i_{k,\ell}\Phi_\ell\big{)}^T\mathbb{E}_\ell\big{(} Y_{\ell+1}^{k,\Phi}w^i_\ell\big{)}+\big{(}{{\bar{C}}^i_{k,\ell}+{\bar{D}}^i_{k,\ell}\Phi_\ell}\big{)}^T\mathbb{E}_k\big{(}Y^{k,\Phi}_{\ell+1}w^i_\ell\big{)}\big{]} \\[1mm]
\hphantom{Z_{\ell}^{k,\Phi}=} +\Phi_{\ell}^T\rho_{k,\ell} +q_{k,\ell},\\[1mm]
X^{k,\Phi}_k=X^{t,x,*}_k,~~Y_N^{k,\Phi}=G_kX^{k,\Phi}_N+\bar{G}_k\mathbb{E}_kX^{k,\Phi}_N+F_kX^{t,x,*}_k+g_k,~~~~~\ell\in \mathbb{T}_k,
\end{array}
\right.
\end{eqnarray}
and $\widetilde{J}(k,0;\bar{u}_k)$ is given in (\ref{Sec3:lemma-J}). In (\ref{system-adjoint-closed}), $X^{t,x,*}_k$ is computed via
\begin{eqnarray*}
\left\{\begin{array}{l}
X^{t,x,*}_{k+1} =\big{[}\big{(}\mathcal{A}_{k,k}+\mathcal{B}_{k,k}\Phi_k\big{)} X^{t,x,*}_{k}+\mathcal{B}_{k,k}v^{t,x}_k+f_{k,k}\big{]}\\[1mm]
\hphantom{X^{t,x,*}_{k+1} =}+\sum_{i=1}^p\big{[}\big{(}\mathcal{C}^i_{k,k}+\mathcal{D}^i_{k,k}\Phi_k\big{)} X^{t,x,*}_{k}+\mathcal{D}^i_{k,k}v^{t,x}_k+d^i_{k,k}\big{]} w^i_k,\\[1mm]
X^{t,x,*}_{t} = x,~~k\in \mathbb{T}_t.
\end{array}
\right.
\end{eqnarray*}

\end{itemize}

Furthermore, under any of the above conditions, $(\Phi, v^{t,x})$ given in ii) is a mixed equilibrium solution of Problem (LQ)$_{tx}$.

\end{theorem}

\emph{Proof.} This follows from the definition and Lemma \ref{Lemma-difference}. \hfill$\square$

To proceed, we first study the expression of $Y^{k,\Phi}$ of (\ref{system-adjoint-closed}) under some additional condition.

\begin{lemma}\label{Sec3:Lemma-Z}
If for $k\in \mathbb{T}_t$, $v^{t,x}_k=\Gamma_kX^{t,x,*}_k+\bar{v}^{t,x}_k$ with $\Gamma_k, \bar{v}^{t,x}_k$ being deterministic, then the backward state $Y^{k,\Phi}$ of (\ref{system-adjoint-closed}) has the following expression:
\begin{eqnarray*}\label{Sec3:Lemma-Z--0}
\hspace{-0.5em}Y^{k,\Phi}_{\ell}=S_{k,\ell}X^{k,\Phi}_{\ell}+\bar{S}_{k,\ell}\mathbb{E}_{k}X^{k,\Phi}_{\ell}+T_{k,\ell}X^{t,x,*}_\ell
+\bar{T}_{k,\ell}\mathbb{E}_kX^{t,x,*}_\ell+U_{k,\ell}X^{t,x,*}_k+\pi_{k,\ell},~\,\ell\in \mathbb{T}_k,~\,k\in \mathbb{T}_t,
\end{eqnarray*}
where
\begin{eqnarray*}\label{Sec3:Lemma-Z--0-2}
\left\{
\begin{array}{l}
S_{k,\ell}=Q_{k,\ell}+\Phi_{\ell}^TR_{k,\ell}\Phi_{\ell}+\big{(}A_{k,\ell} +B_{k,\ell}\Phi_{\ell}\big{)}^TS_{k,\ell+1}\big{(}A_{k,\ell}+B_{k,\ell}\Phi_{\ell}\big{)}\\[1mm]
\hphantom{S_{k,\ell}=}+\sum_{i,j=1}^p\delta_\ell^{ij}\big{(}C^i_{k,\ell} +D^i_{k,\ell}\Phi_\ell\big{)}^TS_{k,\ell+1}\big{(}C^j_{k,\ell}+D^j_{k,\ell}\Phi_\ell\big{)},\\[1mm]
\bar{S}_{k,\ell}=\bar{Q}_{k,\ell}+\Phi_\ell^T\bar{R}_{k,\ell}\Phi_{\ell}+\big{(}A_{k,\ell} +B_{k,\ell}\Phi_\ell\big{)}^T\big{[}S_{k,\ell+1}(\bar{A}_{k,\ell}+\bar{B}_{k,\ell}\Phi_{\ell}) \\[1mm]
\hphantom{\bar{S}_{k,\ell}=}+\bar{S}_{k,\ell+1}({\mathcal{A}}_{k,\ell}+{\mathcal{B}}_{k,\ell}\Phi_{\ell}) \big{]}+\big{(}\bar{A}_{k,\ell}+\bar{B}_{k,\ell}\Phi_\ell\big{)}^T(S_{k,\ell+1}+\bar{S}_{k,\ell+1})(\mathcal{A}_{k,\ell}+\mathcal{B}_{k,\ell}\Phi_{\ell})\\[1mm]
\hphantom{\bar{S}_{k,\ell}=}+\sum_{i,j=1}^p\delta_\ell^{ij}\big{[}\big{(}C^i_{k,\ell}+D^i_{k,\ell}\Phi_\ell\big{)}^TS_{k,\ell+1}(\bar{C}^j_{k,\ell} +\bar{D}^j_{k,\ell}\Phi_{\ell})\\[1mm]
\hphantom{\bar{S}_{k,\ell}=}+\big{(}{\bar{C}}^i_{k,\ell}+{\bar{D}}^i_{k,\ell}\Phi_{\ell}\big{)}^TS_{k,\ell+1}(\mathcal{C}^j_{k,\ell}+\mathcal{D}^j_{k,\ell}\Phi_{\ell})\big{]},\\[1mm] %
S_{k,N}=G_k,~~\bar{S}_{k,N}=\bar{G}_k,~~\ell\in \mathbb{T}_k,
\end{array}
\right.
\end{eqnarray*}
\begin{eqnarray*}
\left\{
\begin{array}{l}
T_{k,\ell}=\Big{\{}\Phi_\ell^TR_{k,\ell}+\big{(}A_{k,\ell}+B_{k,\ell}\Phi_\ell\big{)}^TS_{k,\ell+1}B_{k,\ell}+\sum_{i,j=1}^p\delta_\ell^{ij}\big{(}C^i_{k,\ell}+D^i_{k,\ell}\Phi_\ell\big{)}^TS_{k,\ell+1}D^j_{k,\ell}\Big{\}}\Gamma_\ell\\[1mm]
\hphantom{T_{k,\ell}=}+\big{(}A_{k,\ell}+B_{k,\ell}\Phi_{\ell}\big{)}^TT_{k,\ell+1}\big{(}\mathcal{A}_{\ell,\ell} +\mathcal{B}_{\ell,\ell}\Phi_\ell+\mathcal{B}_{\ell,\ell}\Gamma_{\ell} \big{)}\\[1mm]
\hphantom{T_{k,\ell}=}+\sum_{i,j=1}^p\delta_\ell^{ij}\big{(}C^i_{k,\ell}+D^i_{k,\ell}\Phi_{\ell}\big{)}^TT_{k,\ell+1}\big{(}\mathcal{C}^j_{\ell,\ell} +\mathcal{D}^j_{\ell,\ell}\Phi_\ell+\mathcal{D}^j_{\ell,\ell}\Gamma_{\ell} \big{)},\\[1mm]
\bar{T}_{k,\ell}=\Big{\{}\Phi_\ell^T{\bar{R}}_{k,\ell}+\big{(}A_{k,\ell}+B_{k,\ell}\Phi_\ell\big{)}^T \big{(}S_{k,\ell+1}\bar{B}_{k,\ell}+\bar{S}_{k,\ell+1}\mathcal{B}_{k,\ell}\big{)}+\big{(}\bar{A}_{k,\ell}+\bar{B}_{k,\ell}\Phi_\ell\big{)}^T\\[1mm]
\hphantom{\bar{T}_{k,\ell}=}\times \big{(}S_{k,\ell+1}+\bar{S}_{k,\ell+1}\big{)}\mathcal{B}_{k,\ell}+\sum_{i,j=1}^p \delta_\ell^{ij}\big{[}\big{(}C^i_{k,\ell}+D^i_{k,\ell}\Phi_\ell\big{)}^TS_{k,\ell+1}\bar{D}^j_{k,\ell}\\[1mm]
\hphantom{\bar{T}_{k,\ell}=}+\big{(}{\bar{C}}^i_{k,\ell}+{\bar{D}}^i_{k,\ell}\Phi_{\ell}\big{)}^T S_k\mathcal{D}^j_{k,\ell} \big{]}\Big{\}}\Gamma_{\ell}+\big{(}A_{k,\ell}+B_{k,\ell}\Phi_\ell\big{)}^T\bar{T}_{k,\ell+1}\big{(}\mathcal{A}_{\ell,\ell}+\mathcal{B}_{\ell,\ell}\Phi_\ell +\mathcal{B}_{\ell,\ell}\Gamma_\ell \big{)}\\[1mm]
\hphantom{\bar{T}_{k,\ell}=}+\big{(}\bar{A}_{k,\ell}+\bar{B}_{k,\ell}\Phi_\ell\big{)}^T\big{(}T_{k,\ell+1}+\bar{T}_{k,\ell+1}\big{)}\big{(}\mathcal{A}_{\ell,\ell}+\mathcal{B}_{\ell,\ell}\Phi_\ell +\mathcal{B}_{\ell,\ell}\Gamma_\ell \big{)}\\[1mm]
\hphantom{\bar{T}_{k,\ell}=}+\sum_{i,j=1}^p\delta_\ell^{ij}\big{(}\bar{C}^i_{k,\ell}+\bar{D}^i_{k,\ell}\Phi_\ell\big{)}^TT_{k,\ell+1}\big{(}\mathcal{C}^j_{\ell,\ell} +\mathcal{D}^j_{\ell,\ell}\Phi_\ell +\mathcal{D}^j_{\ell,\ell}\Gamma_\ell \big{)}\\[1mm]
T_{k,N}=0,~~\bar{T}_{k,N}=0,~~\ell\in \mathbb{T}_k,
\end{array}
\right.
\end{eqnarray*}
\begin{eqnarray*}
\left\{
\begin{array}{l}
U_{k,\ell}=(\mathcal{A}_{k,\ell}+\mathcal{B}_{k,\ell}\Phi_\ell)U_{k,\ell+1},\quad\\
U_{k,N}=F_k,\quad \ell\in \mathbb{T}_{k},
\end{array}
\right.
\end{eqnarray*}
and
\begin{eqnarray*}
\left\{
\begin{array}{l}
\pi_{k,\ell}=\beta_{k,\ell}\bar{v}^{t,x}_{\ell}+\big{(}\mathcal{A}_{k,\ell}+\mathcal{B}_{k,\ell}\Phi_\ell\big{)}^T\big{(}(S_{k,\ell+1}+\bar{S}_{k,\ell+1})f_{k,\ell}+\pi_{k,\ell+1}\big{)} \\[1mm]
\hphantom{\pi_{k,\ell}=}+\big{(}\mathcal{A}_{k,\ell}+\mathcal{B}_{k,\ell}\Phi_\ell\big{)}^T\big{(}T_{k,\ell+1}+\bar{T}_{k,\ell+1}\big{)} f_{\ell,\ell}+\sum_{i,j=1}^p\delta_\ell^{ij}\big{[}\big{(}\mathcal{C}^i_{k,\ell}+\mathcal{D}^i_{k,\ell}\Phi_\ell\big{)}^TS_{k,\ell+1}d^j_{k,\ell}\\[1mm]
\hphantom{\pi_{k,\ell}=}+\big{(}\mathcal{C}^i_{k,\ell}+\mathcal{D}^i_{k,\ell}\Phi_\ell\big{)}^TT_{k,\ell+1}d^j_{\ell,\ell}\big{]}+\Phi_\ell^T\rho_{k,\ell} +q_{k,\ell},\quad\\
\pi_{k,N}=g_k,~~\ell\in \mathbb{T}_k
\end{array}
\right.
\end{eqnarray*}
with
\begin{eqnarray*}
\begin{array}{l}
\beta_{k,\ell}=\Phi_\ell^T\mathcal{R}_{k,\ell}+\big{(}\mathcal{A}_{k,\ell}+\mathcal{B}_{k,\ell}\Phi_\ell\big{)}^T\big{[} (S_{k,\ell+1}+\bar{S}_{k,\ell+1})\mathcal{B}_{k,\ell}+(T_{k,\ell+1}+\bar{T}_{k,\ell+1})\mathcal{B}_{\ell,\ell} \big{]}\\[1mm]
\hphantom{\beta_{k,\ell}=}+\sum_{i,j=1}^p\big{(}\mathcal{C}^i_{k,\ell}+\mathcal{D}^i_{k,\ell}\Phi_\ell\big{)}^T \big{[}S_{k,\ell+1}\mathcal{D}^j_{k,\ell}+T_{k,\ell+1}\mathcal{D}^j_{\ell,\ell}\big{]},\quad~~
\ell\in \mathbb{T}_k.
\end{array}
\end{eqnarray*}

\end{lemma}

\emph{Proof}. See Appendix \ref{appendix-Sec2:lemma-Z}. \hfill $\square$
%

For a matrix $M\in \mathbb{R}^{n\times m}$, let $M^\dagger$ be its Moore-Penrose inverse. Then, we have
the following lemma \cite{Ait-Chen-Zhou-2002}.

\begin{lemma}\label{Lemma-matrix-equation}
Let matrices $L$, $M$, and $N$ be given with appropriate size. Then,
$LXM=N$ has a solution $X$ if and only if $LL^\dagger NMM^\dagger=N$.
Moreover, the solution of $LXM=N$ can be expressed as
$X=L^\dagger NM^\dagger+V-L^\dagger LVMM^\dagger$,
where $V$ is a matrix with appropriate size.
\end{lemma}

If $M=I$ in Lemma \ref{Lemma-matrix-equation}, then $LL^\dagger N=N$ is equivalent to $\mbox{Ran}(N)\subset \mbox{Ran}(L)$. Here, $\mbox{Ran}(N)$ is the range of $N$. %
The following theorem is concerned with the necessary and sufficient condition for the existence of a mixed equilibrium solution.

\begin{theorem}\label{Sec3:Theorem-Necessary}
The following statements are equivalent:

\begin{itemize}

\item[i)] Problem (LQ)$_{tx}$ admits a mixed equilibrium solution.

\item[ii)] There exists $\Phi\in l^2(\mathbb{T}_t;\mathbb{R}^{m\times n})$ such that the following assertions hold.
\begin{itemize}

\item[a)] The coupled equations
\begin{eqnarray}\label{Sec3:Theorem-Necessary-S}
\left\{
\begin{array}{l}
\left\{
\begin{array}{l}
S_{k,\ell}=Q_{k,\ell}+\Phi_\ell^TR_{k,\ell}\Phi_{\ell}+\big{(}A_{k,\ell} +B_{k,\ell}\Phi_\ell\big{)}^TS_{k,\ell+1}\big{(}A_{k,\ell}+B_{k,\ell}\Phi_\ell\big{)}\\[1mm]
\hphantom{S_{k,\ell}=}+\sum_{i,j=1}^p\delta_\ell^{ij}\big{(}C^i_{k,\ell} +D^i_{k,\ell}\Phi_\ell\big{)}^TS_{k,\ell+1}\big{(}C^j_{k,\ell}+D^j_{k,\ell}\Phi_\ell\big{)},\\[1mm]
{\mathcal{S}}_{k,\ell}
={\mathcal{Q}}_{k,\ell}+\Phi_\ell^T{\mathcal{R}}_{k,\ell}\Phi_{\ell}+\big{(}\mathcal{A}_{k,\ell} +\mathcal{B}_{k,\ell}\Phi_\ell\big{)}^T\mathcal{S}_{k,\ell+1}({\mathcal{A}}_{k,\ell}+{\mathcal{B}}_{k,\ell}\Phi_{\ell}) \\[1mm]
\hphantom{\bar{S}_{k,\ell}=}+\sum_{i,j=1}^p\delta_\ell^{ij}\big{(}\mathcal{C}^i_{k,\ell}+\mathcal{D}^i_{k,\ell}\Phi_\ell\big{)}^T S_{k,\ell+1}({\mathcal{C}}^j_{k,\ell}+{\mathcal{D}}^j_{k,\ell}\Phi_{\ell}),\\[1mm]
%
%
%
S_{k,N}=G_k,~~{\mathcal{S}}_{k,N}=G_k+\bar{G}_k,~~\ell\in \mathbb{T}_k, \quad

\end{array}
\right.\\[1mm]
{\mathbb{O}}_k=\mathcal{R}_{k,k}+\mathcal{B}^T_{k,k}\mathcal{S}_{k,k+1}\mathcal{B}_{k,k}+\sum_{i,j=1}^p\delta_k^{ij}(\mathcal{D}^i_{k,k})^T{S}_{k,k+1} \mathcal{D}^j_{k,k}\succeq 0,\quad k\in \mathbb{T}_t
\end{array}
\right.
\end{eqnarray}
are solvable in the sense of
$
{\mathbb{O}}_k\succeq 0, k\in \mathbb{T}_t
$, namely, $\mathbb{O}_k, k\in \mathbb{T}_t$, are all nonnegative definite.

\item[b)] The condition
\begin{eqnarray}\label{stationary-condition-closed-1}
\mathcal{L}_{k}X^{t,x,*}_{k}+\theta_{k}\in \mbox{Ran}\big{(}\mathcal{O}_{k}\big{)},~~~k\in \mathbb{T}_t
\end{eqnarray}
is satisfied. Here, $X^{t,x,*}$ is computed via
\begin{eqnarray*}\label{equili-state-closed-2}
\left\{\begin{array}{l}
X^{t,x,*}_{k+1} =\big{[}\big{(}\mathcal{A}_{k,k}-\mathcal{B}_{k,k}\mathcal{O}_k^{\dagger}\mathcal{L}_k \big{)}X^{t,x,*}_{k}-\mathcal{B}_{k,k}\mathcal{O}^\dagger_k\theta_k +f_{k,k}\big{]}\\[1mm]
\hphantom{X^{t,x,*}_{k+1} =}+\sum_{i=1}^p\big{[}\big{(}\mathcal{C}^i_{k,k}-\mathcal{D}^i_{k,k}\mathcal{O}_k^{\dagger}\mathcal{L}_k \big{)} X^{t,x,*}_{k}-\mathcal{D}^i_{k,k}\mathcal{O}^\dagger_k\theta_k+d^i_{k,k}\big{]} w^i_k,\\[1mm]
X^{t,x,*}_{t} = x,~~k\in \mathbb{T}_t,
\end{array}
\right.
\end{eqnarray*}
and $\mathcal{O}_k, \mathcal{L}_k, \theta_{k}, k\in \mathbb{T}_t$ are given by
\begin{eqnarray*}
\left\{
\begin{array}{l}
\mathcal{O}_{k}=\mathcal{R}_{k,k}+\mathcal{B}^T_{k,k}\big{(}\mathcal{S}_{k,k+1}+\mathcal{T}_{k,k+1}\big{)}\mathcal{B}_{k,k} +\sum_{i,j=1}^p\delta_k^{ij}(\mathcal{D}^i_{k,k})^T\big{(}{S}_{k,k+1}+T_{k,k+1}\big{)}\mathcal{D}^j_{k,k},\\[1mm]
\mathcal{L}_{k}=\mathcal{B}^T_{k,k}\big{(}\mathcal{S}_{k,k+1}+\mathcal{T}_{k,k+1}\big{)}\mathcal{A}_{k,k}+\sum_{i,j=1}^p\delta_k^{ij}(\mathcal{D}^i_{k,k})^T \big{(}S _{k,k+1}+T_{k,k+1}\big{)}\mathcal{C}^j_{k,k}\\
\hphantom{\mathcal{L}_{k}=}+\mathcal{B}_{k,k}^TU_{k,k+1},\\[1mm]
\theta_{k}=\mathcal{B}^T_{k,k}\big{(}\mathcal{S}_{k,k+1}+\mathcal{T}_{k,k+1}\big{)}f_{k,k}+\sum_{i,j=1}^p\delta_k^{ij}(\mathcal{D}^i_{k,k})^T\big{(}S_{k,+1}+T_{k,k+1}\big{)}d^j_{k,k}\\[1mm]
\hphantom{\theta_{k}=}+\mathcal{B}^T_{k,k}\pi_{k,k+1}+\rho_{k,k},\quad\\[1mm]
k\in \mathbb{T}_t,
\end{array}
\right.
\end{eqnarray*}
where
\begin{eqnarray*}
\left\{
\begin{array}{l}
T_{k,\ell}=\Big{\{}\Phi_\ell^TR_{k,\ell}+\big{(}A_{k,\ell}+B_{k,\ell}\Phi_\ell\big{)}^TS_{k,\ell+1}B_{k,\ell} \\[1mm]
\hphantom{T_{k,\ell}=} +\sum_{i,j=1}^p\delta_\ell^{ij}\big{(}C^i_{k,\ell}+D^i_{k,\ell}\Phi_\ell\big{)}^TS_{k,\ell+1}D^j_{k,\ell}\Big{\}}\Gamma_\ell\\[1mm]
\hphantom{T_{k,\ell}=}+\big{(}A_{k,\ell}+B_{k,\ell}\Phi_{\ell}\big{)}^TT_{k,\ell+1}\big{(}\mathcal{A}_{\ell,\ell} +\mathcal{B}_{\ell,\ell}\Phi_\ell+\mathcal{B}_{\ell,\ell}\Gamma_{\ell} \big{)}\\[1mm]
\hphantom{T_{k,\ell}=}+\sum_{i,j=1}^p\delta_\ell^{ij}\big{(}C^i_{k,\ell}+D^i_{k,\ell}\Phi_{\ell}\big{)}^TT_{k,\ell+1}\big{(}\mathcal{C}^j_{\ell,\ell} +\mathcal{D}^j_{\ell,\ell}\Phi_\ell+\mathcal{D}^j_{\ell,\ell}\Gamma_{\ell} \big{)},\\[1mm]
{\mathcal{T}}_{k,\ell}=\Big{\{}\Phi_\ell^T{{\mathcal{R}}}_{k,\ell}+\big{(}\mathcal{A}_{k,\ell}+\mathcal{B}_{k,\ell}\Phi_\ell\big{)}^T \mathcal{S}_{k,\ell+1}{\mathcal{B}}_{k,\ell}\\[1mm]
\hphantom{\bar{T}_{k,\ell}=}+\sum_{i,j=1}^p\delta_\ell^{ij} \big{(}\mathcal{C}^i_{k,\ell}+\mathcal{D}^i_{k,\ell}\Phi_\ell\big{)}^TS_{k,\ell+1}{\mathcal{D}}^j_{k,\ell}\Big{\}}\Gamma_{\ell}\\[1mm]
\hphantom{\bar{T}_{k,\ell}=}+\big{(}\mathcal{A}_{k,\ell}+\mathcal{B}_{k,\ell}\Phi_\ell\big{)}^T{\mathcal{T}}_{k,\ell+1}\big{(}\mathcal{A}_{\ell,\ell}+\mathcal{B}_{\ell,\ell}\Phi_\ell +\mathcal{B}_{\ell,\ell}\Gamma_\ell \big{)}\\[1mm]
%
%
\hphantom{\bar{T}_{k,\ell}=}+\sum_{i,j=1}^p\delta_\ell^{ij}\big{(}{\mathcal{C}}^i_{k,\ell}+{\mathcal{D}}^i_{k,\ell}\Phi_\ell\big{)}^TT_{k,\ell+1} \big{(}\mathcal{C}^j_{\ell,\ell} +\mathcal{D}^j_{\ell,\ell}\Phi_\ell +\mathcal{D}^j_{\ell,\ell}\Gamma_\ell \big{)}\\[1mm]
T_{k,N}=0,~~{\mathcal{T}}_{k,N}=0,~~\ell\in \mathbb{T}_k, \quad k\in \mathbb{T}_t,
\end{array}
\right.
\end{eqnarray*}
\begin{eqnarray*}
\left\{
\begin{array}{l}
U_{k,\ell}=(\mathcal{A}_{k,\ell}+\mathcal{B}_{k,\ell}\Phi_\ell)U_{k,\ell+1},\quad\\
U_{k,N}=F_k,~~\ell\in \mathbb{T}_{k}, \quad
k\in \mathbb{T}_t,
\end{array}
\right.
\end{eqnarray*}
and
\begin{eqnarray*}
\left\{
\begin{array}{l}
\pi_{k,\ell}=-\beta_{k,\ell}\mathcal{O}^\dagger_\ell\theta_\ell+\big{(}\mathcal{A}_{k,\ell}+\mathcal{B}_{k,\ell}\Phi_\ell\big{)}^T\big{(}\mathcal{S}_{k,\ell+1}f_{k,\ell}+\pi_{k,\ell+1}\big{)} \\[1mm]
\hphantom{\pi_{k,\ell}=}+\sum_{i,j=1}^p\delta_\ell^{ij}\big{[}\big{(}\mathcal{C}^i_{k,\ell}+\mathcal{D}^i_{k,\ell}\Phi_\ell\big{)}^TS_{k,\ell+1}d^j_{k,\ell} +\big{(}\mathcal{C}^i_{k,\ell}+\mathcal{D}^i_{k,\ell}\Phi_\ell\big{)}^TT_{k,\ell+1}d^j_{\ell,\ell}\\[1mm]
\hphantom{\pi_{k,\ell}=}+\big{(}\mathcal{A}_{k,\ell}+\mathcal{B}_{k,\ell}\Phi_\ell\big{)}^T\mathcal{T}_{k,\ell+1} f_{\ell,\ell}+\Phi_\ell^T\rho_{k,\ell} +q_{k,\ell},\\[1mm]
\pi_{k,N}=g_k,~~\ell\in \mathbb{T}_k, \quad k\in \mathbb{T}_t
\end{array}
\right.
\end{eqnarray*}
with
\begin{eqnarray}\label{Phi+Gamma}
\Gamma_k=-\mathcal{O}_k^{\dagger}\mathcal{L}_k-\Phi_k,~~~k\in \mathbb{T}_t
\end{eqnarray}
and
\begin{eqnarray*}
%
\begin{array}{l}
\beta_{k,\ell}=\Phi_\ell^T\mathcal{R}_{k,\ell}+\big{(}\mathcal{A}_{k,\ell}+\mathcal{B}_{k,\ell}\Phi_\ell\big{)}^T\big{[} \mathcal{S}_{k,\ell+1}\mathcal{B}_{k,\ell}+\mathcal{T}_{k,\ell+1}\mathcal{B}_{\ell,\ell} \big{]}\\[1mm]
\hphantom{\beta_{k,\ell}=}+\sum_{i,j=1}^p\delta_\ell^{ij}\big{(}\mathcal{C}^i_{k,\ell}+\mathcal{D}^i_{k,\ell}\Phi_\ell\big{)}^T \big{[}S_{k,\ell+1}\mathcal{D}^j_{k,\ell}+T_{k,\ell+1}\mathcal{D}^j_{\ell,\ell}\big{]},\quad
\ell\in \mathbb{T}_k.
\end{array}
\end{eqnarray*}

\end{itemize}
\end{itemize}

Furthermore, under condition ii), let
\begin{eqnarray}\label{v}
v^{t,x}_k=\Gamma_kX^{t,x,*}_k-\mathcal{O}^\dagger_k\theta_k, ~~k\in \mathbb{T}_t,
\end{eqnarray}
and $\Phi_k, \Gamma_k, k\in \mathbb{T}_t$, are given in ii); then, $(\Phi, v^{t,x})$ is a mixed equilibrium solution of Problem (LQ)$_{tx}$.

\end{theorem}

\emph{Proof.} See Appendix \ref{appendix:theorem-necessary}.
\hfill $\square$

\begin{remark}\label{Remark-Phi}
In Theorem \ref{Sec3:Theorem-Necessary}, the solvability of (\ref{Sec3:Theorem-Necessary-S}) is to characterize the convexity (\ref{convex-closed}), while (\ref{stationary-condition-closed-1}) is to characterize the stationary condition (\ref{stationary-condition-closed}).
If $\Phi_{\ell}$, $\Gamma_\ell$, $\ell\in \mathbb{T}_{k+1}$ have been determined, then $-\mathcal{O}_k^\dagger \mathcal{L}_k$ can be further constructed. Noting (\ref{Phi+Gamma}), it is impossible to determine the value of $\Phi_k$ by using the property $\Phi_k+\Gamma_k=-\mathcal{O}_k^\dagger \mathcal{L}_k$, and any $(\Phi, v^{t,x})$ that satisfies condition ii) of Theorem \ref{Sec3:Theorem-Necessary} is a mixed equilibrium solution. Nevertheless, the freedom of selecting $\Phi$ could enable us to deal with the open-loop equilibrium control and linear feedback equilibrium strategy in a unified way.

\end{remark}

From Theorem \ref{Sec3:Theorem-Necessary}, the following two corollaries are straightforward. The first concerns the open-loop equilibrium control, which is obtained by letting $\Phi=0$ in Theorem \ref{Sec3:Theorem-Necessary}.

\begin{corollary}\label{Corollary--open}
The following statements are equivalent:

\begin{itemize}

\item[i)] Problem (LQ)$_{tx}$ admits an open-loop equilibrium control.

\item[ii)] The following assertions hold.
\begin{itemize}

\item[a)]  The coupled equations
\begin{eqnarray}\label{Sec3:Theorem-Necessary-S-open}
\left\{
\begin{array}{l}
\left\{
\begin{array}{l}
\widehat{S}_{k,\ell}=Q_{k,\ell}+A_{k,\ell}^T\widehat{S}_{k,\ell+1}A_{k,\ell}+\sum_{i,j=1}^p\delta_{\ell}^{ij}(C^i_{k,\ell})^T\widehat{S}_{k,\ell+1}C^j_{k,\ell},\\[1mm]
\widehat{\mathcal{S}}_{k,\ell}
={\mathcal{Q}}_{k,\ell}+\mathcal{A}_{k,\ell}^T\widehat{\mathcal{S}}_{k,\ell+1}{\mathcal{A}}_{k,\ell}+\sum_{i,j=1}^p\delta_{\ell}^{ij}(\mathcal{C}^i_{k,\ell})^T\widehat{S}_{k,\ell+1}{\mathcal{C}}^j_{k,\ell},\\[1mm]
%
%
%
\widehat{S}_{k,N}=G_k,~~\widehat{\mathcal{S}}_{k,N}=G_k+\bar{G}_k,
\quad \ell\in \mathbb{T}_k,
\end{array}
\right.\\
\widehat{\mathbb{O}}_k\succeq 0, \quad k\in \mathbb{T}_t
\end{array}
\right.
\end{eqnarray}
are solvable in the sense of $\widehat{\mathbb{O}}_k\succeq 0, k\in \mathbb{T}_t$.

\item[b)] The condition
\begin{eqnarray}\label{station-all}
\widehat{\mathcal{L}}_{k}\widehat{X}^{t,x,*}_{k}+\widehat{\theta}_{k}\in \mbox{Ran}\big{(}\widehat{\mathcal{O}}_{k}\big{)},~~~k\in \mathbb{T}_t
\end{eqnarray}
is satisfied. Here, $\widehat{X}^{t,x,*}$ is computed via
\begin{eqnarray*}
\left\{\begin{array}{l}
\widehat{X}^{t,x,*}_{k+1} =\big{[}\big{(}\mathcal{A}_{k,k}-\mathcal{B}_{k,k}\widehat{\mathcal{O}}_k^{\dagger}\widehat{\mathcal{L}}_k \big{)}\widehat{X}^{t,x,*}_{k}-\mathcal{B}_{k,k}\widehat{\mathcal{O}}^\dagger_k\widehat{\theta}_k +f_{k,k}\big{]}\\[1mm]
\hphantom{\widehat{X}^{t,x,*}_{k+1} =}+\sum_{i=1}^p\big{[}\big{(}\mathcal{C}^i_{k,k}-\mathcal{D}^i_{k,k}\widehat{\mathcal{O}}_k^{\dagger}\widehat{\mathcal{L}}_k\big{)} \widehat{X}^{t,x,*}_{k}-\mathcal{D}^i_{k,k}\widehat{\mathcal{O}}^\dagger_k\widehat{\theta}_k+d^i_{k,k}\big{]} w^i_k,\\[1mm]
\widehat{X}^{t,x,*}_{t} = x,~~k\in \mathbb{T}_t,
\end{array}
\right.
\end{eqnarray*}
and  $\widehat{\mathcal{O}}_k, \widehat{\mathcal{L}}_k, \theta_k, k\in \mathbb{T}_t$ are given by
\begin{eqnarray}\label{O-L-pi-open-fixed}
\left\{
\begin{array}{l}
\widehat{\mathcal{O}}_{k}=\mathcal{R}_{k,k}+\mathcal{B}^T_{k,k}\big{(}\widehat{\mathcal{S}}_{k,k+1}+\widehat{\mathcal{T}}_{k,k+1}\big{)}\mathcal{B}_{k,k} +\sum_{i,j=1}^p\delta_{k}^{ij}(\mathcal{D}^i_{k,k})^T\big{(}\widehat{S}_{k,k+1}+\widehat{T}_{k,k+1}\big{)}\mathcal{D}^j_{k,k},\\[1mm]
\widehat{\mathcal{L}}_{k}=\mathcal{B}^T_{k,k}\big{(}\widehat{\mathcal{S}}_{k,k+1}+\widehat{\mathcal{T}}_{k,k+1}\big{)}\mathcal{A}_{k,k} +\sum_{i,j=1}^p\delta_{k}^{ij}(\mathcal{D}^i_{k,k})^T \big{(}\widehat{S} _{k,k+1}+\widehat{T}_{k,k+1}\big{)}\mathcal{C}^j_{k,k}\\[1mm]
\hphantom{\widehat{\mathcal{L}}_{k}=}+\mathcal{B}_{k,k}^T\widehat{U}_{k,k+1},\\[1mm]
\widehat{\theta}_{k}=\mathcal{B}^T_{k,k}\big{(}\widehat{\mathcal{S}}_{k,k+1}+\widehat{\mathcal{T}}_{k,k+1}\big{)}f_{k,k} +\sum_{i,j=1}^p\delta_k^{ij}(\mathcal{D}^i_{k,k})^T\big{(}\widehat{S}_{k,+1} +\widehat{T}_{k,k+1}\big{)}d^j_{k,k}\\[1mm]
\hphantom{\widehat{\theta}_{k}=}+\mathcal{B}^T_{k,k}\widehat{\pi}_{k,k+1}+\rho_{k,k},\quad\\[1mm]
k\in \mathbb{T}_t,
\end{array}
\right.
\end{eqnarray}
where
\begin{eqnarray*}
\left\{
\begin{array}{l}
\widehat{T}_{k,\ell}=A_{k,\ell}^T\widehat{T}_{k,\ell+1}\mathcal{A}_{\ell,\ell}+\sum_{i,j=1}^p\delta_\ell^{ij}(C^i_{k,\ell})^T\widehat{T}_{k,\ell+1}\mathcal{C}^j_{\ell,\ell} -\Big{\{}A_{k,\ell}^T\widehat{S}_{k,\ell+1}B_{k,\ell}\\[1mm]
\hphantom{\widehat{T}_{k,\ell}=}+\sum_{i,j=1}^p(C^i_{k,\ell})^T\widehat{S}_{k,\ell+1}D^j_{k,\ell} +A_{k,\ell}^T\widehat{T}_{k,\ell+1} \mathcal{B}_{\ell,\ell}+\sum_{i,j=1}^p\delta_\ell^{ij}(C^i_{k,\ell})^T\widehat{T}_{k,\ell+1}\mathcal{D}^j_{\ell,\ell}\Big{\}}\widehat{\mathcal{O}}_\ell^{\dagger}\widehat{\mathcal{L}}_\ell,\\[1mm]
\widehat{\mathcal{T}}_{k,\ell}=\mathcal{A}_{k,\ell}^T{\widehat{\mathcal{T}}}_{k,\ell+1}\mathcal{A}_{\ell,\ell}+ \sum_{i,j=1}^p\delta_\ell^{ij}({\mathcal{C}}^i_{k,\ell})^T\widehat{T}_{k,\ell+1}\mathcal{C}^j_{\ell,\ell}-\Big{\{}\mathcal{A}_{k,\ell}^T \widehat{\mathcal{S}}_{k,\ell+1}{\mathcal{B}}_{k,\ell}\\[1mm]
\hphantom{\widehat{\mathcal{T}}_{k,\ell}=}+\sum_{i,j=1}^p(\mathcal{C}^i_{k,\ell})^T\widehat{S}_{k,\ell+1}{\mathcal{D}}^j_{k,\ell} +\mathcal{A}_{k,\ell}^T{\widehat{\mathcal{T}}}_{k,\ell+1}\mathcal{B}_{\ell,\ell}+ \sum_{i,j=1}^p\delta_\ell^{ij}({\mathcal{C}}^i_{k,\ell})^T\widehat{T}_{k,\ell+1}\mathcal{D}^j_{\ell,\ell}\Big{\}}\widehat{\mathcal{O}}_\ell^{\dagger}\widehat{\mathcal{L}}_\ell,\\[1mm]
\widehat{T}_{k,N}=0,~~\widehat{\mathcal{T}}_{k,N}=0,~~\ell\in \mathbb{T}_k, \quad k\in \mathbb{T}_t,
\end{array}
\right.
\end{eqnarray*}
\begin{eqnarray*}
\left\{
\begin{array}{l}
\widehat{U}_{k,\ell}=\mathcal{A}_{k,\ell}\widehat{U}_{k,\ell+1},\quad\\
\widehat{U}_{k,N}=F_k,\quad \ell\in \mathbb{T}_{k}, \quad
k\in \mathbb{T}_t,
\end{array}
\right.
\end{eqnarray*}
and
\begin{eqnarray*}
\left\{
\begin{array}{l}
\widehat{\pi}_{k,\ell}=-\widehat{\beta}_{k,\ell}\widehat{\mathcal{O}}^\dagger_\ell\widehat{\theta}_\ell+\mathcal{A}_{k,\ell}^T\big{(}\widehat{\mathcal{S}}_{k,\ell+1}f_{k,\ell}+\widehat{\pi}_{k,\ell+1}\big{)} +\mathcal{A}_{k,\ell}^T\widehat{\mathcal{T}}_{k,\ell+1} f_{\ell,\ell}\\[1mm]
\hphantom{\widehat{\pi}_{k,\ell}=}+\sum_{i,j=1}^p\delta_\ell^{ij}\big{[}(\mathcal{C}^i_{k,\ell})^T\widehat{S}_{k,\ell+1}d^j_{k,\ell} +(\mathcal{C}^i_{k,\ell})^T\widehat{T}_{k,\ell+1}d^j_{\ell,\ell}\big{]}+q_{k,\ell},\\[1mm]
\widehat{\pi}_{k,N}=g_k,\quad \ell\in \mathbb{T}_k, \quad k\in \mathbb{T}_t,
\end{array}
\right.
\end{eqnarray*}
with
$$
\widehat{\beta}_{k,\ell}=\mathcal{A}_{k,\ell}^T\big{[} \widehat{\mathcal{S}}_{k,\ell+1}\mathcal{B}_{k,\ell}+\widehat{\mathcal{T}}_{k,\ell+1}\mathcal{B}_{\ell,\ell} \big{]}+\sum_{i,j=1}^p\delta_\ell^{ij}(\mathcal{C}^i_{k,\ell})^T \big{[}\widehat{S}_{k,\ell+1}\mathcal{D}^j_{k,\ell}+\widehat{T}_{k,\ell+1}\mathcal{D}^j_{\ell,\ell}\big{]}, ~~
\ell\in \mathbb{T}_k.
$$

\end{itemize}

\end{itemize}

Furthermore, under condition ii), the control
$$\widehat{v}^{t,x}_k=-\widehat{\mathcal{O}}^\dagger_k\widehat{\mathcal{L}}_k\widehat{X}^{t,x,*}_k
-\widehat{\mathcal{O}}^\dagger_k\widehat{\theta}_k,~~~
k\in \mathbb{T}_t$$
is an open-loop equilibrium control of Problem (LQ)$_{tx}$.

\end{corollary}

Note that the linear feedback equilibrium strategy has nothing to do with the initial state $x$. The second corollary is concerned with the existence of a linear feedback equilibrium strategy, which is obtained by letting $\Gamma_k=0$, $k\in \mathbb{T}_t$ in Theorem \ref{Sec3:Theorem-Necessary}.

\begin{corollary}\label{Corollary--feedback}
The following statements are equivalent:

\begin{itemize}

\item[i)] Problem (LQ)$_{tx}$ admits a linear feedback equilibrium strategy.

\item[ii)] The following assertions hold.
\begin{itemize}

\item[a)] The coupled equations
\begin{eqnarray}\label{Sec3:Theorem-Necessary-S-feedback}
\left\{
\begin{array}{l}
\left\{
\begin{array}{l}
\widetilde{S}_{k,\ell}=Q_{k,\ell}+\widetilde{\Phi}_{\ell}^TR_{k,\ell}\widetilde{\Phi}_{\ell}+\big{(}A_{k,\ell} +B_{k,\ell}\widetilde{\Phi}_{\ell}\big{)}^T\widetilde{S}_{k,\ell+1}\big{(}A_{k,\ell}+B_{k,\ell}\widetilde{\Phi}_{\ell}\big{)}\\[1mm]
\hphantom{\widetilde{S}_{k,\ell}=}+\sum_{i,j=1}^p\delta_\ell^{ij}\big{(}C^i_{k,\ell} +D^i_{k,\ell}\widetilde{\Phi}_{\ell}\big{)}^T\widetilde{S}_{k,\ell+1}\big{(}C^j_{k,\ell}+D^j_{k,\ell}\widetilde{\Phi}_{\ell}\big{)},\\[1mm]
\widetilde{\mathcal{S}}_{k,\ell}
={\mathcal{Q}}_{k,\ell}+\widetilde{\Phi}_{\ell}^T{\mathcal{R}}_{k,\ell}\widetilde{\Phi}_{\ell}+\big{(}\mathcal{A}_{k,\ell} +\mathcal{B}_{k,\ell}\widetilde{\Phi}_{\ell}\big{)}^T\widetilde{\mathcal{S}}_{k,\ell+1}({\mathcal{A}}_{k,\ell}+{\mathcal{B}}_{k,\ell}\widetilde{\Phi}_{\ell}) \\[1mm]
\hphantom{\widetilde{\bar{S}}_{k,\ell}=}+\sum_{i,j=1}^p\delta_\ell^{ij}\big{(}\mathcal{C}^i_{k,\ell}+\mathcal{D}^i_{k,\ell}\widetilde{\Phi}_{\ell}\big{)}^T\widetilde{S}_{k,\ell+1} ({\mathcal{C}}^j_{k,\ell}+{\mathcal{D}}^j_{k,\ell}\widetilde{\Phi}_{\ell}),\\[1mm]
%
%
%
\widetilde{S}_{k,N}=G_k,~~\widetilde{\mathcal{S}}_{k,N}=G_k+\bar{G}_k, \quad \ell\in \mathbb{T}_k,
\end{array}
\right.\\
\widetilde{\mathbb{O}}_k=\mathcal{R}_{k,k}+\mathcal{B}^T_{k,k}\widetilde{\mathcal{S}}_{k,k+1}\mathcal{B}_{k,k} +\sum_{i,j=1}^p\delta_k^{ij}(\mathcal{D}^i_{k,k})^T\widetilde{S}_{k,k+1}\mathcal{D}^j_{k,k}\succeq 0, \quad k\in \mathbb{T}_t
\end{array}
\right.
\end{eqnarray}
are solvable in the sense of $
\widetilde{\mathbb{O}}_k\succeq 0, k\in \mathbb{T}_t.
$

\item[b)] The condition
\begin{eqnarray}\label{station-strategy}
\widetilde{{\mathbb{L}}}_{k}\widetilde{X}^{t,x,*}_{k}+\widetilde{\theta}_{k}\in \mbox{Ran}\big{(}\widetilde{\mathbb{O}}_k\big{)},~~~k\in \mathbb{T}_t
\end{eqnarray}
is satisfied. Here, $\widetilde{X}^{t,x,*}$ is computed via
\begin{eqnarray*}\label{equili-state-closed-2}
\left\{\begin{array}{l}
\widetilde{X}^{t,x,*}_{k+1} =\big{[}\big{(}\mathcal{A}_{k,k}-\mathcal{B}_{k,k}{\widetilde{\mathbb{O}}}_k^{\dagger}{\widetilde{\mathbb{L}}}_k \big{)}\widetilde{X}^{t,x,*}_{k}-\mathcal{B}_{k,k}{\widetilde{\mathbb{O}}}^\dagger_k\widetilde{\theta}_k +f_{k,k}\big{]}\\[1mm]
\hphantom{X^{t,x,*}_{k+1} =}+\sum_{i=1}^p\big{[}\big{(}\mathcal{C}^i_{k,k}-\mathcal{D}^i_{k,k}{\widetilde{\mathbb{O}}}_k^{\dagger}{\widetilde{\mathbb{L}}}_k \big{)} \widetilde{X}^{t,x,*}_{k}-\mathcal{D}^i_{k,k}{\widetilde{\mathbb{O}}}^\dagger_k\widetilde{\theta}_k+d^i_{k,k}\big{]} w^i_k,\\[1mm]
\widetilde{X}^{t,x,*}_{t} = x,~~k\in \mathbb{T}_t,
\end{array}
\right.
\end{eqnarray*}
and $\widetilde{{\mathbb{L}}}_k, \widetilde{\theta}_k, k\in \mathbb{T}_t$ are given by
\begin{eqnarray*}
\left\{
\begin{array}{l}
%
%
\widetilde{{\mathbb{L}}}_{k}=\mathcal{B}^T_{k,k}\widetilde{\mathcal{S}}_{k,k+1}\mathcal{A}_{k,k}+\sum_{i,j=1}^p\delta_k^{ij}(\mathcal{D}^i_{k,k})^T \widetilde{S}_{k,k+1}\mathcal{C}^j_{k,k}+\mathcal{B}_{k,k}^T\widetilde{U}_{k,k+1},\\[1mm]
\widetilde{\theta}_{k}=\mathcal{B}^T_{k,k}\widetilde{\mathcal{S}}_{k,k+1}f_{k,k}+\sum_{i,j=1}^p\delta_k^{ij}(\mathcal{D}^i_{k,k})^T \widetilde{S}_{k,k+1}d^j_{k,k}+\mathcal{B}^T_{k,k}\widetilde{\pi}_{k,k+1}+\rho_{k,k},\quad
k\in \mathbb{T}_t,
\end{array}
\right.
\end{eqnarray*}
where
\begin{eqnarray*}
\left\{
\begin{array}{l}
\widetilde{U}_{k,\ell}=(\mathcal{A}_{k,\ell}+\mathcal{B}_{k,\ell}\widetilde{\Phi}_{\ell})\widetilde{U}_{k,\ell+1},\\
\widetilde{U}_{k,N}=F_k,~~ \ell\in \mathbb{T}_{k}, ~~k\in \mathbb{T}_t,
\end{array}
\right.
\end{eqnarray*}
\begin{eqnarray*}
\left\{
\begin{array}{l}
\widetilde{\pi}_{k,\ell}=-\widetilde{\beta}_{k,\ell}\widetilde{\mathbb{O}}^\dagger_\ell\widetilde{\theta}_\ell+\big{(}\mathcal{A}_{k,\ell} +\mathcal{B}_{k,\ell}\widetilde{\Phi}_{\ell}\big{)}^T \big{(}\widetilde{\mathcal{S}}_{k,\ell+1}f_{k,\ell}+\widetilde{\pi}_{k,\ell+1}\big{)} \\[1mm]
\hphantom{\widetilde{\pi}_{k,\ell}=}+\sum_{i,j=1}^p\delta_\ell^{ij}\big{(}\mathcal{C}^i_{k,\ell}+\mathcal{D}^i_{k,\ell}\widetilde{\Phi}_{\ell}\big{)}^T \widetilde{S}_{k,\ell+1}d^j_{k,\ell} +\widetilde{\Phi}_{\ell}^T\rho_{k,\ell} +q_{k,\ell},\\[1mm]
\pi_{k,N}=g_k,~~\ell\in \mathbb{T}_k,\quad
k\in \mathbb{T}_t
\end{array}
\right.
\end{eqnarray*}
with
$$
\widetilde{\beta}_{k,\ell}=\widetilde{\Phi}_{\ell}^T\mathcal{R}_{k,\ell}+\big{(}\mathcal{A}_{k,\ell}+\mathcal{B}_{k,\ell}\widetilde{\Phi}_{\ell}\big{)}^T \widetilde{\mathcal{S}}_{k,\ell+1}\mathcal{B}_{k,\ell}
+\sum_{i,j=1}^p\delta_\ell^{ij}\big{(}\mathcal{C}^i_{k,\ell}+\mathcal{D}^i_{k,\ell}\widetilde{\Phi}_{\ell}\big{)}^T \widetilde{S}_{k,\ell+1}\mathcal{D}^j_{k,\ell},~~
\ell\in \mathbb{T}_k.$$

\item[c)] $\widetilde{{\Phi}}$ above is given by
$\{\widetilde{{\Phi}}_k=-\widetilde{\mathbb{O}}_k^\dagger \widetilde{\mathbb{L}}_k, k\in \mathbb{T}_t\}$.

\end{itemize}

\end{itemize}

Furthermore, under condition ii), $\{(-\widetilde{\mathbb{O}}_k^\dagger \widetilde{\mathbb{L}}_k, -\widetilde{\mathbb{O}}^\dagger_k\widetilde{\theta}_k), k\in \mathbb{T}_t\}$ is a linear feedback equilibrium strategy of Problem (LQ)$_{tx}$.

\end{corollary}

We now consider the unique existence of open-loop equilibrium control and linear feedback equilibrium strategy.

\begin{theorem}\label{Theorem-unique}
Let $t\in \mathbb{T}$ and $x\in l^2_{\mathcal{F}}(t;\mathbb{R}^n)$. The following is true.

\begin{itemize}
\item[i)] The following statements are equivalent:

\begin{itemize}

\item[a)] Problem (LQ)$_{tx}$ admits a unique open-loop equilibrium control.

\item[b)] (\ref{Sec3:Theorem-Necessary-S-open}) is solvable, and $\widehat{\mathcal{O}}_k, k\in \mathbb{T}_t$, are invertible, which are given in (\ref{O-L-pi-open-fixed}).

\item[c)]For any $k\in \mathbb{T}_t$ and any $\xi\in l^2_{\mathcal{F}}(k;\mathbb{R}^n)$, Problem (LQ)$_{k\xi}$ admits a unique open-loop equilibrium control.

\end{itemize}

\item[ii)] The following statements are equivalent:

\begin{itemize}

\item[d)] Problem (LQ)$_{tx}$ admits a unique linear feedback equilibrium strategy.

\item[e)] $\widetilde{\mathbb{O}}_k\succ 0, k\in \mathbb{T}_t$, namely, $\widetilde{\mathbb{O}}_k, k\in \mathbb{T}_t$, are all positive definite, which are given in (\ref{Sec3:Theorem-Necessary-S-feedback}).

\item[f)]For any $k\in \mathbb{T}_t$ and any $\xi\in l^2_{\mathcal{F}}(k;\mathbb{R}^n)$, Problem (LQ)$_{k\xi}$ admits a unique linear feedback equilibrium strategy.

\end{itemize}

\end{itemize}

\end{theorem}

\emph{Proof.} i). a)$\Leftrightarrow$b). Let $\widehat{v}^{t,x}$ be an open-loop equilibrium control of Problem (LQ)$_{tx}$. In this case, (\ref{stationary-condition-closed}) becomes
\begin{eqnarray}\label{stationary-condition-open}
0=\mathcal{R}_{k,k}\widehat{v}^{t,x}_k+\mathcal{B}^T_{k,k}\mathbb{E}_kY_{k+1}^{k}+\sum_{i=1}^p(\mathcal{D}^i_{k,k})^T\mathbb{E}_k\big{(}Y_{k+1}^{k}w^i_k\big{)} +\rho_{k,k},~~~k\in \mathbb{T}_t.
\end{eqnarray}
Mimicking the proof of Theorem \ref{Sec3:Theorem-Necessary} and based on Lemma \ref{Lemma-matrix-equation}, a control of the following form
\begin{eqnarray}\label{v-open}
\widehat{v}^{t,x}_k=-\widehat{\mathcal{O}}^\dagger_k\widehat{\mathcal{L}}_k\widehat{X}^{t,x,*}_k-\widehat{\mathcal{O}}^\dagger_k\widehat{\theta}_k+(I-\widehat{\mathcal{O}}^\dagger_k\widehat{\mathcal{O}}_k)V_k, ~~k\in \mathbb{T}_t
\end{eqnarray}
also satisfies (\ref{stationary-condition-open}), where $V_k\in \mathbb{R}^m$ is deterministic and $\widehat{X}^{t,x,*}$, $\widehat{\theta}_k$ are given by %
\begin{eqnarray*}
\left\{\begin{array}{l}
\widehat{X}^{t,x,*}_{k+1} =\big{[}\big{(}\mathcal{A}_{k,k}-\mathcal{B}_{k,k}\widehat{\mathcal{O}}_k^{\dagger}\widehat{\mathcal{L}}_k \big{)}\widehat{X}^{t,x,*}_{k}-\mathcal{B}_{k,k}\big{(}\widehat{\mathcal{O}}^\dagger_k\widehat{\theta}_k-(I-\widehat{\mathcal{O}}^\dagger_k\widehat{\mathcal{O}}_k)V_k\big{)} +f_{k,k}\big{]}\\[1mm]
\hphantom{\widehat{X}^{t,x,*}_{k+1} =}+\sum_{i=1}^p\big{[}\big{(}\mathcal{C}^i_{k,k}-\mathcal{D}^i_{k,k}\widehat{\mathcal{O}}_k^{\dagger}\widehat{\mathcal{L}}_k\big{)} \widehat{X}^{t,x,*}_{k}-\mathcal{D}^i_{k,k}\big{(}\widehat{\mathcal{O}}^\dagger_k\widehat{\theta}_k-(I-\widehat{\mathcal{O}}^\dagger_k\widehat{\mathcal{O}}_k)V_k\big{)} +d^i_{k,k}\big{]} w^i_k,\\[1mm]
\widehat{X}^{t,x,*}_{t} = x,~~k\in \mathbb{T}_t,
\end{array}
\right.
\end{eqnarray*}
and
\begin{eqnarray*}
\widehat{\theta}_{k}=\mathcal{B}^T_{k,k}\big{(}\widehat{\mathcal{S}}_{k,k+1}+\widehat{\mathcal{T}}_{k,k+1}\big{)}f_{k,k} +\sum_{i,j=1}^p\delta_k^{ij}(\mathcal{D}^i_{k,k})^T\big{(}\widehat{S}_{k,+1} +\widehat{T}_{k,k+1}\big{)}d^j_{k,k}+\mathcal{B}^T_{k,k}\widehat{\pi}_{k,k+1}+\rho_{k,k},
\end{eqnarray*}
with $\widehat{\pi}_{k,k+1}$ computed via
\begin{eqnarray*}
\left\{
\begin{array}{l}
\widehat{\pi}_{k,\ell}=-\widehat{\beta}_{k,\ell}\big{[}\widehat{\mathcal{O}}^\dagger_\ell\widehat{\theta}_\ell -(I-\widehat{\mathcal{O}}^\dagger_\ell\widehat{\mathcal{O}}_\ell)V_\ell\big{]}+\mathcal{A}_{k,\ell}^T\big{(}\widehat{\mathcal{S}}_{k,\ell+1}f_{k,\ell}+\widehat{\pi}_{k,\ell+1}\big{)} +\mathcal{A}_{k,\ell}^T\widehat{\mathcal{T}}_{k,\ell+1} f_{\ell,\ell}\\[1mm]
\hphantom{\widehat{\pi}_{k,\ell}=}+\sum_{i,j=1}^p\delta_\ell^{ij}\big{[}(\mathcal{C}^i_{k,\ell})^T\widehat{S}_{k,\ell+1}d^j_{k,\ell} +(\mathcal{C}^i_{k,\ell})^T\widehat{T}_{k,\ell+1}d^j_{\ell,\ell}\big{]}+q_{k,\ell},\\[1mm]
\widehat{\pi}_{k,N}=g_k,~~\ell\in \mathbb{T}_k,\quad k\in \mathbb{T}_t.
\end{array}
\right.
\end{eqnarray*}
Combining the solvability of (\ref{Sec3:Theorem-Necessary-S-open}), we know that any control of the form (\ref{v-open}) is an open-loop equilibrium control of Problem (LQ)$_{tx}$. Therefore, Problem (LQ)$_{tx}$ admits a unique open-loop equilibrium control if and only if $\widehat{\mathcal{O}}_k, k\in \mathbb{T}_t$ are invertible, and (\ref{Sec3:Theorem-Necessary-S-open}) is solvable.

b)$\Rightarrow$c). As (\ref{Sec3:Theorem-Necessary-S-open}) is solvable and $\widehat{\mathcal{O}}_k, k\in \mathbb{T}_t$, are invertible, we know from Corollary \ref{Corollary--open} that Problem (LQ)$_{k\xi}$ admits an open-loop equilibrium control. Based on the proof of a)$\Rightarrow$b), it follows that Problem (LQ)$_{k\xi}$ admits a unique open-loop equilibrium control.

c)$\Rightarrow$a). It is obvious.

ii). d)$\Leftrightarrow$e). Let $(\Phi, v)$ be a linear feedback equilibrium strategy of Problem (LQ)$_{tx}$.  $\widetilde{\mathbb{O}}_k\succeq 0, k\in \mathbb{T}_k$ follows from the solvability of (\ref{Sec3:Theorem-Necessary-S-feedback}). We now prove that $\widetilde{\mathbb{O}}_k, k\in \mathbb{T}_k$ are all invertible. Note that the linear feedback equilibrium strategy is independent of $x$. If some of $\widetilde{\mathbb{O}}_k, k\in \mathbb{T}_k$, are singular, then similar to those of a)$\Rightarrow$b), $v$
can be selected as any one of the following forms:
\begin{eqnarray}\label{v-feedback}
v_k=-\widetilde{\mathcal{O}}^\dagger_k\widetilde{\theta}_k+(I-\widetilde{\mathcal{O}}^\dagger_k\widetilde{\mathcal{O}}_k)V_k, ~~k\in \mathbb{T}_t.
\end{eqnarray}
In (\ref{v-feedback}), $V_k\in \mathbb{R}^m, k\in \mathbb{T}_t$, are deterministic, and $\widetilde{\theta}_k$ is given by
\begin{eqnarray*}
\widetilde{\theta}_{k}=\mathcal{B}^T_{k,k}\widetilde{\mathcal{S}}_{k,k+1}f_{k,k} +\sum_{i,j=1}^p\delta_k^{ij}(\mathcal{D}^i_{k,k})^T \widetilde{S}_{k,+1}d^j_{k,k}+\mathcal{B}^T_{k,k}\widetilde{\pi}_{k,k+1}+\rho_{k,k},
\end{eqnarray*}
with $\widetilde{\pi}_{k,k+1}$ computed via
\begin{eqnarray*}
\left\{
\begin{array}{l}
\widetilde{\pi}_{k,\ell}=-\widetilde{\beta}_{k,\ell}\big{[}\widetilde{\mathcal{O}}^\dagger_\ell\widetilde{\theta}_\ell -(I-\widetilde{\mathcal{O}}^\dagger_\ell \widetilde{\mathcal{O}}_\ell)V_\ell\big{]}+\big{(}\mathcal{A}_{k,\ell} +\mathcal{B}_{k,\ell}\widetilde{\Phi}_{\ell}\big{)}^T \big{(}\widetilde{\mathcal{S}}_{k,\ell+1}f_{k,\ell}+\widetilde{\pi}_{k,\ell+1}\big{)} \\[1mm]
\hphantom{\widetilde{\pi}_{k,\ell}=}+\sum_{i,j=1}^p\delta_\ell^{ij}\big{(}\mathcal{C}^i_{k,\ell}+\mathcal{D}^i_{k,\ell}\widetilde{\Phi}_{\ell}\big{)}^T \widetilde{S}_{k,\ell+1}d^j_{k,\ell} +\widetilde{\Phi}_{\ell}^T\rho_{k,\ell} +q_{k,\ell},\\[1mm]
\pi_{k,N}=g_k,~~\ell\in \mathbb{T}_k,\quad k\in \mathbb{T}_t.
\end{array}
\right.
\end{eqnarray*}
Therefore, Problem (LQ)$_{tx}$ admits a unique linear feedback equilibrium strategy if and only if $\widetilde{\mathbb{O}}_k, k\in \mathbb{T}_k$ are all invertible, and thus are positive definite.

e)$\Rightarrow$f) and f)$\Rightarrow$d) are obvious. This completes the proof. \hfill $\square$

\begin{remark}
Problem (LQ)$_{tx}$ admitting a unique open-loop equilibrium control is a local property,  which is only of the unique existence  for the fixed initial pair $(t,x)$. Interestingly, this local property could ensure a semi-global property, namely, for any $k\in \mathbb{T}_t$ (after $t$) and any $\xi\in l_\mathcal{F}^2(k; \mathbb{R}^n)$, Problem (LQ)$_{k\xi}$ also admits a unique open-loop equilibrium control. A similar property also holds for the linear feedback equilibrium strategy.

\end{remark}

\subsection{The case with all of the initial pairs}

Simply knowing that Problem (LQ)$_{tx}$ admits an open-loop equilibrium control or a linear feedback equilibrium strategy, it is hard or generally impossible to derive sharp results like those of Theorem \ref{Theorem-unique}. Alternatively, in this section, we consider the case that the initial pair is allowed to vary. To begin, we first state the results for the open-loop equilibrium control.

\begin{theorem}\label{Theorem-open-all}
The following statements are equivalent:

\begin{itemize}

\item[i)] For any $(t,x)$ with $t\in \mathbb{T}$ and $x\in l_{\mathcal{F}}^2(t;\mathbb{R}^n)$, Problem (LQ)$_{tx}$ admits an open-loop equilibrium control.

\item[ii)] The coupled equations
\begin{eqnarray}\label{S-all}
\left\{
\begin{array}{l}
\left\{
\begin{array}{l}
\widehat{S}_{k,\ell}=Q_{k,\ell}+A_{k,\ell}^T\widehat{S}_{k,\ell+1}A_{k,\ell}+\sum_{i,j=1}^p\delta_\ell^{ij}(C^i_{k,\ell})^T\widehat{S}_{k,\ell+1}C^j_{k,\ell},\\[1mm]
\widehat{\mathcal{S}}_{k,\ell}
={\mathcal{Q}}_{k,\ell}+\mathcal{A}_{k,\ell}^T\widehat{\mathcal{S}}_{k,\ell+1}{\mathcal{A}}_{k,\ell}+\sum_{i,j=1}^p\delta_\ell^{ij}(\mathcal{C}^i_{k,\ell})^T\widehat{S}_{k,\ell+1}{\mathcal{C}}^j_{k,\ell},\\[1mm]
\widehat{S}_{k,N}=G_k,~~\widehat{\mathcal{S}}_{k,N}=G_k+\bar{G}_k,\quad \ell\in \mathbb{T}_k,
\end{array}
\right.\\[1mm]
\widehat{\mathbb{O}}_k\succeq  0, \quad k\in \mathbb{T},
\end{array}
\right.
\end{eqnarray}
\begin{eqnarray}\label{T-all}
\left\{
\begin{array}{l}
\left\{
\begin{array}{l}
\widehat{T}_{k,\ell}=A_{k,\ell}^T\widehat{T}_{k,\ell+1} \mathcal{A}_{\ell,\ell}+\sum_{i,j=1}^p\delta_\ell^{ij}(C^i_{k,\ell})^T\widehat{T}_{k,\ell+1}\mathcal{C}^j_{\ell,\ell}- \Big{\{}A_{k,\ell}^T\widehat{S}_{k,\ell+1}B_{k,\ell} \\[1mm]
\hphantom{\widehat{T}_{k,\ell}=}+A_{k,\ell}^T\widehat{T}_{k,\ell+1}\mathcal{B}_{\ell,\ell}+\sum_{i,j=1}^p\delta_\ell^{ij}\big{[}(C^i_{k,\ell})^T\widehat{S}_{k,\ell+1}D^j_{k,\ell} +(C^i_{k,\ell})^T\widehat{T}_{k,\ell+1} \mathcal{D}^j_{\ell,\ell}\big{]} \Big{\}}\widehat{\mathcal{O}}_\ell^{\dagger}\widehat{\mathcal{L}}_\ell,\\[1mm]
\widehat{\mathcal{T}}_{k,\ell}=\mathcal{A}_{k,\ell}^T{\widehat{\mathcal{T}}}_{k,\ell+1} \mathcal{A}_{\ell,\ell}+\sum_{i,j=1}^p\delta_\ell^{ij}({\mathcal{C}}^i_{k,\ell})^T\widehat{T}_{k,\ell+1}\mathcal{C}^j_{\ell,\ell}-\Big{\{}\mathcal{A}_{k,\ell}^T \widehat{\mathcal{S}}_{k,\ell+1}{\mathcal{B}}_{k,\ell}\\[1mm]
\hphantom{\widehat{\mathcal{T}}_{k,\ell}=}+\mathcal{A}_{k,\ell}^T\widehat{\mathcal{T}}_{k,\ell+1}\mathcal{B}_{\ell,\ell} + \sum_{i,j=1}^p\delta_\ell^{ij}\big{[}(\mathcal{C}^i_{k,\ell})^T\widehat{S}_{k,\ell+1}{\mathcal{D}}^j_{k,\ell} +(\mathcal{C}^i_{k,\ell})^T\widehat{T}_{k,\ell+1}\mathcal{D}^j_{\ell,\ell}\big{]}\Big{\}}\widehat{\mathcal{O}}_\ell^{\dagger}\widehat{\mathcal{L}}_\ell,\\[1mm]
\widehat{T}_{k,N}=0,~~\widehat{\mathcal{T}}_{k,N}=0, \quad \ell\in \mathbb{T}_k,
\end{array}
\right.\\
\widehat{\mathcal{O}}_k\widehat{\mathcal{O}}_k^\dagger\widehat{\mathcal{L}}_k=\widehat{\mathcal{L}}_k, \quad k\in \mathbb{T},
\end{array}
\right.
\end{eqnarray}
and
\begin{eqnarray}\label{pi-all}
\left\{
\begin{array}{l}
\left\{
\begin{array}{l}
\widehat{\pi}_{k,\ell}=-\widehat{\beta}_{k,\ell}\widehat{\mathcal{O}}^\dagger_\ell\widehat{\theta}_\ell+\mathcal{A}_{k,\ell}^T\big{(}\widehat{\mathcal{S}}_{k,\ell+1}f_{k,\ell}+\widehat{\pi}_{k,\ell+1}\big{)} +\mathcal{A}_{k,\ell}^T\widehat{\mathcal{T}}_{k,\ell+1} f_{\ell,\ell}\\[1mm]
\hphantom{\widehat{\pi}_{k,\ell}=}+\sum_{i,j=1}^p\delta_\ell^{ij}\big{[}(\mathcal{C}^i_{k,\ell})^T\widehat{S}_{k,\ell+1}d^j_{k,\ell} +(\mathcal{C}^i_{k,\ell})^T\widehat{T}_{k,\ell+1}d^j_{\ell,\ell}\big{]}+q_{k,\ell},\\[1mm]
\widehat{\pi}_{k,N}=g_k,\quad
\ell\in \mathbb{T}_k,
\end{array}
\right.\\
\widehat{\mathcal{O}}_k\widehat{\mathcal{O}}_k^\dagger\widehat{\theta}_k=\widehat{\theta}_k, \quad k\in \mathbb{T}
\end{array}
\right.
\end{eqnarray}
are solvable in the sense of
\begin{eqnarray*}
\widehat{\mathbb{O}}_k\succeq 0,~~~
\widehat{\mathcal{O}}_k\widehat{\mathcal{O}}_k^\dagger\widehat{\mathcal{L}}_k-\widehat{\mathcal{L}}_k=0,~~~
\widehat{\mathcal{O}}_k\widehat{\mathcal{O}}_k^\dagger\widehat{\theta}_k-\widehat{\theta}_k=0,~~~
k\in \mathbb{T},
\end{eqnarray*}
where
\begin{eqnarray*}
\left\{
\begin{array}{l}
\widehat{\mathbb{O}}_k=\mathcal{R}_{k,k}+\mathcal{B}^T_{k,k}\widehat{\mathcal{S}}_{k,k+1}\mathcal{B}_{k,k}+\sum_{i,j=1}^p\delta_k^{ij}(\mathcal{D}^i_{k,k})^T\widehat{S}_{k,k+1}\mathcal{D}^j_{k,k},\\[1mm]
\widehat{\mathcal{O}}_{k}=\mathcal{R}_{k,k}+\mathcal{B}^T_{k,k}\big{(}\widehat{\mathcal{S}}_{k,k+1}+\widehat{\mathcal{T}}_{k,k+1}\big{)}\mathcal{B}_{k,k} +\sum_{i,j=1}^p\delta_k^{ij}(\mathcal{D}^i_{k,k})^T\big{(}\widehat{S}_{k,k+1}+\widehat{T}_{k,k+1}\big{)}\mathcal{D}^j_{k,k},\\[1mm]
\widehat{\mathcal{L}}_{k}=\mathcal{B}^T_{k,k}\big{(}\widehat{\mathcal{S}}_{k,k+1}+\widehat{\mathcal{T}}_{k,k+1}\big{)}\mathcal{A}_{k,k} +\sum_{i,j=1}^p\delta_k^{ij}(\mathcal{D}^i_{k,k})^T \big{(}\widehat{S} _{k,k+1}+\widehat{T}_{k,k+1}\big{)}\mathcal{C}^j_{k,k}+\mathcal{B}_{k,k}^T\widehat{U}_{k,k+1},\\[1mm]
\widehat{\theta}_{k}=\mathcal{B}^T_{k,k}\big{(}\widehat{\mathcal{S}}_{k,k+1}+\widehat{\mathcal{T}}_{k,k+1}\big{)}f_{k,k} +\sum_{i,j=1}^p\delta_k^{ij}(\mathcal{D}^i_{k,k})^T\big{(}\widehat{S}_{k,+1} +\widehat{T}_{k,k+1}\big{)}d^j_{k,k}+\mathcal{B}^T_{k,k}\widehat{\pi}_{k,k+1}+\rho_{k,k},\\[1mm]
k\in \mathbb{T}
\end{array}
\right.
\end{eqnarray*}
with
\begin{eqnarray*}
\left\{
\begin{array}{l}
\widehat{U}_{k,\ell}=\mathcal{A}_{k,\ell}\widehat{U}_{k,\ell+1},\\
\widehat{U}_{k,N}=F_k,~~\ell\in \mathbb{T}_{k}, ~~k\in \mathbb{T},
\end{array}
\right.
\end{eqnarray*}
and
$$
\widehat{\beta}_{k,\ell}=\mathcal{A}_{k,\ell}^T\big{[} \widehat{\mathcal{S}}_{k,\ell+1}\mathcal{B}_{k,\ell}+\widehat{\mathcal{T}}_{k,\ell+1}\mathcal{B}_{\ell,\ell} \big{]}+\sum_{i,j=1}^p\delta_\ell^{ij}(\mathcal{C}^i_{k,\ell})^T \big{[}\widehat{S}_{k,\ell+1}\mathcal{D}^j_{k,\ell}+\widehat{T}_{k,\ell+1}\mathcal{D}^j_{\ell,\ell}\big{]},
~~\ell\in \mathbb{T}_k.$$

\end{itemize}

\end{theorem}

\emph{Proof}. \emph{ii)$\Rightarrow$i)}. From the solvability of (\ref{T-all}) and (\ref{pi-all}), we know that (\ref{station-all}) holds for any $(t,x)$.
Therefore, i) holds.

\emph{i)$\Rightarrow$ii)}. Note that (\ref{station-all}) is equivalent to
$\widehat{\mathcal{O}}_k \widehat{\mathcal{O}}_k^\dagger\big{(}\widehat{\mathcal{L}}_k\widehat{X}^{t,x,*}_k+\widehat{\theta}_k \big{)}=\widehat{\mathcal{L}}_k\widehat{X}^{t,x,*}_k+\widehat{\theta}_k$, $k\in \mathbb{T}_t$.
Letting $k=t$ and taking different $x's$, we have
$\widehat{\mathcal{O}}_t \widehat{\mathcal{O}}_t^\dagger\widehat{\mathcal{L}}_t=\widehat{\mathcal{L}}_t$,
$\widehat{\mathcal{O}}_t \widehat{\mathcal{O}}_t^\dagger\widehat{\theta}_t=\widehat{\theta}_t$.
As for any $(t,x)$ with $t\in \mathbb{T}$ and $x\in l^2_\mathcal{F}(t;\mathbb{R}^n)$ Problem (LQ)$_{tx}$ admits an open-loop equilibrium control, we must have the solvability of
(\ref{S-all})-(\ref{pi-all}). \hfill$\square$

The following result is for the feedback equilibrium strategy.

\begin{theorem}\label{Sec3:Theorem-N-S}
The following statements are equivalent:

\begin{itemize}
\item[i)] For any $(t,x)$ with $t\in \mathbb{T}$ and $x\in l^2_{\mathcal{F}}(t; \mathbb{R}^n)$, Problem (LQ)$_{tx}$ admits a linear feedback equilibrium strategy.

\item[ii)] The coupled equations
\begin{eqnarray}\label{S-all-feedback}
\left\{
\begin{array}{l}
\left\{
\begin{array}{l}
\widetilde{S}_{k,\ell}=Q_{k,\ell}+A_{k,\ell}^T\widetilde{S}_{k,\ell+1}A_{k,\ell}+\sum_{i,j=1}^p\delta_\ell^{ij}(C^i_{k,\ell})^T\widetilde{S}_{k,\ell+1}C^j_{k,\ell}\\[1mm]
\hphantom{\widetilde{S}_{k,\ell}=}-\big{(}A_{k,\ell}^T\widetilde{S}_{k,\ell+1}B_{k,\ell}+\sum_{i,j=1}^p\delta_\ell^{ij}(C^i_{k,\ell})^T\widetilde{S}_{k,\ell+1}D^j_{k,\ell} \big{)}\widetilde{\mathbb{O}}^{\dagger}_{\ell} \widetilde{\mathbb{L}}_{\ell}\\[1mm]
\hphantom{\widetilde{S}_{k,\ell}=}-\widetilde{\mathbb{L}}_{\ell}^T\widetilde{\mathbb{O}}^{\dagger}_{\ell}\big{(}B_{k,\ell}^T\widetilde{S}_{k,\ell+1}A_{k,\ell} +\sum_{i,j=1}^p\delta_\ell^{ij}(D^i_{k,\ell})^T\widetilde{S}_{k,\ell+1}C^j_{k,\ell} \big{)}\\[1mm]
\hphantom{\widetilde{S}_{k,\ell}=}+\widetilde{\mathbb{L}}_{\ell}^T\widetilde{\mathbb{O}}^{\dagger}_{\ell}\big{(}R_{k,\ell}+B_{k,\ell}^T\widetilde{S}_{k,\ell+1}B_{k,\ell} +\sum_{i,j=1}^p\delta_\ell^{ij}(D^i_{k,\ell})^T\widetilde{S}_{k,\ell+1}D^j_{k,\ell} \big{)}\widetilde{\mathbb{O}}^{\dagger}_{\ell}\widetilde{\mathbb{L}}_\ell,
\\[1mm]
\widetilde{\mathcal{S}}_{k,\ell}=
%
\mathcal{Q}_{k,\ell}+\mathcal{A}_{k,\ell}^T\widetilde{\mathcal{S}}_{k,\ell+1}\mathcal{A}_{k,\ell} +\sum_{i,j=1}^p\delta_\ell^{ij}(\mathcal{C}^i_{k,\ell})^T\widetilde{S}_{k,\ell+1}\mathcal{C}^j_{k,\ell}\\[1mm]
\hphantom{\widetilde{\mathcal{S}}_{k,\ell}=}-\big{(}\mathcal{A}_{k,\ell}^T\widetilde{\mathcal{S}}_{k,\ell+1}\mathcal{B}_{k,\ell} +\sum_{i,j=1}^p\delta_\ell^{ij}(\mathcal{C}^i_{k,\ell})^T\widetilde{S}_{k,\ell+1}\mathcal{D}^j_{k,\ell} \big{)}\widetilde{\mathbb{O}}^{\dagger}_{\ell} \widetilde{\mathbb{L}}_{\ell}\\[1mm]
\hphantom{\widetilde{\mathcal{S}}_{k,\ell}=}-\widetilde{\mathbb{L}}_{\ell}^T\widetilde{\mathbb{O}}^{\dagger}_{\ell}\big{(}\mathcal{B}_{k,\ell}^T\widetilde{\mathcal{S}}_{k,\ell+1}\mathcal{A}_{k,\ell} +\sum_{i,j=1}^p\delta_\ell^{ij}(\mathcal{D}^i_{k,\ell})^T\widetilde{S}_{k,\ell+1}\mathcal{C}^j_{k,\ell} \big{)}\\[1mm]
\hphantom{\mathcal{S}_{k,\ell}=}+\widetilde{\mathbb{L}}_{\ell}^T\widetilde{\mathbb{O}}^{\dagger}_{\ell}\big{(}\mathcal{R}_{k,\ell} +\mathcal{B}_{k,\ell}^T\widetilde{\mathcal{S}}_{k,\ell+1}\mathcal{B}_{k,\ell}+\sum_{i,j=1}^p\delta_\ell^{ij}(\mathcal{D}^i_{k,\ell})^T\widetilde{S}_{k,\ell+1}\mathcal{D}^j_{k,\ell} \big{)}\widetilde{\mathbb{O}}^{\dagger}_{\ell}\widetilde{\mathbb{L}}_\ell,\\[1mm]
%
%
%
%
\widetilde{S}_{k,N}=G_k,~~\widetilde{\mathcal{S}}_{k,N}=G_k+\bar{G}_k, \quad
\ell\in \mathbb{T}_k,
\end{array}
\right.\\
\widetilde{\mathbb{O}}_k\succeq  0, \quad\widetilde{\mathbb{O}}_k\widetilde{\mathbb{O}}_k^\dagger\widetilde{\mathbb{L}}_k=\widetilde{\mathbb{L}}_k,\quad k\in \mathbb{T},
\end{array}
\right.
\end{eqnarray}
and
\begin{eqnarray}\label{pi-all-feedback}
\left\{
\begin{array}{l}
\left\{
\begin{array}{l}
\widetilde{\pi}_{k,\ell}=-\widetilde{\beta}_{k,\ell}\widetilde{{\mathbb{O}}}^\dagger_k\widetilde{\theta}_k+\big{(}\mathcal{A}_{k,\ell} -\mathcal{B}_{k,\ell}\widetilde{\mathbb{O}}^{\dagger}_{\ell}\widetilde{\mathbb{L}}_\ell\big{)}^T \big{(}\widetilde{\mathcal{S}}_{k,\ell+1}f_{k,\ell}+\widetilde{\pi}_{k,\ell+1}\big{)} \\[1mm]
\hphantom{\widetilde{\pi}_{k,\ell}=}+\sum_{i,j=1}^p\delta_\ell^{ij}\big{(}\mathcal{C}^i_{k,\ell}-\mathcal{D}^i_{k,\ell}\widetilde{\mathbb{O}}^{\dagger}_{\ell}\widetilde{\mathbb{L}}_\ell\big{)}^T \widetilde{S}_{k,\ell+1}d^j_{k,\ell} -\widetilde{\mathbb{L}}_\ell^T\widetilde{\mathbb{O}}^{\dagger}_{\ell}\rho_{k,\ell} +q_{k,\ell},\\[1mm]
\pi_{k,N}=g_k, \quad
\ell\in \mathbb{T}_k,
\end{array}
\right.\\
\widetilde{\mathbb{O}}_k\widetilde{\mathbb{O}}_k^\dagger\widetilde{\theta}_k=\widetilde{\theta}_k,\quad k\in \mathbb{T}
\end{array}
\right.
\end{eqnarray}
are solvable in the sense of
\begin{eqnarray*}
\widetilde{\mathbb{O}}_k\succeq  0,~~~\widetilde{\mathbb{O}}_k\widetilde{\mathbb{O}}_k^\dagger\widetilde{\mathbb{L}}_k=\widetilde{\mathbb{L}}_k, ~~~
\widetilde{\mathbb{O}}_k\widetilde{\mathbb{O}}_k^\dagger\widetilde{\theta}_k=\widetilde{\theta}_k,~~~k\in \mathbb{T},
\end{eqnarray*}
where
\begin{eqnarray*}\label{Sec3:Theorem-Necessary-O-L}
\left\{
\begin{array}{l}
\widetilde{\mathbb{O}}_{k}=\mathcal{R}_{k,k}+\mathcal{B}^T_{k,k}\widetilde{\mathcal{S}}_{k,k+1}\mathcal{B}_{k,k}+\sum_{i,j=1}^p\delta_k^{ij}(\mathcal{D}^i_{k,k})^T\widetilde{S}_{k,k+1} \mathcal{D}^j_{k,k},\\[1mm]
\widetilde{\mathbb{L}}_{k}=\mathcal{B}^T_{k,k}\widetilde{\mathcal{S}}_{k,k+1}\mathcal{A}_{k,k}+\sum_{i,j=1}^p\delta_k^{ij}(\mathcal{D}^i_{k,k})^T\widetilde{S}_{k,k+1}\mathcal{C}^j_{k,k} +\mathcal{B}_{k,k}^T\widetilde{U}_{k,k+1},\\[1mm]
\widetilde{\theta}_{k}=\mathcal{B}^T_{k,k}\widetilde{\mathcal{S}}_{k,k+1}f_{k,k} +\sum_{i,j=1}^p\delta_k^{ij}(\mathcal{D}^i_{k,k})^T\widetilde{S}_{k,k+1}d^j_{k,k}
+\mathcal{B}^T_{k,k}\widetilde{\pi}_{k,k+1}+\rho_{k,k},\quad \\[1mm]
k\in \mathbb{T}
\end{array}
\right.
\end{eqnarray*}
with
\begin{eqnarray*}
\left\{
\begin{array}{l}
\widetilde{U}_{k,\ell}=(\mathcal{A}_{k,\ell}-\mathcal{B}_{k,\ell}
\widetilde{\mathbb{O}}^{\dagger}_{\ell}\widetilde{\mathbb{L}}_\ell)\widetilde{U}_{k,\ell+1},\\
\widetilde{U}_{k,N}=F_k,~~\ell\in \mathbb{T}_{k}, ~~k\in \mathbb{T},
\end{array}
\right.
\end{eqnarray*}
and
\begin{eqnarray*}
\widetilde{\beta}_{k,\ell}=-\widetilde{\mathbb{L}}_\ell^T\widetilde{\mathbb{O}}^{\dagger}_{\ell}\mathcal{R}_{k,\ell}+\big{(}\mathcal{A}_{k,\ell}-\mathcal{B}_{k,\ell}\widetilde{\mathbb{O}}_\ell^\dagger\widetilde{\mathbb{L}}_\ell\big{)}^T \widetilde{\mathcal{S}}_{k,\ell+1}\mathcal{B}_{k,\ell}+\sum_{i,j=1}^p\delta_\ell^{ij}\big{(}\mathcal{C}^i_{k,\ell}-\mathcal{D}^i_{k,\ell}\widetilde{\mathbb{O}}_\ell^\dagger\widetilde{\mathbb{L}}_\ell\big{)}^T \widetilde{S}_{k,\ell+1}\mathcal{D}^j_{k,\ell}, ~~ \ell\in \mathbb{T}_k.
\end{eqnarray*}

\item[iii)] There exists a pair $(\Phi, v)\in l^2(\mathbb{T}; \mathbb{R}^{m\times n}) \times  l^2_{\mathcal{F}}(\mathbb{T};\mathbb{R}^m)$ such that for any $(t,x)$ with $t\in \mathbb{T}$ and $x\in l^2_{\mathcal{F}}(t; \mathbb{R}^n)$, $(\Phi, v)|_{\mathbb{T}_t}$ is a linear feedback equilibrium strategy of Problem (LQ)$_{tx}$. Here, $(\Phi, v)|_{\mathbb{T}_t}$ is the restriction of $(\Phi, v)$ on ${\mathbb{T}_t}$.

\item[iv)] There exists a $\psi\in \mathbb{F}_\mathbb{T}$ such that for any $(t,x)$ with $t\in \mathbb{T}$ and $x\in l^2_{\mathcal{F}}(t;\mathbb{R}^n)$, $\psi|_{\mathbb{T}_t}$ is a feedback equilibrium strategy of Problem (LQ)$_{tx}$. Here, $\psi|_{\mathbb{T}_t}$ is the restriction of $\psi$ on ${\mathbb{T}_t}$.

\end{itemize}

Furthermore, under any of the above conditions, the pair $(\Phi^t, v^t)$ with
$
\Phi^t=\{-\widetilde{\mathbb{O}}^{\dagger}_{k} \widetilde{\mathbb{L}}_{k},~k\in \mathbb{T}_t\}$,
$v^t=\{-\widetilde{\mathbb{O}}^{\dagger}_{k} \widetilde{\theta}_{k}$, $k\in \mathbb{T}_t\}$
is a feedback equilibrium strategy of Problem (LQ)$_{tx}$.
\end{theorem}

\emph{Proof.} See Appendix \ref{appendix:Theorem-N-S}. \hfill$\square$

We now consider the mixed equilibrium solution. If it exists, we have some freedom to select the pure-feedback-strategy part of the mixed equilibrium solution, as pointed out in Remark \ref{Remark-Phi}. In Theorem \ref{Theorem-all-mixed}, we have the necessary and sufficient condition to ensure the existence of a mixed equilibrium solution for all of the initial pairs. Because different initial pairs may correspond to different pure-feedback-strategy parts of the mixed equilibrium solution, the condition of Theorem \ref{Theorem-all-mixed} is for the case that specifies the pure-feedback-strategy part $\Phi$.

\begin{theorem}\label{Theorem-all-mixed}
The following statements are equivalent:

\begin{itemize}
\item[i)] For any $(t,x)$ with $t\in \mathbb{T}$ and $x\in l^2_{\mathcal{F}}(t; \mathbb{R}^n)$, Problem (LQ)$_{tx}$ admits a mixed equilibrium solution and the pure-feedback-strategy part is $\Phi|_{\mathbb{T}_t}$. Here, $\Phi \in l^2(\mathbb{T};\mathbb{R}^{m\times n})$ and $\Phi|_{\mathbb{T}_t}$ is the restriction of $\Phi$ on $\mathbb{T}_t$.

\item[ii)] There exists  $\Phi \in l^2(\mathbb{T};\mathbb{R}^{m\times n})$ such that the following difference equations
 \begin{eqnarray}\label{S-all-mixed}
\left\{
\begin{array}{l}
\left\{
\begin{array}{l}
S_{k,\ell}=Q_{k,\ell}+\Phi_\ell^TR_{k,\ell}\Phi_{\ell}+\big{(}A_{k,\ell} +B_{k,\ell}\Phi_\ell\big{)}^TS_{k,\ell+1}\big{(}A_{k,\ell}+B_{k,\ell}\Phi_\ell\big{)}\\[1mm]
\hphantom{S_{k,\ell}=}+\sum_{i,j=1}^p\delta_\ell^{ij}\big{(}C^i_{k,\ell} +D^i_{k,\ell}\Phi_\ell\big{)}^TS_{k,\ell+1}\big{(}C^j_{k,\ell}+D^j_{k,\ell}\Phi_\ell\big{)},\\[1mm]
{\mathcal{S}}_{k,\ell}
={\mathcal{Q}}_{k,\ell}+\Phi_\ell^T{\mathcal{R}}_{k,\ell}\Phi_{\ell}+\big{(}\mathcal{A}_{k,\ell} +\mathcal{B}_{k,\ell}\Phi_\ell\big{)}^T\mathcal{S}_{k,\ell+1}({\mathcal{A}}_{k,\ell}+{\mathcal{B}}_{k,\ell}\Phi_{\ell}) \\[1mm]
\hphantom{\bar{S}_{k,\ell}=}+\sum_{i,j=1}^p\delta_\ell^{ij}\big{(}\mathcal{C}^i_{k,\ell}+\mathcal{D}^i_{k,\ell}\Phi_\ell\big{)}^T S_{k,\ell+1}({\mathcal{C}}^j_{k,\ell}+{\mathcal{D}}^j_{k,\ell}\Phi_{\ell}),\\[1mm]
%
%
%
S_{k,N}=G_k,~~{\mathcal{S}}_{k,N}=G_k+\bar{G}_k,~~\ell\in \mathbb{T}_k, \quad
\end{array}
\right.\\[1mm]
{\mathbb{O}}_k=\mathcal{R}_{k,k}+\mathcal{B}^T_{k,k}\mathcal{S}_{k,k+1}\mathcal{B}_{k,k}+\sum_{i,j=1}^p\delta_k^{ij}(\mathcal{D}^i_{k,k})^T{S}_{k,k+1} \mathcal{D}^j_{k,k}\succeq 0,\quad k\in \mathbb{T},
\end{array}
\right.
\end{eqnarray}
\begin{eqnarray}\label{T-all-mixed}
\left\{
\begin{array}{l}
\left\{
\begin{array}{l}
T_{k,\ell}=\Big{\{}\Phi_{\ell}^TR_{k,\ell}+\big{(}A_{k,\ell}+B_{k,\ell}\Phi_{\ell}\big{)}^TS_{k,\ell+1}B_{k,\ell}\\[1mm]
\hphantom{T_{k,\ell}=}  +\sum_{i,j=1}^p\delta_\ell^{ij}\big{(}C^i_{k,\ell}+D^i_{k,\ell}\Phi_{\ell}\big{)}^TS_{k,\ell+1}D^j_{k,\ell}\Big{\}}\Gamma_\ell\\[1mm]
\hphantom{T_{k,\ell}=}+\big{(}A_{k,\ell}+B_{k,\ell}\Phi_{\ell}\big{)}^TT_{k,\ell+1}\big{(}\mathcal{A}_{\ell,\ell} +\mathcal{B}_{\ell,\ell}\Phi_{\ell}+\mathcal{B}_{\ell,\ell}\Gamma_{\ell} \big{)}\\[1mm]
\hphantom{T_{k,\ell}=}+\sum_{i,j=1}^p\delta_\ell^{ij}\big{(}C^i_{k,\ell}+D^i_{k,\ell}\Phi_{\ell}\big{)}^TT_{k,\ell+1}\big{(}\mathcal{C}^j_{\ell,\ell} +\mathcal{D}^j_{\ell,\ell}\Phi_{\ell}+\mathcal{D}^j_{\ell,\ell}\Gamma_{\ell} \big{)},\\[1mm]
{\mathcal{T}}_{k,\ell}=\Big{\{}\Phi_{\ell}^T{{\mathcal{R}}}_{k,\ell}+\big{(}\mathcal{A}_{k,\ell}+\mathcal{B}_{k,\ell}\Phi_{\ell}\big{)}^T \mathcal{S}_{k,\ell+1}{\mathcal{B}}_{k,\ell}\\[1mm]
\hphantom{\bar{T}_{k,\ell}=}+\sum_{i,j=1}^p\delta_\ell^{ij} \big{(}\mathcal{C}^i_{k,\ell}+\mathcal{D}^i_{k,\ell}\Phi_{\ell}\big{)}^TS_{k,\ell+1}{\mathcal{D}}^j_{k,\ell}\Big{\}}\Gamma_{\ell}\\[1mm]
\hphantom{\bar{T}_{k,\ell}=}+\big{(}\mathcal{A}_{k,\ell}+\mathcal{B}_{k,\ell}\Phi_\ell\big{)}^T{\mathcal{T}}_{k,\ell+1}\big{(}\mathcal{A}_{\ell,\ell} +\mathcal{B}_{\ell,\ell}\Phi_\ell +\mathcal{B}_{\ell,\ell}\Gamma_\ell \big{)}\\[1mm]
%
%
\hphantom{\bar{T}_{k,\ell}=}+\sum_{i,j=1}^p\delta_\ell^{ij}\big{(}{\mathcal{C}}^i_{k,\ell}+{\mathcal{D}}^i_{k,\ell}\Phi_\ell\big{)}^TT_{k,\ell+1} \big{(}\mathcal{C}^j_{\ell,\ell} +\mathcal{D}^j_{\ell,\ell}\Phi_\ell +\mathcal{D}^j_{\ell,\ell}\Gamma_\ell \big{)}\\[1mm]
T_{k,N}=0,~~{\mathcal{T}}_{k,N}=0, \quad
\ell\in \mathbb{T}_k,
\end{array}
\right.\\
\mathcal{O}_k\mathcal{O}_k^\dagger \mathcal{L}_k=\mathcal{L}_k, \quad k\in \mathbb{T},
\end{array}
\right.
\end{eqnarray}
and
\begin{eqnarray}\label{pi-all-mixed}
\left\{
\begin{array}{l}
\left\{
\begin{array}{l}
\pi_{k,\ell}=-\beta_{k,\ell}\mathcal{O}^\dagger_\ell\theta_\ell+\big{(}\mathcal{A}_{k,\ell}+\mathcal{B}_{k,\ell}\Phi_{\ell}\big{)}^T\big{(}\mathcal{S}_{k,\ell+1}f_{k,\ell}+\pi_{k,\ell+1}\big{)} \\[1mm]
\hphantom{\pi_{k,\ell}=}+\sum_{i,j=1}^p\delta_\ell^{ij}\big{[}\big{(}\mathcal{C}^i_{k,\ell}+\mathcal{D}^i_{k,\ell}\Phi_{\ell}\big{)}^TS_{k,\ell+1}d^j_{k,\ell} +\big{(}\mathcal{C}^i_{k,\ell}+\mathcal{D}^i_{k,\ell}\Phi_\ell\big{)}^TT_{k,\ell+1}d^j_{\ell,\ell}\\[1mm]
\hphantom{\pi_{k,\ell}=}+\big{(}\mathcal{A}_{k,\ell}+\mathcal{B}_{k,\ell}\Phi_\ell\big{)}^T\mathcal{T}_{k,\ell+1} f_{\ell,\ell}+\Phi_{\ell}^T\rho_{k,\ell} +q_{k,\ell},\\[1mm]
\pi_{k,N}=g_k, \quad
\ell\in \mathbb{T}_k,
\end{array}
\right.\\
\mathcal{O}_k\mathcal{O}_k^\dagger \theta_k=\theta_k,\quad k\in \mathbb{T}
\end{array}
\right.
\end{eqnarray}
are solvable in the sense of
${\mathbb{O}}_k\succeq  0,~~~\mathcal{O}_k\mathcal{O}_k^\dagger \mathcal{L}_k=\mathcal{L}_k,~~~\mathcal{O}_k\mathcal{O}_k^\dagger \theta_k=\theta_k$, $k\in \mathbb{T}$,
where
\begin{eqnarray*}
\left\{
\begin{array}{l}
\mathcal{O}_{k}=\mathcal{R}_{k,k}+\mathcal{B}^T_{k,k}\big{(}\mathcal{S}_{k,k+1}+\mathcal{T}_{k,k+1}\big{)}\mathcal{B}_{k,k} +\sum_{i,j=1}^p\delta_k^{ij}(\mathcal{D}^i_{k,k})^T\big{(}{S}_{k,k+1}+T_{k,k+1}\big{)}\mathcal{D}^j_{k,k},\\[1mm]
\mathcal{L}_{k}=\mathcal{B}^T_{k,k}\big{(}\mathcal{S}_{k,k+1}
+\mathcal{T}_{k,k+1}\big{)}\mathcal{A}_{k,k}+\sum_{i,j=1}^p\delta_k^{ij}(\mathcal{D}^i_{k,k})^T \big{(}S _{k,k+1}+T_{k,k+1}\big{)}\mathcal{C}^j_{k,k}+\mathcal{B}_{k,k}^TU_{k,k+1},\\[1mm]
\theta_{k}=\mathcal{B}^T_{k,k}\big{(}\mathcal{S}_{k,k+1}+\mathcal{T}_{k,k+1}\big{)}f_{k,k}
+\sum_{i,j=1}^p\delta_k^{ij}(\mathcal{D}^i_{k,k})^T\big{(}S_{k,k+1}+T_{k,k+1}\big{)}d^j_{k,k}
+\mathcal{B}^T_{k,k}\pi_{k,k+1}+\rho_{k,k},\\[1mm]
k\in \mathbb{T}
\end{array}
\right.
\end{eqnarray*}
and $\Gamma_k=-\mathcal{O}_k^{\dagger}\mathcal{L}_k-\Phi_k$, $k\in \mathbb{T}$
with
\begin{eqnarray*}
\left\{
\begin{array}{l}
U_{k,\ell}=(\mathcal{A}_{k,\ell}+\mathcal{B}_{k,\ell}\Phi_{\ell})U_{k,\ell+1},\\
U_{k,N}=F_k, ~~\ell\in \mathbb{T}_{k}, ~~k\in \mathbb{T},
\end{array}
\right.
\end{eqnarray*}
and
\begin{eqnarray*}
\begin{array}{l}
\beta_{k,\ell}=\Phi_{\ell}^T\mathcal{R}_{k,\ell}+\big{(}\mathcal{A}_{k,\ell}+\mathcal{B}_{k,\ell}\Phi_{\ell}\big{)}^T\big{[} \mathcal{S}_{k,\ell+1}\mathcal{B}_{k,\ell}+\mathcal{T}_{k,\ell+1}\mathcal{B}_{\ell,\ell} \big{]}\\[1mm]
\hphantom{\beta_{k,\ell}=}+\sum_{i,j=1}^p\delta_\ell^{ij}\big{(}\mathcal{C}^i_{k,\ell}+\mathcal{D}^i_{k,\ell}\Phi_{\ell}\big{)}^T \big{[}S_{k,\ell+1}\mathcal{D}^j_{k,\ell}+T_{k,\ell+1}\mathcal{D}^j_{\ell,\ell}\big{]},\quad \ell\in \mathbb{T}_k.
\end{array}
\end{eqnarray*}
For any $t\in \mathbb{T}$ and $x\in l^2_\mathcal{F}(t; \mathbb{R}^{n})$, let
$v^{t,x}_k=-\big{(}\mathcal{O}_k^{\dagger}\mathcal{L}_k+\Phi_k\big{)}X^{t,x,*}_k-\mathcal{O}^\dagger_k\theta_k$,
$k\in \mathbb{T}_t$,
where
\begin{eqnarray*}
\left\{\begin{array}{l}
X^{t,x,*}_{k+1} =\big{[}\big{(}\mathcal{A}_{k,k}-\mathcal{B}_{k,k}\mathcal{O}_k^{\dagger}\mathcal{L}_k \big{)}X^{t,x,*}_{k}-\mathcal{B}_{k,k}\mathcal{O}^\dagger_k\theta_k +f_{k,k}\big{]}\\[1mm]
\hphantom{X^{t,x,*}_{k+1} =}+\sum_{i=1}^p\big{[}\big{(}\mathcal{C}^i_{k,k}-\mathcal{D}^i_{k,k}\mathcal{O}_k^{\dagger}\mathcal{L}_k \big{)} X^{t,x,*}_{k}-\mathcal{D}^i_{k,k}\mathcal{O}^\dagger_k\theta_k+d^i_{k,k}\big{]} w^i_k,\\[1mm]
X^{t,x,*}_{t} = x,~~k\in \mathbb{T}_t.
\end{array}
\right.
\end{eqnarray*}

\end{itemize}

Under condition ii) and for any $(t,x)$ with $t\in \mathbb{T}$ and $x\in l^2_{\mathcal{F}}(t; \mathbb{R}^n)$, $(\Phi|_{\mathbb{T}_t}, v^{t,x})$ is a mixed equilibrium solution of Problem (LQ)$_{tx}$.

\end{theorem}

\emph{Proof}. \emph{ii)$\Rightarrow$i)}. This follows from Theorem \ref{Sec3:Theorem-Necessary}.

\emph{i)$\Rightarrow$ii)}. In this case, for any $(t,x)$ with $t\in \mathbb{T}$ and $x\in l^2_{\mathcal{F}}(t; \mathbb{R}^n)$, the pure-feedback-strategy part of the mixed equilibrium solution is $\Phi|_{\mathbb{T}_t}$. Note that (\ref{stationary-condition-closed-1}) is equivalent to ${\mathcal{O}}_k {\mathcal{O}}_k^\dagger\big{(}{\mathcal{L}}_k{X}^{t,x,*}_k+{\theta}_k \big{)}={\mathcal{L}}_k{X}^{t,x,*}_k+{\theta}_k$, $k\in \mathbb{T}_t$.
Letting $k=t$ and taking different $x's$, we have
${\mathcal{O}}_t {\mathcal{O}}_t^\dagger{\mathcal{L}}_t={\mathcal{L}}_t$,
${\mathcal{O}}_t {\mathcal{O}}_t^\dagger{\theta}_t={\theta}_t$.
Therefore, (\ref{Sec3:Theorem-Necessary-S}), (\ref{T-all-mixed}), and (\ref{pi-all-mixed}) are solvable.   \hfill $\square$

To end this section, we pose the following assumption.

\textbf{(H)} \quad $Q_{t,k}, Q_{t,k}+\bar{Q}_{t,k}, G_t, G_t+\bar{G}_k\succeq 0, R_{t,k}, R_{t,k}+\bar{R}_{t,k}\succ0, t\in \mathbb{T},  k\in \mathbb{T}_t$.

The following result is straightforward.

\begin{theorem}

Letting (H) hold, then $\widetilde{\mathbb{O}}_k, k\in \mathbb{T}$, are all positive definite. Furthermore, for any $t\in \mathbb{T}$ and any $x\in l^2_{\mathcal{F}}(t; \mathbb{R}^n)$, Problem (LQ)$_{tx}$ admits a unique feedback equilibrium strategy $(\Phi, v)$ with
\begin{eqnarray*}
%
\Phi=\{-\widetilde{\mathbb{O}}^{-1}_{k} \widetilde{\mathbb{L}}_{k},~k\in \mathbb{T}_t\},~~~
v=\{-\widetilde{\mathbb{O}}^{-1}_{k} \widetilde{\theta}_{k},~k\in \mathbb{T}_t\}.
\end{eqnarray*}

\end{theorem}

\emph{Proof}. A simple calculation shows that (\ref{S-all-feedback}) is equal to
\begin{eqnarray}\label{Sec3:corollary-S}
\left\{
\begin{array}{l}
\left\{
\begin{array}{l}
\widetilde{S}_{k,\ell}=Q_{k,\ell}+\widetilde{\mathbb{L}}_\ell^T\widetilde{\mathbb{O}}_\ell^\dagger R_{k,\ell}\widetilde{\mathbb{O}}_\ell^\dagger\widetilde{\mathbb{L}}_\ell+\big{(}A_{k,\ell} -B_{k,\ell}\widetilde{\mathbb{O}}_\ell^\dagger\widetilde{\mathbb{L}}_\ell\big{)}^TS_{k,\ell+1}\big{(}A_{k,\ell}-B_{k,\ell}\widetilde{\mathbb{O}}_\ell^\dagger \widetilde{\mathbb{L}}_\ell\big{)}\\[1mm]
\hphantom{\widetilde{S}_{k,\ell}=}+\sum_{i,j=1}^p\delta_\ell^{ij}\big{(}C^i_{k,\ell} -D^i_{k,\ell}\widetilde{\mathbb{O}}_\ell^\dagger\widetilde{\mathbb{L}}_\ell\big{)}^TS_{k,\ell+1}\big{(}C^j_{k,\ell}-D^j_{k,\ell}\widetilde{\mathbb{O}}_\ell^\dagger \widetilde{\mathbb{L}}_\ell\big{)},\\[1mm]
\widetilde{\mathcal{S}}_{k,\ell}
={\mathcal{Q}}_{k,\ell}+\widetilde{\mathbb{L}}_\ell^T\widetilde{\mathbb{O}}_\ell^\dagger {\mathcal{R}}_{k,\ell}\widetilde{\mathbb{O}}_\ell^\dagger\widetilde{\mathbb{L}}_\ell+\big{(}\mathcal{A}_{k,\ell} -\mathcal{B}_{k,\ell}\widetilde{\mathbb{O}}_\ell^\dagger\widetilde{\mathbb{L}}_\ell\big{)}^T\widetilde{\mathcal{S}}_{k,\ell+1}({\mathcal{A}}_{k,\ell} -{\mathcal{B}}_{k,\ell}\widetilde{\mathbb{O}}_\ell^\dagger\widetilde{\mathbb{L}}_\ell) \\[1mm]
\hphantom{\bar{S}_{k,\ell}=}+\sum_{i,j=1}^p\delta_\ell^{ij}\big{(}\mathcal{C}^i_{k,\ell}-\mathcal{D}^i_{k,\ell}\widetilde{\mathbb{O}}_\ell^\dagger\widetilde{\mathbb{L}}_\ell \big{)}^TS_{k,\ell+1}({\mathcal{C}}^j_{k,\ell}- {\mathcal{D}}^j_{k,\ell}\widetilde{\mathbb{O}}_\ell^\dagger\widetilde{\mathbb{L}}_\ell),\\[1mm]
%
%
%
\widetilde{S}_{k,N}=G_k,~~\widetilde{\mathcal{S}}_{k,N}=G_k+\bar{G}_k,
\quad
\ell\in \mathbb{T}_k,
\\[1mm]
%
%
\end{array}
\right.\\
\widetilde{\mathbb{O}}_k\succeq 0, \quad
\widetilde{\mathbb{O}}_k\widetilde{\mathbb{O}}_k^\dagger\widetilde{\mathbb{L}}_k=\widetilde{\mathbb{L}}_k, \quad k\in \mathbb{T}_t.
\end{array}
\right.
\end{eqnarray}
Due to (H), we have that $\widetilde{S}_{k,\ell}, \widetilde{\mathbb{S}}_{k,\ell}\succeq 0, \widetilde{\mathbb{O}}_k\succ0, k\in \mathbb{T}, \ell\in \mathbb{T}_k$. Therefore, (\ref{S-all-feedback}) and (\ref{pi-all-feedback}) are solvable. This completes the proof. \hfill $\square$

\section{An example}\label{Section-example}

Consider a discrete-time stochastic LQ problem, whose system dynamics and cost functional are given, respectively, by
\begin{eqnarray*}
\left\{
\begin{array}{l}
X_{k+1}=(A_kX_k+B_ku_k)+D_ku_kw_k, \\
X_t=x,~~t\in \{0,1,2,3\},~~k\in \{t,...,3\},
\end{array}
\right.
\end{eqnarray*}
and
\begin{eqnarray*}
J(t,x;u)=\sum_{k=t}^{3}\mathbb{E}_t\big{[}X_k^TQ_{k}X^t_k+ u_k^TR_{k}u_k\big{]}+\mathbb{E}_t\big{[}X_4^TGX_4\big{]}
+(\mathbb{E}_tX_4)^T \bar{G}\mathbb{E}_tX_{4}, 
\end{eqnarray*}
where
\begin{eqnarray*}
&&\hspace{-2em}A_0=\left[
\begin{array}{cc}
1& 0.4\\0.3 &2
\end{array}
\right],~~A_1=\left[
\begin{array}{cc}
1.102& -0.24\\ 0.53& 1.89
\end{array}
\right],~~A_2=\left[
\begin{array}{cc}
1.89& 0.49\\0& 1.75
\end{array}
\right],~~A_3=\left[
\begin{array}{cc}
0.8& -0.4\\0.2& 0.7
\end{array}
\right],\\[1mm]
&&\hspace{-2em}B_0=\left[
\begin{array}{c}
1.2\\ -0.5
\end{array}
\right],~~B_1=\left[
\begin{array}{c}
1\\ 1
\end{array}
\right],~~B_2=\left[
\begin{array}{c}
1.2\\ 0.2
\end{array}
\right],~~B_3=\left[
\begin{array}{c}
1\\ 0.3
\end{array}
\right],~~D_0=\left[
\begin{array}{c}
1\\ 0.3
\end{array}
\right],~~D_1=\left[
\begin{array}{c}
1\\ 0.4
\end{array}
\right],\\[1mm]
&&\hspace{-2em}D_2=\left[
\begin{array}{c}
0.45\\ 0.25
\end{array}
\right],~D_3=\left[
\begin{array}{c}
0.52\\ 0
\end{array}
\right],~
Q_0=\left[
\begin{array}{cc}
3& 0.5\\0.5 &-2
\end{array}
\right],~
Q_1=\left[
\begin{array}{cc}
2& -0.65\\-0.65 &0
\end{array}
\right],
~
Q_2=\left[
\begin{array}{cc}
0.5& 0.5\\0.5& -2
\end{array}
\right],\\[1mm]
&&\hspace{-2em}Q_3=\left[
\begin{array}{cc}
-0.1& 0\\0 &-0.75
\end{array}
\right],~~R_0=0,~~R_1=-2.5,~~R_2=1,~~R_3=-0.5,\\
&&\hspace{-2em}G=\left[
\begin{array}{cc}
1& -0.1\\-0.1& 1
\end{array}
\right],~~\bar{G}=\left[
\begin{array}{cc}
-0.3 &0\\0& -0.3
\end{array}
\right],
\end{eqnarray*}
and $\{w_k, k=0,1,2,3\}$ is a martingale difference with constant second-order conditional moment $\mathbb{E}_{k}(w_k^2)=1, k=0,1,2,3$.

{\underline{Open-loop equilibrium control}}

For this LQ problem, by performing the iteration (\ref{S-all}), we have
$\widehat{\mathbb{O}}_0=8.7645$, $\widehat{\mathbb{O}}_1=-0.4783$, $\widehat{\mathbb{O}}_2=1.6935$, $\widehat{\mathbb{O}}_3=0.7193$.
Because $\widehat{\mathbb{O}}_1=-0.4783<0$, for $(t,x)$ with $t=0,1$ and $x\in l^2_{\mathcal{F}}(0; \mathbb{R}^2)$ or $x\in l^2_{\mathcal{F}}(1; \mathbb{R}^2)$ and based on Corollary \ref{Corollary--open}, the open-loop equilibrium control of this LQ problem must not exist.

{\underline{Feedback equilibrium strategy}}

By performing the iteration (\ref{S-all-feedback}), we have
$\widetilde{\mathbb{O}}_0=-11.0590$, $\widetilde{\mathbb{O}}_1=20.5335$, $\widetilde{\mathbb{O}}_2=-0.5593$, $\widetilde{\mathbb{O}}_3=0.4734$.
Because $\widetilde{\mathbb{O}}_0<0, \widetilde{\mathbb{O}}_2<0$, for $(t,x)$ with $t=0,1,2$ and $x\in l^2_{\mathcal{F}}(0; \mathbb{R}^2)$ or $x\in l^2_{\mathcal{F}}(1; \mathbb{R}^2)$, and based on Corollary \ref{Corollary--feedback} and Theorem \ref{Sec3:Theorem-N-S}, the feedback equilibrium strategy of this LQ problem must not exist.

{\underline{Mixed equilibrium solution}}

We use the command ``randn'' of MATLAB to randomly generate $\Phi=\{\Phi_k, k=0,1,2,3\}$. Note that $\Phi_k\in \mathbb{R}^{1\times 2}, \mathbb{O}_k, \mathcal{O}_k\in \mathbb{R}^1, k=0,1,2,3$, and let
\begin{eqnarray*}
\psi=[\Phi_0^T,~ \Phi_1^T,~ \Phi_2^T, ~ \Phi_3^T]^T,~~~\mathbb{O}=(\mathbb{O}_0, \mathbb{O}_1, \mathbb{O}_2, \mathbb{O}_3),~~~\mathcal{O}=(\mathcal{O}_0, \mathcal{O}_1, \mathcal{O}_2, \mathcal{O}_3).
\end{eqnarray*}
By performing the iterations (\ref{Sec3:Theorem-Necessary-S})-(\ref{T-all-mixed})-(\ref{pi-all-mixed}), we select 10 $\psi$s and get the corresponding $\mathbb{O}$s and $\mathcal{O}$s,
\begin{eqnarray*}
&&\hspace{-3em}\psi=\left[
\begin{array}{cc}
1.4090  &  1.4172\\ -0.1241  &  1.4897 \\ 0.7147 &  -0.2050 \\ 0.7254 &  -0.0631
\end{array}
\right],~~~
\begin{array}{l}
\mathbb{O}=(42.1215,~ 21.2758,~ 3.1578,~ 0.4734),\\[2mm]
\mathcal{O}=(-2.1680,~ -10.6485,~ 0.4740, ~ 0.4734),
\end{array}\\
&&\hspace{-3em}\psi=\left[
\begin{array}{cc}
 0.7269 &  -0.3034\\
 0.4889  &  1.0347 \\
0.7172   & 1.6302 \\
 0.6715  & -1.2075
\end{array}
\right],~~~
\begin{array}{l}
\mathbb{O}=(106.9951,~ 28.5844,~ 2.3227,~ 0.4734),\\[2mm]
\mathcal{O}=(-2.4665,~ -10.5353,~ 0.4860, ~ 0.4734),
\end{array}\\
&&\hspace{-3em}\psi=\left[
\begin{array}{cc}
 0.3192  &  0.3129\\
   -0.1022 &  -0.2414 \\
  1.3703 &  -1.7115\\
0.3252 &  -0.7549
\end{array}
\right],~~~
\begin{array}{l}
\mathbb{O}=(35.1212,~ 1.8350,~ 1.7640,~ 0.4734),\\[2mm]
\mathcal{O}=(-0.8786,~ -9.8337,~ 0.4876, ~ 0.4734),
\end{array}\\
&&\hspace{-3em}\psi=\left[
\begin{array}{cc}
   -0.7648  & -1.4023\\
   -0.1924  &  0.8886\\
   -0.6156  &  0.7481\\
   -1.0616  &  2.3505
\end{array}
\right],~~~
\begin{array}{l}
\mathbb{O}=(20.1218,~ 2.2184,~ 0.8268,~ 0.4734),\\[2mm]
\mathcal{O}=(-1.3929,~ -9.2281,~ 0.4817, ~ 0.4734),
\end{array}\\
&&\hspace{-3em}\psi=\left[
\begin{array}{cc}
   -0.4390 &  -1.7947\\
   -0.0825 &  -1.9330\\
   -0.6669 &   0.1873\\
    0.7223 &   2.5855
\end{array}
\right],~~~
\begin{array}{l}
\mathbb{O}=(52.1877,~ 31.3899,~ 5.4614,~ 0.4734),\\[2mm]
\mathcal{O}=(-1.9877,~ -10.6047,~ 0.4485, ~ 0.4734),
\end{array}\\
&&\hspace{-3em}\psi=\left[
\begin{array}{cc}
    0.4900  &  0.7394\\
    0.3035  & -0.6003\\
    0.1001  & -0.5445\\
    0.8404  & -0.8880
\end{array}
\right],~~~
\begin{array}{l}
\mathbb{O}=(31.4336,~ 2.6274,~ 2.9500,~ 0.4734),\\[2mm]
\mathcal{O}=(-1.0423,~ -9.8977,~ 0.4799, ~ 0.4734),
\end{array}\\
&&\hspace{-3em}\psi=\left[
\begin{array}{cc}
    0.9610  &  0.1240\\
    1.3546  & -1.0722\\
   -2.1384  & -0.8396\\
    1.7119  & -0.1941
\end{array}
\right],~~~
\begin{array}{l}
\mathbb{O}=(429.0833,~ 38.2114,~ 6.7849,~ 0.4734),\\[2mm]
\mathcal{O}=(1.1514,~ -8.0070,~ 0.4581, ~ 0.4734),
\end{array}\\
&&\hspace{-3em}\psi=\left[
\begin{array}{cc}
    1.3790  & -1.0582\\
    2.9080  &  0.8252\\
   -0.1977  & -1.2078\\
    1.4367  & -1.9609
\end{array}
\right],~~~
\begin{array}{l}
\mathbb{O}=(112.4586,~ 3.2533,~ 4.0958,~ 0.4734),\\[2mm]
\mathcal{O}=(1.7922,~ -9.5504,~ 0.4799, ~ 0.4734),
\end{array}\\
&&\hspace{-3em}\psi=\left[
\begin{array}{cc}
   -1.1564 &  -0.5336\\
   -0.8314 &  -0.9792\\
   -1.7502 &  -0.2857\\
    0.0229 &  -0.2620
\end{array}
\right],~~~
\begin{array}{l}
\mathbb{O}=(7.6517,~ 5.3349,~ 1.3968,~ 0.4734),\\[2mm]
\mathcal{O}=(-0.8077,~ -8.7128,~ 0.4881, ~ 0.4734),
\end{array}\\
&&\hspace{-3em}\psi=\left[
\begin{array}{cc}
    0.0513  &  0.8261\\
   -0.3031  &  0.0230\\
   -0.1952  & -0.2176\\
    0.6601  & -0.0679
\end{array}
\right],~~~
\begin{array}{l}
\mathbb{O}=(11.0638,~ 4.5685,~ 2.9632,~ 0.4734),\\[2mm]
\mathcal{O}=(-1.2944,~ -9.9027,~ 0.4752, ~ 0.4734).
\end{array}
\end{eqnarray*}
For all 10 cases, $\mathbb{O}_k, k=0,1,2,3$, are all positive, and $\mathcal{O}_k, k=0,1,2,3$ are all invertible. Then, due to Theorem \ref{Theorem-all-mixed}, for any $(t,x)$ with $t\in \{0,1,2,3\}$ and $x\in l^2_{\mathcal{F}}(t; \mathbb{R}^2)$, the above 10 cases correspond to 10 mixed equilibrium solutions of the considered LQ problem, which can be easily constructed from Theorem \ref{Theorem-all-mixed}. For example, with the last $\psi$ given above, the mixed equilibrium solution is as follows. Let
\begin{eqnarray*}
\Phi_0=[0.0513~~  0.8261],~\Phi_1=[-0.3031 ~~  0.0230],~\Phi_2=[-0.1952 ~~ -0.2176],~\Phi_3=[0.6601 ~~ -0.0679],
\end{eqnarray*}
and
$$v^{0,x}_k=-\big{(}\mathcal{O}_k^{\dagger}\mathcal{L}_k+\Phi_k\big{)}X^{0,x,*}_k,~~~k\in \{0,1,2,3\}$$
with
$$
\left\{
\begin{array}{l}
X^{0,x,*}_{k+1} =\big{(}{A}_{k}-{B}_{k}\mathcal{O}_k^{\dagger}\mathcal{L}_k \big{)}X^{0,x,*}_{k}-{D}_{k}\mathcal{O}_k^{\dagger}\mathcal{L}_k  X^{0,x,*}_{k} w_k,\\[1mm]
X^{0,x,*}_{0} = x,~~~ k\in \{0,1,2,3\},
\end{array}
\right.
$$
and
\begin{eqnarray*}
&&-\mathcal{O}_0^{\dagger}\mathcal{L}_0=[1.4347,~ 4.2547], -\mathcal{O}_1^{\dagger}\mathcal{L}_1=[-0.3247,~ -0.5193],\\[1mm]
&&-\mathcal{O}_2^{\dagger}\mathcal{L}_2=[1.4568,~ 0.3845],
-\mathcal{O}_3^{\dagger}\mathcal{L}_3=[-1.1787,~ 0.4035].
\end{eqnarray*}
Then, $(\Phi, v^{0,x})$ is a mixed equilibrium solution of this LQ problem for the initial pair $(0,x)$, where $\Phi=\{\Phi_k, k=0,1,2,3\}$.

\section{Summary}\label{Section-Summary}

In this paper, the notion of mixed equilibrium solution is introduced for the time-inconsistent discrete-time
mean-field stochastic LQ optimal control. For a pair of pure-feedback strategy and open-loop control, necessary and sufficient conditions are given to ensure that such a pair is a mixed equilibrium solution. On this basis, the open-loop equilibrium control and feedback equivalent strategy can be dealt with in a unified way.

Although we provide some relevant results, the theory for mixed equilibrium solution is far from mature. From the example in Section \ref{Section-example}, we know that a remarkable property of mixed equilibrium solution is its non-uniqueness. Thus, we propose that the following topics warrant further study:
\begin{itemize}

\item[i)] Characterize the set of all of the mixed equilibrium solutions of Problem (LQ).

\item[ii)] Find the ``best'' mixed equilibrium solution, which should be the one under which the equilibrium value function attains its extreme.

\item[iii)] As a test, the multi-period mean-variance portfolio selection must be thoroughly investigated.

\item[iv)] Finally, the analysis should be extended beyond the realm of the LQ controls and to a continuous-time setting.
\end{itemize}



\appendix

\appendix

\section{Proof of Lemma \ref{Lemma-difference}}\label{appendix-0}         

From (\ref{Sec2:Defini-closed-X-1}) and (\ref{Sec3:barX}), we have
\begin{eqnarray*}
\left\{
\begin{array}{l}
\frac{\bar{X}^{k,\bar{u}_k,\lambda}_{\ell+1}-X^{k,\Phi}_{\ell+1}}{\lambda}=\big{(}A_{k,\ell}+B_{k,\ell}\Phi_\ell\big{)}\frac{\bar{X}_\ell^{k,\bar{u}_k,\lambda}-X^{k,\Phi}_\ell}{\lambda}+\big{(}\bar{A}_{k,\ell}+\bar{B}_{k,\ell}\Phi_\ell\big{)}\frac{\mathbb{E}_k{X}_\ell^{k,\lambda}-\mathbb{E}_kX^{k,u_k}_\ell}{\lambda}\\[1mm]
\hphantom{\frac{\bar{X}^{k,\bar{u}_k,\lambda}_{\ell+1}-X^{k,\Phi}_{\ell+1}}{\lambda}=}+\sum_{i=1}^p\Big{[}\big{(}C^i_{k,\ell}+D^i_{k,\ell}\Phi_\ell\big{)}\frac{\bar{X}_\ell^{k,\bar{u}_k,\lambda}-X^{k,\Phi}_\ell}{\lambda} +\big{(}\bar{C}^i_{k,\ell}+\bar{D}^i_{k,\ell}\Phi_\ell\big{)}\frac{\mathbb{E}_k{X}_\ell^{k,\lambda}-\mathbb{E}_kX^{k,u_k}_\ell}{\lambda}\Big{]}w^i_\ell,\\[1mm]
\frac{\bar{X}^{k,\bar{u}_k,\lambda}_{k+1}-{X}^{k,\Phi}_{k+1}}{\lambda}=\mathcal{B}_{k,k}\bar{u}_k+\sum_{i=1}^p\mathcal{D}^i_{k,k}\bar{u}_kw^i_k,\\[1mm]
\frac{\bar{X}_k^{k,\bar{u}_k,\lambda}-X^{k,\Phi}}{\lambda}=0,~~~\ell\in \mathbb{T}_{k+1}.
\end{array}
\right.
\end{eqnarray*}
Denote $\frac{\bar{X}_\ell^{k,\bar{u}_k,\lambda}-X^{k,\Phi}_\ell}{\lambda}$ by $\alpha^{k,\bar{u}_k}_\ell$. Then $\alpha^{k,\bar{u}_k}=\{\alpha^{k,\bar{u}_k}_\ell, \ell\in \widetilde{\mathbb{T}}_k\}$ satisfies (\ref{system-y-k}).
Obviously, we have $\bar{X}^{k, \bar{u}_k, \lambda}_\ell=X^{k,\Phi}_\ell
+\lambda \alpha^{k,\bar{u}_k}_\ell$, $\ell \in \mathbb{T}_k$. Then, we obtain
\begin{eqnarray}\label{Lemma-difference:J}
&&\hspace{-1.75em}J\big{(}k, X_k^{t,x,*}; (\Phi_k\cdot \bar{X}^{k,\bar{u}_k,\lambda}_k+v^{t,x}_k +\lambda \bar{u}_k,(\Phi \bar{X}^{k,\bar{u}_k,\lambda}+v^{t,x})|_{\mathbb{T}_{k+1}})\big{)}-J\big{(}k, X_k^{t,x,*}; (\Phi\cdot X^{k,\Phi}+v^{t,x})|_{\mathbb{T}_k}\big{)}\nonumber\\[1mm]
&&\hspace{-1.75em}=2\lambda \mathbb{E}_k\Big{\{} \sum_{\ell=k}^{N-1}\Big{[}(X^{k,\Phi}_\ell)^TQ_{k,\ell}\alpha^{k,\bar{u}_k}_\ell+(\mathbb{E}_kX^{k,\Phi}_\ell)^T \bar{Q}_{k,\ell}\mathbb{E}_k\alpha^{k,\bar{u}_k}_\ell
+(\Phi_\ell X^{k,\Phi}_\ell+v^{t,x}_\ell)^TR_{k,\ell}\Phi_\ell \alpha^{k,\bar{u}_k}_\ell  \nonumber\\
&&\hspace{-1.75em}\hphantom{=}+\mathbb{E}_k(\Phi_\ell X^{k,\Phi}_\ell+v^{t,x}_\ell)^T\bar{R}_{k,\ell}\Phi_\ell \mathbb{E}_k\alpha^{k,\bar{u}_k}_\ell+q_{k,\ell}^T\alpha^{k,\bar{u}_k}_\ell+\rho^T_{k,\ell}\Phi_\ell \alpha^{k,\bar{u}_k}_\ell
\Big{]}+\big{[}\mathcal{R}_{k,k}(\Phi_kX^{k,\Phi}_k+v^{t,x}_k)+\rho_{k,k}
\big{]}^T\bar{u}_k\nonumber\\[1mm]
&&\hspace{-1.75em}\hphantom{=}+\big{[}G_kX_N^{k,\Phi}+F_kX^{t,x,*}_k+g_k \big{]}^T \alpha^{k,\bar{u}_k}_N+(\mathbb{E}_kX^{k,\Phi}_{N})^T\bar{G}_{k}\mathbb{E}_k\alpha^{k,\bar{u}_k}_N
\Big{\}}\nonumber\\[1mm]
&&\hspace{-1.75em}\hphantom{=}+\lambda^2\Big{\{}\sum_{\ell=k}^{N-1}\mathbb{E}_k\Big{[}(\alpha^{k,\bar{u}_k}_\ell)^T Q_{k,\ell}\alpha^{k,\bar{u}_k}_\ell+ (\mathbb{E}_k\alpha^{k,\bar{u}_k}_\ell)^T\bar{Q}_{k,\ell}\mathbb{E}_k\alpha^{k,\bar{u}_k}_\ell+(\alpha^{k,\bar{u}_k}_\ell)^T \Phi_\ell^TR_{k,\ell}\Phi_\ell \alpha^{k,\bar{u}_k}_\ell\nonumber\\
&&\hspace{-1.75em}\hphantom{=}+(\mathbb{E}_k\alpha^{k,\bar{u}_k}_\ell)^T\Phi_\ell^T\bar{R}_{k,\ell}\Phi_\ell \mathbb{E}_k\alpha^{k,\bar{u}_k}_\ell  \Big{]}+\mathbb{E}_k\big{[}\bar{u}_k^T\mathcal{R}_{k,k}\bar{u}_k
\big{]}+\mathbb{E}_k\big{[}(\alpha^{k,\bar{u}_k}_N)^TG_k \alpha_N^{k,\bar{u}_k} \big{]}+(\mathbb{E}_k\alpha^{k,\bar{u}_k}_N)^T\bar{G}_k \mathbb{E}_k\alpha_N^{k,\bar{u}_k}  \Big{\}}\nonumber\\
%
%
%
%
%
%
%
%
%
%
%
%
&&\hspace{-1.75em}=2\lambda \mathbb{E}_k\Big{\{} \sum_{\ell=k}^{N-1}\Big{[}\Big{(}Q_{k,\ell}X^{k,\Phi}_\ell+\Phi_\ell^T R_{k,\ell}(\Phi_\ell X^{k,\Phi}_\ell+v^{t,x}_\ell)+q_{k,\ell}+\Phi_\ell^T\rho_{k,\ell}\Big{)}^T\alpha^{k,\bar{u}_k}_\ell\nonumber
 \\
&&\hspace{-1.75em}\hphantom{=}+\Big{(}\bar{Q}_{k,\ell}\mathbb{E}_kX^{k,\Phi}_\ell+\Phi_\ell^T\bar{R}_{k,\ell}(\Phi_\ell \mathbb{E}_k X^{k,\Phi}_\ell+\mathbb{E}_kv^{t,x}_\ell)\Big{)}^T\mathbb{E}_k\alpha^{k,\bar{u}_k}_\ell
\Big{]}+\big{[}\mathcal{R}_{k,k}(\Phi_kX^{k,\Phi}_k+v^{t,x}_k)+\rho_{k,k}
\big{]}^T\bar{u}_k\nonumber\\[1mm]
&&\hspace{-1.75em}\hphantom{=}+\big{[}G_kX_N^{k,\Phi}+F_kX^{t,x,*}_k+g_k \big{]}^T \alpha^{k,\bar{u}_k}_N+(\mathbb{E}_kX^{k,\Phi}_{N})^T\bar{G}_{k}\mathbb{E}_k\alpha^{k,\bar{u}_k}_N
\Big{\}}+\lambda^2\widetilde{J}(k,0;\bar{u}_k).
\end{eqnarray}
From (\ref{system-y-k}) and (\ref{Sec3:bsde}), it follows that
\begin{eqnarray*}
&&\hspace{-1.5em} \mathbb{E}_k\Big{\{} \sum_{\ell=k}^{N-1}\Big{[}\Big{(}Q_{k,\ell}X^{k,\Phi}_\ell+(\Phi_\ell )^TR_{k,\ell}(\Phi_\ell X^{k,\Phi}_\ell+v^{t,x}_\ell)+q_{k,\ell}+\Phi_\ell^T\rho_{k,\ell}\Big{)}^T\alpha^{k,\bar{u}_k}_\ell
 \\
&&\hspace{-1.5em} \hphantom{=}+\Big{(}\bar{Q}_{k,\ell}\mathbb{E}_kX^{k,\Phi}_\ell+(\Phi_\ell )^T\bar{R}_{k,\ell}(\Phi_\ell \mathbb{E}_k X^{k,\Phi}_\ell+\mathbb{E}_kv^{t,x}_\ell)\Big{)}^T\mathbb{E}_k\alpha^{k,\bar{u}_k}_\ell
\Big{]}\\[1mm]
&&\hspace{-1.5em} \hphantom{=}+\big{[}\mathcal{R}_{k,k}(\Phi_kX^{k,\Phi}_k+v^{t,x}_k)+\rho_{k,k}
\big{]}^T\bar{u}_k+\big{[}G_kX_N^{k,\Phi}+F_kX^{t,x,*}_k+g_k \big{]}^T \alpha^{k,\bar{u}_k}_N+(\mathbb{E}_kX^{k,\Phi}_{N})^T\bar{G}_{k}\mathbb{E}_k\alpha^{k,\bar{u}_k}_N\Big{\}}\\
%
%
%
%
%
&&\hspace{-1.5em} =\sum_{\ell=k}^{N-1}\mathbb{E}_k\Big{\{}\Big{[}\big{(}Q_{k,\ell}+\Phi_\ell^TR_{k,\ell}\Phi_\ell\big{)} (X^{k,\Phi}_\ell-\mathbb{E}_kX^{k,\Phi}_\ell)+\Phi_\ell^TR_{k,\ell}(v^{t,x}_\ell-\mathbb{E}_kv^{t,x}_\ell)\\
&&\hspace{-1.5em} \hphantom{=}+(A_{k,\ell}+B_{k,\ell}\Phi_\ell)^T(\mathbb{E}_\ell Y_{\ell+1}^{k,\Phi}-\mathbb{E}_k Y_{\ell+1}^{k,\Phi})
+\sum_{i=1}^p(C^i_{k,\ell}+D^i_{k,\ell}\Phi_{\ell})^T\big{(}\mathbb{E}_\ell(Y_{\ell+1}^{k,\Phi}w^i_\ell) -\mathbb{E}_k(Y_{\ell+1}^{k,\Phi}w^i_\ell)\big{)}\\
&&\hspace{-1.5em} \hphantom{=}-(Y^{k,\Phi}_\ell-\mathbb{E}_kY^{k,\Phi}_\ell)\Big{]}^T(\alpha^{k,\bar{u}_k}_\ell-\mathbb{E}_k\alpha^{k,\bar{u}_k}_\ell )
+\Big{[}\big{(}{\mathcal{Q}}_{k,\ell}+\Phi_\ell^T{\mathcal{R}}_{k,\ell}\Phi_\ell\big{)} \mathbb{E}_kX^{k,\Phi}_\ell+\Phi_\ell^T\mathcal{R}_{k,\ell}\mathbb{E}_kv^{t,x}_\ell+q_{k,\ell}+\Phi_\ell^T\rho_{k,\ell}\\[1mm]
&&\hspace{-1.5em} \hphantom{=}+\sum_{i=1}^p(\mathcal{C}^i_{k,\ell}+\mathcal{D}^i_{k,\ell}\Phi_\ell)^T\mathbb{E}_{k}(Y^{k,\Phi}_{\ell+1}w^i_\ell) +(\mathcal{A}_{k,\ell}+\mathcal{B}_{k,\ell}\Phi_\ell)^T\mathbb{E}_kY^{k,\Phi}_{\ell+1}-\mathbb{E}_kY^{k,\Phi}_\ell\Big{]}^T \mathbb{E}_k\alpha^{k,\bar{u}_k}_\ell
\Big{\}}\\[1mm]
&&\hspace{-1.5em} \hphantom{=}+\big{[}\mathcal{R}_{k,k}(\Phi_kX^{k,\Phi}_k+v^{t,x}_k)+\mathcal{B}^T_{k,k}\mathbb{E}_kY_{k+1}^{k,\Phi} +\sum_{i=1}^p(\mathcal{D}^i_{k,k})^T\mathbb{E}_k(Y_{k+1}^{k,\Phi}w^i_{k})+\rho_{k,k}
\big{]}^T\bar{u}_k\\[1mm]
&&\hspace{-1.5em} =\big{[}\mathcal{R}_{k,k}(\Phi_kX^{k,\Phi}_k+v^{t,x}_k)+\mathcal{B}^T_{k,k}\mathbb{E}_kY_{k+1}^{k,\Phi} +\sum_{i=1}^p(\mathcal{D}^i_{k,k})^T\mathbb{E}_k(Y_{k+1}^{k,\Phi}w^i_{k})+\rho_{k,k}
\big{]}^T\bar{u}_k.
\end{eqnarray*}
From (\ref{Lemma-difference:J}), we can complete the proof. \hfill $\square$

\section{Proof of Lemma \ref{Sec3:Lemma-Z}}\label{appendix-Sec2:lemma-Z}

By simple calculation, we have
\begin{eqnarray*}
&&\hspace{-2em}\mathbb{E}_{N-1}Y^{k,\Phi}_{N}=\mathbb{E}_{N-1}\big{[}G_kX^{k,\Phi}_N+\bar{G}_k\mathbb{E}_kX^{k,\Phi}_N+F_kX^{t,x,*}_k+g_k\big{]}\nonumber\\[1mm]
&&\hspace{-2em}\hphantom{\mathbb{E}_{N-1}Y^{k,\Phi}_{N}}=G_k(A_{k,N-1}+B_{k,N-1}\Phi_{N-1})X^{k,\Phi}_{N-1}+\big{[}G_k(\bar{A}_{k,N-1}+\bar{B}_{k,N-1}\Phi_{N-1}) \nonumber\\[1mm]
&&\hspace{-2em}\hphantom{\mathbb{E}_{N-1}Y^{k,\Phi}_{N}=}+\bar{G}_k({\mathcal{A}}_{k,N-1}+{\mathcal{B}}_{k,N-1}\Phi_{N-1})\big{]} \mathbb{E}_kX^{k,\Phi}_{N-1}+G_kB_{k,N-1}\Gamma_{N-1}X^{t,x,*}_{N-1}\nonumber\\[1mm]
&&\hspace{-2em}\label{Sec3:Lemma-Z--1}\hphantom{\mathbb{E}_{N-1}Y^{k,\Phi}_{N}=}+\big{[}G_k\bar{B}_{k,N-1} +\bar{G}_k\mathcal{B}_{k,N-1}\big{]}\Gamma_{N-1}\mathbb{E}_kX^{t,x,*}_{N-1} +{\mathcal{G}}_k\mathcal{B}_{k,N-1}\bar{v}^{t,x}_{N-1}+\mathcal{G}_kf_{k,N-1}+F_kX^{t,x,*}_k+g_k,\\[1mm]
&&\hspace{-2em}\mathbb{E}_{k}Y^{k,\Phi}_{N}=\mathcal{G}_k(\mathcal{A}_{k,N-1}+\mathcal{B}_{k,N-1}\Phi_{N-1})\mathbb{E}_kX^{k,\Phi}_{N-1} +\mathcal{G}_k\mathcal{B}_{k,N-1}\big{(}\Gamma_{N-1}\mathbb{E}_kX^{t,x,*}_{N-1}+\bar{v}^{t,x}_{N-1}\big{)}\nonumber\\[1mm]
&&\hspace{-2em}\hphantom{\mathbb{E}_{k}Y^{k,\Phi}_{N}=}+\mathcal{G}_kf_{k,N-1}+F_kX^{t,x,*}_k+g_k,\nonumber\\[1mm]
&&\hspace{-2em}\mathbb{E}_{N-1}[Y^{k,\Phi}_{N}w^i_{N-1}]=\mathbb{E}_{N-1}\big{[}\big{(}G_kX^{k,\Phi}_N +\bar{G}_k\mathbb{E}_kX^{k,\Phi}_N+F_kX^{t,x,*}_{k}+g_k\big{)}w^i_{N-1}\big{]}\nonumber\\[1mm]
&&\hspace{-2em}\hphantom{\mathbb{E}_{N-1}Y^{k,\Phi}_{N}}=G_k\sum_{j=1}^p\delta_{N-1}^{ij}\Big{[}(C^j_{k,N-1}+D^j_{k,N-1}\Phi_{N-1})X^{k,\Phi}_{N-1} +(\bar{C}^j_{k,N-1}+\bar{D}^j_{k,N-1}\Phi_{N-1})\mathbb{E}_{k}X^{k,\Phi}_{N-1}\nonumber\\
&&\hspace{-2em}\label{Sec3:Lemma-Z--2}\hphantom{\mathbb{E}_{N-1}Y^{k,\Phi}_{N}=}+D^j_{k,N-1}\Gamma_{N-1}X^{t,x,*}_{N-1}+\bar{D}^j_{k,N-1}\Gamma_{N-1} \mathbb{E}_{k}X^{t,x,*}_{N-1}+\mathcal{D}^j_{k,N-1}\bar{v}^{t,x}_{N-1}+d^j_{k,N-1}\Big{]},\\[1mm]
&&\hspace{-2em}\mathbb{E}_{k}[Y^{k,\Phi}_{N}w^i_{N-1}]=G_k\sum_{j=1}^p\delta_{N-1}^{ij}\Big{[}(\mathcal{C}^j_{k,N-1} +\mathcal{D}^j_{k,N-1}\Phi_{N-1})\mathbb{E}_kX^{k,\Phi}_{N-1} +\mathcal{D}^j_{k,N-1}\big{(}\Gamma_{N-1}\mathbb{E}_{k}X^{t,x,*}_{N-1}
+\bar{v}^{t,x}_{N-1}\big{)}+d^j_{k,N-1}\Big{]}\nonumber.
\end{eqnarray*}
From (\ref{system-adjoint-closed}), we have
\begin{eqnarray*}
&&\hspace{-1em}Y_{N-1}^{k,\Phi}=Q_{k,N-1}X^{k,\Phi}_{N-1}+\bar{Q}_{k,N-1}\mathbb{E}_kX^{k,\Phi}_{N-1}+\Phi_{N-1}^TR_{k,{N-1}}\Phi_{N-1} X^{k,\Phi}_{N-1}
\\[1mm]
&&\hspace{-1em}\hphantom{Z_{{N-1}}^{k,\Phi}=}
+\Phi_{N-1}^T\bar{R}_{k,{N-1}}\Phi_{N-1} \mathbb{E}_kX^{k,\Phi}_{N-1}+\Phi_{N-1}^TR_{k,{N-1}}\big{(}\Gamma_{N-1}X^{t,x,*}_{N-1}+\bar{v}^{t,x}_{N-1}\big{)}\\[1mm]
&&\hspace{-1em}\hphantom{Z_{{N-1}}^{k,\Phi}=}+\Phi_{N-1}^T{\bar{R}}_{k,{N-1}}\big{(}\Gamma_{N-1}\mathbb{E}_kX^{t,x,*}_{N-1} +\bar{v}^{t,x}_{N-1}\big{)}+\big{(}A_{k,{N-1}}+B_{k,{N-1}}\Phi_{N-1}\big{)}^T \mathbb{E}_{N-1} Y_{N}^{k,\Phi}\\[1mm]
&&\hspace{-1em}\hphantom{Z_{{N-1}}^{k,\Phi}=}+\big{(}\bar{A}_{k,{N-1}}+\bar{B}_{k,{N-1}}\Phi_{N-1}\big{)}^T\mathbb{E}_kY_{N}^{k,\Phi} +\sum_{i=1}^p\big{[}\big{(}C^i_{k,{N-1}}+D^i_{k,{N-1}}\Phi_{N-1}\big{)}^T\mathbb{E}_{N-1}\big{(} Y_{N}^{k,\Phi}w^i_{N-1}\big{)}\\[1mm]
&&\hspace{-1em}\hphantom{Z_{{N-1}}^{k,\Phi}=}+\big{(}{\bar{C}}^i_{k,{N-1}}+{\bar{D}}^i_{k,{N-1}}\Phi_{N-1}\big{)}^T \mathbb{E}_k\big{(}Y^{k,\Phi}_{N}w^i_{N-1}\big{)}\big{]} +\Phi_{{N-1}}^T\rho_{k,{N-1}} +q_{k,{N-1}}\\[1mm]
&&\hspace{-1em}\hphantom{Z_{N-1}^{k,\Phi}}=\Big{\{}Q_{k,N-1}+\Phi_{N-1}^TR_{k,N-1}\Phi_{N-1}+\big{(}A_{k,N-1}+B_{k,N-1}\Phi_{N-1}\big{)}^TG_{k} \big{(}A_{k,N-1}+B_{k,N-1}\Phi_{N-1}\big{)}\\[1mm]
&&\hspace{-1em}\hphantom{Z_{N-1}^{k,\Phi}=}+\sum_{i,j=1}^p\delta_{N-1}^{ij}\big{(}C^i_{k,N-1} +D^i_{k,N-1}\Phi_{N-1}\big{)}^TG_{k}\big{(}C^j_{k,N-1}+D^j_{k,N-1}\Phi_{N-1}\big{)}  \Big{\}}X^{k,\Phi}_{N-1}\\[1mm]
&&\hspace{-1em}\hphantom{Z_{N-1}^{k,\Phi}=}+\Big{\{}\bar{Q}_{k,N-1}+\Phi_{N-1}^T\bar{R}_{k,N-1}\Phi_{N-1}+\big{(}A_{k,N-1}+B_{k,N-1}\Phi_{N-1}\big{)}^T\big{[}G_k(\bar{A}_{k,N-1}+\bar{B}_{k,N-1}\Phi_{N-1}) \\[1mm]
&&\hspace{-1em}\hphantom{Z_{N-1}^{k,\Phi}=}+\bar{G}_k({\mathcal{A}}_{k,N-1}+{\mathcal{B}}_{k,N-1}\Phi_{N-1}) \big{]}+\big{(}\bar{A}_{k,N-1}+\bar{B}_{k,N-1}\Phi_{N-1}\big{)}^T\mathcal{G}_k(\mathcal{A}_{k,N-1}+\mathcal{B}_{k,N-1}\Phi_{N-1})\\[1mm]
&&\hspace{-1em}\hphantom{Z_{N-1}^{k,\Phi}=}+\sum_{i,j=1}^j\delta_{N-1}^{ij}\big{[}\big{(}C^i_{k,N-1}+D^i_{k,N-1} \Phi_{N-1}\big{)}^TG_k(\bar{C}^j_{k,N-1}+\bar{D}^j_{k,N-1}\Phi_{N-1})\\[1mm]
&&\hspace{-1em}\hphantom{Z_{N-1}^{k,\Phi}=}+\big{(}{\bar{C}}^i_{k,N-1}+{\bar{D}}^i_{k,N-1}\Phi_{N-1}\big{)}^TG_k(\mathcal{C}^j_{k,N-1} +\mathcal{D}^j_{k,N-1}\Phi_{N-1})\big{]}\Big{\}}\mathbb{E}_kX^{k,\Phi}_{N-1}\\[1mm]
&&\hspace{-1em}\hphantom{Z_{N-1}^{k,\Phi}=}+\Big{\{}\Phi_{N-1}^TR_{k,N-1}+\big{(}A_{k,N-1}+B_{k,N-1}\Phi_{N-1}\big{)}^TG_kB_{k,N-1}\\[1mm]
&&\hspace{-1em}\hphantom{Z_{N-1}^{k,\Phi}=}+\sum_{i,j=1}^p\delta_{N-1}^{ij}\big{(}C^i_{k,N-1}+D^i_{k,N-1}\Phi_{N-1}\big{)}^TG_kD^j_{k,N-1}\Big{\}} \big{(}\Gamma_{N-1}X^{t,x,*}_{N-1}+\bar{v}^{t,x}_{N-1}\big{)}\\[1mm]
&&\hspace{-1em}\hphantom{Z_{N-1}^{k,\Phi}=}+\Big{\{}\Phi_{N-1}^T{\bar{R}}_{k,N-1}+\big{(}A_{k,N-1}+B_{k,N-1}\Phi_{N-1}\big{)}^T \big{(}G_k\bar{B}_{k,N-1}+\bar{G}_k\mathcal{B}_{k,N-1}\big{)}\\[1mm]
&&\hspace{-1em}\hphantom{Z_{N-1}^{k,\Phi}=}+\big{(}\bar{A}_{k,N-1}+\bar{B}_{k,N-1}\Phi_{N-1}\big{)}^T\mathcal{G}_k\mathcal{B}_{k,N-1} +\sum_{i,j=1}^p\delta_{N-1}^{ij}\big{[}\big{(}C^i_{k,N-1}+D^i_{k,N-1}\Phi_{N-1}\big{)}^TG_k\bar{D}^j_{k,N-1}\\[1mm]
&&\hspace{-1em}\hphantom{Z_{N-1}^{k,\Phi}=}+\big{(}{\bar{C}}^i_{k,N-1}+{\bar{D}}^i_{k,N-1}\Phi_{N-1}\big{)}^TG_k\mathcal{D}^j_{k,N-1} \big{]} \Big{\}}\big{(}\Gamma_{N-1}\mathbb{E}_kX^{t,x,*}_{N-1}+\bar{v}^{t,x}_{N-1}\big{)}\\[1mm]
&&\hspace{-1em}\hphantom{Z_{N-1}^{k,\Phi}=}+\big{(}A_{k,N-1}+B_{k,N-1}\Phi_{N-1}\big{)}^T\big{(}\mathcal{G}_kf_{k,N-1}+F_kX^{t,x,*}_k+g_k\big{)}\\[1mm]
&&\hspace{-1em}\hphantom{Z_{N-1}^{k,\Phi}=}+\big{(}\bar{A}_{k,N-1}+\bar{B}_{k,N-1}\Phi_{N-1}\big{)}^T\big{(}\mathcal{G}_kf_{k,N-1}+F_kX^{t,x,*}_k+g_k\big{)}\\[1mm]
&&\hspace{-1em}\hphantom{Z_{N-1}^{k,\Phi}=} +\sum_{i,j=1}^p\delta_{N-1}^{ij}\big{[}\big{(}C^i_{k,N-1}+D^i_{k,N-1}\Phi_{N-1}\big{)}^TG_kd^j_{k,N-1}+\big{(}\bar{C}^i_{k,N-1}+\bar{D}^i_{k,N-1}\Phi_{N-1}\big{)}^TG_kd^j_{k,N-1}\big{]}\\[1mm]
&&\hspace{-1em}\hphantom{Z_{N-1}^{k,\Phi}=}+\Phi_{N-1}^T\rho_{k,N-1} +q_{k,N-1}\\[1mm]
&&\hspace{-1em}\hphantom{Z_{N-1}^{k,\Phi}}= S_{k,N-1}X^{k,\Phi}_{N-1}+\bar{S}_{k,N-1}\mathbb{E}_{k}X^{k,\Phi}_{N-1}+T_{k,N-1}X^{t,x,*}_{N-1}+\bar{T}_{N-1}\mathbb{E}_{k}X^{t,x,x*}_{N-1}+U_{k,N-1}X^{t,x,*}_{k}+\pi_{k,N-1}.
\end{eqnarray*}
In the above, we apply the property $\bar{v}^{t,x}\in l^2(\mathbb{T}_t; \mathbb{R}^{m})$. By deduction, we can get the desired result. \hfill $\square$

\section{Proof of Theorem \ref{Sec3:Theorem-Necessary}}\label{appendix:theorem-necessary}

{\emph{i)$\Rightarrow$ii)}}. Let $(\Phi, v^{t,x})$ be a mixed equilibrium solution of Problem (LQ)$_{tx}$, which satisfies (\ref{stationary-condition-closed}) and (\ref{convex-closed}).
By simple calculation, we have
\begin{eqnarray}\label{convex-closed-2}
&&\widetilde{J}(k,0;\bar{u}_k)=\sum_{\ell=k}^{N-1}\mathbb{E}_k\Big{[}
(\alpha^{k,\bar{u}_k}_{\ell+1}-\mathbb{E}_k\alpha^{k,\bar{u}_k}_{\ell+1})^TS_{k,\ell+1} (\alpha^{k,\bar{u}_k}_{\ell+1}-\mathbb{E}_k\alpha^{k,\bar{u}_k}_{\ell+1})\nonumber \\[1mm]
&&\hphantom{\widetilde{J}(k,0;\bar{u}_k)=\Big{\{}}-(\alpha^{k,\bar{u}_k}_{\ell}-\mathbb{E}_k\alpha^{k,\bar{u}_k}_{\ell})^TS_{k,\ell} (\alpha^{k,\bar{u}_k}_\ell-\mathbb{E}_k\alpha^{k,\bar{u}_k}_\ell)\nonumber \\[1mm]
&&\hphantom{\widetilde{J}(k,0;\bar{u}_k)=\Big{\{}}+(\alpha^{k,\bar{u}_k}_\ell-\mathbb{E}_k\alpha^{k,\bar{u}_k}_\ell)^T \big{(}Q_{k,\ell}+\Phi_\ell^TR_{k,\ell}\Phi_\ell \big{)}(\alpha^{k,\bar{u}_k}_\ell-\mathbb{E}_k\alpha^{k,\bar{u}_k}_\ell)\nonumber \\[1mm]
&&\hphantom{\widetilde{J}(k,0;\bar{u}_k)=\Big{\{}}+(\mathbb{E}_k\alpha^{k,\bar{u}_k}_{\ell+1})^T\mathcal{S}_{k,\ell+1} \mathbb{E}_k\alpha^{k,\bar{u}_k}_{\ell+1} -(\mathbb{E}_k\alpha^{k,\bar{u}_k}_\ell)^T\mathcal{S}_{k,\ell}\mathbb{E}_k\alpha^{k,\bar{u}_k}_\ell\nonumber
\\[1mm]
&&\hphantom{\widetilde{J}(k,0;\bar{u}_k)=\Big{\{}}+(\mathbb{E}_k\alpha^{k,\bar{u}_k}_\ell)^T\big{(}{\mathcal{Q}}_{k,\ell} +\Phi_\ell^T{\mathcal{R}}_{k,\ell}\Phi_\ell\big{)} \mathbb{E}_k\alpha^{k,\bar{u}_k}_\ell  \Big{]}+\mathbb{E}_k\big{[}\bar{u}_k^T\mathcal{R}_{k,k}\bar{u}_k
\big{]}
\nonumber \\
&&\hphantom{\widetilde{J}(k,0;\bar{u}_k)}=\mathbb{E}_k\big{[}\bar{u}_k^T{\mathbb{O}}_k\bar{u}_k\big{]}=\bar{u}_k^T{\mathbb{O}}_k\bar{u}_k.
\end{eqnarray}
From (\ref{convex-closed}) and (\ref{convex-closed-2}), it holds that
\begin{eqnarray*}
%
&&\inf_{\bar{u}_k\in L^2_\mathcal{F}(k; \mathbb{R}^m)}\widetilde{J}(k,0;\bar{u}_k)= \inf_{\bar{u}_k\in L^2_\mathcal{F}(k; \mathbb{R}^m)}\big{[}\bar{u}_k^T{\mathbb{O}}_k\bar{u}_k\big{]}\geq 0,
\end{eqnarray*}
which implies ${\mathbb{O}}_k\succeq 0$. Then, (\ref{Sec3:Theorem-Necessary-S}) is solvable.

We now prove b) and c).
Letting $k=N-1$ in (\ref{stationary-condition-closed}) and noting
\begin{eqnarray}
&&\hspace{-1em}\mathbb{E}_{N-1}Y^{N-1,\Phi}_{N}=\mathcal{G}_{N-1}(\mathcal{A}_{N-1,N-1}+\mathcal{B}_{N-1,N-1}\Phi_{N-1})X^{N-1,\Phi}_{N-1}+\mathcal{G}_{N-1}\mathcal{B}_{N-1,N-1}v^{t,x}_{N-1} \nonumber \\[1mm]
&&\hspace{-1em}\label{Z-(N-1)}\hphantom{\mathbb{E}_{N-1}Y^{N-1,\Phi}_{N}=} +\mathcal{G}_{N-1}f_{N-1,N-1}+F_{N-1}X^{t,x,*}_{N-1}+g_{N-1},\\
&&\hspace{-1em}\mathbb{E}_{N-1}(Y^{N-1,\Phi}_{N}w^i_{N-1})=G_{N-1}\sum_{j=1}^p\delta_k^{ij}\Big{[}(\mathcal{C}^j_{N-1,N-1}+\mathcal{D}^j_{N-1,N-1}\Phi_{N-1}) X^{N-1,\Phi}_{N-1}\nonumber \\
&&\hspace{-1em}\hphantom{\mathbb{E}_{N-1}(Y^{N-1,\Phi}_{N}w^i_{N-1})=}+\mathcal{D}^j_{N-1,N-1}v^{t,x}_{N-1}+d^j_{N-1,N-1}\Big{]},\nonumber
\end{eqnarray}
we have
\begin{eqnarray}\label{stationary-condition-closed-(N-1)}
&&\hspace{-1em}0=\mathcal{R}_{N-1,N-1} (\Phi_{N-1}X^{t,x,*}_{N-1}+v^{t,x}_{N-1})+\mathcal{B}^T_{N-1,N-1}\mathbb{E}_{N-1}Y_{N}^{N-1,\Phi}\nonumber \\
&&\hspace{-1em}\hphantom{0=}+\sum_{i=1}^p(\mathcal{D}^i_{N-1,N-1})^T\mathbb{E}_{N-1}\big{(}Y_{N}^{N-1,\Phi}w_{N-1}^i\big{)} +\rho_{N-1,N-1}\nonumber \\[1mm]
&&\hspace{-1em}\hphantom{0}=\big{[}\mathcal{R}_{N-1,N-1}+\mathcal{B}^T_{N-1,N-1}\mathcal{G}_{N-1}\mathcal{B}_{N-1,N-1}+\sum_{i,j=1}^p\delta_{N-1}^{ij}(\mathcal{D}^i_{N-1,N-1})^T {G}_{N-1}\mathcal{D}^j_{N-1,N-1} \big{]}\nonumber \\[1mm]
&&\hspace{-1em}\hphantom{0=}\times \big{(}\Phi_{N-1}X^{t,x,*}_{N-1}+v^{t,x}_{N-1}\big{)}+\Big{[}\mathcal{B}^T_{N-1,N-1}\mathcal{G}_{N-1}\mathcal{A}_{N-1,N-1}+ \sum_{i,j=1}^p\delta_{N-1}^{ij}(\mathcal{D}^i_{N-1,N-1})^TG_{N-1}\mathcal{C}^j_{N-1,N-1}\nonumber \\
&&\hspace{-1em}\hphantom{0=}+\mathcal{B}_{N-1,N-1}^TF_{N-1}\Big{]}X^{t,x,*}_{N-1}+\mathcal{B}^T_{N-1,N-1}\mathcal{G}_{N-1}f_{N-1,N-1}+\sum_{i,j=1}^p\delta_{N-1}^{ij} (\mathcal{D}^i_{N-1,N-1})^TG_{N-1}d^j_{N-1,N-1} \nonumber \\[1mm]
&&\hspace{-1em}\hphantom{0=}+\mathcal{B}^T_{N-1,N-1}g_{N-1}+\rho_{N-1,N-1}\nonumber \\[1mm]
&&\hspace{-1em}\hphantom{0}=\mathcal{O}_{N-1}\big{(}\Phi_{N-1}X^{t,x,*}_{N-1} +v^{t,x}_{N-1}\big{)}+\mathcal{L}_{N-1}X^{t,x,*}_{N-1}+\theta_{N-1}.
\end{eqnarray}
Here, $X^{N-1,\Phi}_{N-1} =X^{t,x,*}_{N-1}$ and
\begin{eqnarray*}
\left\{
\begin{array}{l}
\mathcal{O}_{N-1}=\mathcal{R}_{N-1,N-1}+\mathcal{B}^T_{N-1,N-1}\mathcal{G}_{N-1}\mathcal{B}_{N-1,N-1}+\sum_{i,j=1}^p\delta_{N-1}^{ij}(\mathcal{D}^i_{N-1,N-1})^T {G}_{N-1}\mathcal{D}^j_{N-1,N-1},\\[1mm]
\mathcal{L}_{N-1}=\mathcal{B}^T_{N-1,N-1}\mathcal{G}_{N-1}\mathcal{A}_{N-1,N-1}+\sum_{i,j=1}^p\delta_{N-1}^{ij}(\mathcal{D}^i_{N-1,N-1})^TG_{N-1}\mathcal{C}^j_{N-1,N-1}+\mathcal{B}^T_{N-1,N-1}F_{N-1},\\[1mm]
\theta_{N-1}=\mathcal{B}^T_{N-1,N-1}\mathcal{G}_{N-1}f_{N-1,N-1}+\sum_{i,j=1}^p\delta_{N-1}^{ij}(\mathcal{D}^i_{N-1,N-1})^TG_{N-1}d_{N-1,N-1} \\[1mm]
\hphantom{\theta_{N-1}=}+\mathcal{B}^T_{N-1,N-1}g_{N-1}+\rho_{N-1,N-1}.
\end{array}
\right.
\end{eqnarray*}
Note that $(\Phi, v^{t,x})$ is a mixed equilibrium solution and $X^{t,x,*}$ is given in (\ref{equili-state-closed}).
As $\Phi_{N-1}X^{t,x,*}_{N-1} +v^{t,x}_{N-1}$ satisfies (\ref{stationary-condition-closed-(N-1)}), it holds from Lemma \ref{Lemma-matrix-equation} that (\ref{stationary-condition-closed-(N-1)}) is equivalent to
\begin{eqnarray}\label{Ran-(N-1)}
\mathcal{L}_{N-1}X^{t,x,*}_{N-1}+\theta_{N-1}\in \mbox{Ran}\big{(}\mathcal{O}_{N-1}\big{)},
\end{eqnarray}
and for some $\eta_{N-1}\in \mathbb{R}^m$,
\begin{eqnarray}\label{Ran-(N-1)-2}
\Phi_{N-1}X^{t,x,*}_{N-1} +v^{t,x}_{N-1}=-\mathcal{O}^{\dagger}_{N-1} \mathcal{L}_{N-1}X^{t,x,*}_{N-1}-\mathcal{O}^{\dagger}_{N-1}\theta_{N-1}+\big{(}I
-\mathcal{O}^{\dagger}_{N-1}\mathcal{O}_{N-1}\big{)}\eta_{N-1}.
\end{eqnarray}
Clearly, (\ref{Ran-(N-1)-2}) is equivalent to
\begin{eqnarray}\label{v-open-(N-1)}
v^{t,x}_{N-1}=-(\mathcal{O}^{\dagger}_{N-1} \mathcal{L}_{N-1}+\Phi_{N-1})X^{t,x,*}_{N-1}-\mathcal{O}^{\dagger}_{N-1} \theta_{N-1}+\big{(}I-\mathcal{O}^{\dagger}_{N-1}\mathcal{O}_{N-1}\big{)}\eta_{N-1}
\end{eqnarray}
for some $\eta_{N-1}\in \mathbb{R}^m$. If we replace $v^{t,x}_{N-1}$ of (\ref{v-open-(N-1)}) by
\begin{eqnarray}\label{v-open-(N-1)-2}
v^{t,x}_{N-1}=-(\mathcal{O}^{\dagger}_{N-1} \mathcal{L}_{N-1}+\Phi_{N-1})X^{t,x,*}_{N-1}-\mathcal{O}^{\dagger}_{N-1} \theta_{N-1},
\end{eqnarray}
then the new pair $(\Phi, v^{t,x})$ with $v_{N-1}^{t,x}$ given in (\ref{v-open-(N-1)-2}) can also serve as a mixed equilibrium solution. By submitting the pair $(\Phi_{N-1}, v^{t,x}_{N-1})$ (with $v_{N-1}^{t,x}$ given in (\ref{v-open-(N-1)-2})), the equations in (\ref{stationary-condition-closed-(N-1)}) are also satisfied. Therefore, the $v^{t,x}_{N-1}$ (in (\ref{stationary-condition-closed-(N-1)})) is selected as
\begin{eqnarray}\label{v-open-(N-1)-3}
v^{t,x}_{N-1}=-(\mathcal{O}^{\dagger}_{N-1} \mathcal{L}_{N-1}+\Phi_{N-1})X^{t,x,*}_{N-1}-\mathcal{O}^{\dagger}_{N-1} \theta_{N-1}.
\end{eqnarray}

Substituting this $v^{t,x}_{N-1}$ into Lemma \ref{Sec3:Lemma-Z}, we have
\begin{eqnarray*}
&&Y^{N-2,\Phi}_{N-1}=S_{N-2,N-1}X^{N-2,\Phi}_{N-1}+\bar{S}_{N-2,N-1}\mathbb{E}_{N-2}X^{N-2,\Phi}_{N-1}+T_{N-2,N-1}X^{t,x,*}_{N-1} \\
&&\hphantom{Y^{N-2,\Phi}_{N-1}=}+\bar{T}_{N-2,N-1}\mathbb{E}_{N-2}X^{t,x,*}_{N-1}+U_{N-2,N-1}X^{t,x,*}_{N-2}+\pi_{N-2,N-1}.
\end{eqnarray*}
In this case, it holds that
\begin{eqnarray*}
&&\mathbb{E}_{N-2}Y^{N-2,\Phi}_{N-1}=  \big{(}\mathcal{S}_{N-2,N-1} +\mathcal{T}_{N-2,N-1}\big{)}\Big{[}\mathcal{A}_{N-2,N-2}X^{t,x,*}_{N-2}+\mathcal{B}_{N-2,N-2}\big{(}\Phi_{N-2}X^{t,x,*}_{N-2}+v^{t,x}_{N-2}\big{)} \\[1mm]
&&\hphantom{\mathbb{E}_{N-2}Y^{N-2,\Phi}_{N-1}=}+f_{N-2,N-2}\Big{]}+U_{N-2,N-1}X^{t,x,*}_{N-2}+\pi_{N-2,N-1},
\end{eqnarray*}
and
\begin{eqnarray*}
&&\mathbb{E}_{N-2}\big{(}Y^{N-2,\Phi}_{N-1}w^i_{N-2}\big{)}=\big{(}{S}_{N-2,N-1}+T_{N-2,N-1}\big{)}\sum_{j=1}^p\delta_{N-2}^{ij}\Big{[}\mathcal{C}^j_{N-2,N-2}X^{t,x,*}_{N-2}\\[1mm]
&&\hphantom{\mathbb{E}_{N-2}(Y^{N-2,\Phi}_{N-1}w_{N-2})=}+\mathcal{D}^j_{N-2,N-2}\big{(} \Phi_{N-2}X^{t,x,*}_{N-2}+v^{t,x}_{N-2}\big{)}+d^j_{N-2,N-2}\Big{]}.
\end{eqnarray*}
Therefore, we have
\begin{eqnarray}\label{stationary-condition-closed-(N-2)}
&&0=\mathcal{R}_{N-2,N-2} (\Phi_{N-2}X^{t,x,*}_{N-2}+v^{t,x}_{N-2})+\mathcal{B}^T_{N-2,N-2}\mathbb{E}_{N-2}Y_{N-1}^{N-2,\Phi}\nonumber \\[1mm]
&&\hphantom{0=}+\mathcal{D}_{N-2,N-2}^T\mathbb{E}_{N-2}\big{(}Y_{N-1}^{N-2,\Phi}w_{N-2}\big{)} +\rho_{N-2,N-2}\nonumber\\[1mm]
&&\hphantom{0}=\Big{[}\mathcal{R}_{N-2,N-2}+\mathcal{B}^T_{N-2,N-2}\big{(}\mathcal{S}_{N-2,N-1}+\mathcal{T}_{N-2,N-1}\big{)}\mathcal{B}_{N-2,N-2}\nonumber\\[1mm]
&&\hphantom{0=}+\sum_{i,j=1}^p\delta_{N-2}^{ij}(\mathcal{D}^i_{N-2,N-2})^T\big{(}{S}_{N-2,N-1}+T_{N-2,N-1}\big{)}\mathcal{D}^j_{N-2,N-2} \Big{]}\big{(}\Phi_{N-2}X^{t,x,*}_{N-2}+v^{t,x}_{N-2}\big{)}\nonumber\\[1mm]
&&\hphantom{0=}+\Big{[}\mathcal{B}^T_{N-2,N-2}\big{(}\mathcal{S}_{N-2,N-1}+\mathcal{T}_{N-2,N-1}\big{)}\mathcal{A}_{N-2,N-2}\nonumber\\[1mm]
&&\hphantom{0=}+\sum_{i,j=1}^p\delta_{N-2}^{ij}(\mathcal{D}^i_{N-2,N-2})^T \big{(}S _{N-2,N-1}+T_{N-2,N-1}\big{)}\mathcal{C}^j_{N-2,N-2}+\mathcal{B}_{N-2,N-2}^TU_{N-2,N-1}\Big{]}X^{t,x,*}_{N-1}\nonumber \\[1mm]
&&\hphantom{0=}+\mathcal{B}^T_{N-2,N-2}\big{(}\mathcal{S}_{N-2,N-1}+\mathcal{T}_{N-2,N-1}\big{)}f_{N-2,N-2}\nonumber\\[1mm]
&&\hphantom{0=}+\sum_{i,j=1}^p\delta_{N-2}^{ij}(\mathcal{D}^i_{N-2,N-2})^T\big{(}S_{N-2,N-1}+T_{N-2,N-1}\big{)}d^j_{N-2,N-2}+\mathcal{B}^T_{N-2,N-2}\pi_{N-2,N-1}+\rho_{N-2,N-2}\nonumber\\[1mm]
&&\hphantom{0}=\mathcal{O}_{N-2}\big{(}\Phi_{N-2}X^{t,x,*}_{N-2} +v^{t,x}_{N-2}\big{)}+\mathcal{L}_{N-2}X^{t,x,*}_{N-2}+\theta_{N-2},
\end{eqnarray}
where $X^{N-2,\Phi}_{N-2} =X^{t,x,*}_{N-2}$ and
\begin{eqnarray*}
\left\{
\begin{array}{l}
\mathcal{O}_{N-2}=\mathcal{R}_{N-2,N-2}+\mathcal{B}^T_{N-2,N-2}\big{(}\mathcal{S}_{N-2,N-1}+\mathcal{T}_{N-2,N-1}\big{)}\mathcal{B}_{N-2,N-2}\\[1mm]
\hphantom{\mathcal{O}_{N-1}=}+\sum_{i,j=1}^p\delta_{N-2}^{ij}(\mathcal{D}^i_{N-2,N-2})^T\big{(}{S}_{N-2,N-1}+T_{N-2,N-1}\big{)}\mathcal{D}^j_{N-2,N-2},\\[1mm]
\mathcal{L}_{N-2}=\mathcal{B}^T_{N-2,N-2}\big{(}\mathcal{S}_{N-2,N-1}+\mathcal{T}_{N-2,N-1}\big{)}\mathcal{A}_{N-2,N-2}\\[1mm]
\hphantom{\mathcal{L}_{N-1}=}+\sum_{i,j=1}^p\delta_{N-2}^{ij}(\mathcal{D}^i_{N-2,N-2})^T \big{(}S _{N-2,N-1}+T_{N-2,N-1}\big{)}\mathcal{C}^j_{N-2,N-2}+\mathcal{B}_{N-2,N-2}^TU_{N-2,N-1},\\[1mm]
\theta_{N-2}=\mathcal{B}^T_{N-2,N-2}\big{(}\mathcal{S}_{N-2,N-1}+\mathcal{T}_{N-2,N-1}\big{)}f_{N-2,N-2}+\mathcal{B}^T_{N-2,N-2}\pi_{N-2,N-1}+\rho_{N-2,N-2}\\[1mm]
\hphantom{\theta_{N-1}=}+\sum_{i,j=1}^p\delta_{N-2}^{ij}(\mathcal{D}^i_{N-2,N-2})^T\big{(}S_{N-2,N-1}+T_{N-2,N-1}\big{)}d^j_{N-2,N-2}.
\end{array}
\right.
\end{eqnarray*}
The following argument is similar to that between (\ref{Ran-(N-1)}) and (\ref{v-open-(N-1)-3}).
Note that $(\Phi, v^{t,x})$ is a mixed equilibrium solution and $X^{t,x,*}$ is given in (\ref{equili-state-closed}). Because $\Phi_{N-2}X^{t,x,*}_{N-2} +v^{t,x}_{N-2}$ satisfies (\ref{stationary-condition-closed-(N-2)}), we have from Lemma \ref{Lemma-matrix-equation} that (\ref{stationary-condition-closed-(N-2)}) is equivalent to
\begin{eqnarray}\label{Ran-(N-2)}
\mathcal{L}_{N-2}X^{t,x,*}_{N-2}+\theta_{N-2}\in \mbox{Ran}\big{(}\mathcal{O}_{N-2}\big{)},
\end{eqnarray}
and for some $\eta_{N-2}\in \mathbb{R}^m$,
\begin{eqnarray*}
\Phi_{N-2}X^{t,x,*}_{N-2} +v^{t,x}_{N-2}=-\mathcal{O}^{\dagger}_{N-2} \mathcal{L}_{N-2}X^{t,x,*}_{N-2}-\mathcal{O}^{\dagger}_{N-2}\theta_{N-2}
+\big{(}I-\mathcal{O}^{\dagger}_{N-2}\mathcal{O}_{N-2}\big{)}\eta_{N-2},
\end{eqnarray*}
or equivalently,
\begin{eqnarray*}\label{Phi-Gamma-(N-2)}
v^{t,x}_{N-2}=-(\mathcal{O}^{\dagger}_{N-2} \mathcal{L}_{N-2}+\Phi_{N-2})X^{t,x,*}_{N-2}-\mathcal{O}^{\dagger}_{N-2} \theta_{N-2}+\big{(}I-\mathcal{O}^{\dagger}_{N-2}\mathcal{O}_{N-2}\big{)}\eta_{N-2}.
\end{eqnarray*}
If we replace $v^{t,x}_{N-1}, v^{t,x}_{N-2}$ of $(\Phi, v^{t,x})$ by (\ref{v-open-(N-1)-3}) and
\begin{eqnarray}\label{Phi-Gamma-(N-2)-2}
v^{t,x}_{N-2}=-(\mathcal{O}^{\dagger}_{N-2} \mathcal{L}_{N-2}+\Phi_{N-2})X^{t,x,*}_{N-2}-\mathcal{O}^{\dagger}_{N-2} \theta_{N-2},
\end{eqnarray}
then the new pair $(\Phi, v^{t,x})$ is also a mixed equilibrium solution.

By repeating the procedure between (\ref{Z-(N-1)}) and (\ref{Phi-Gamma-(N-2)}), we have the properties b) and c).

\emph{ii)$\Rightarrow$i)}. For $k\in \mathbb{T}$, (\ref{convex-closed-2}) and $\mathbb{O}_k\succeq 0$, we have
\begin{eqnarray*}
%
\inf_{\bar{u}_k\in L^2_\mathcal{F}(k; \mathbb{R}^m)}\widetilde{J}(k,0;\bar{u}_k)= \inf_{\bar{u}_k\in L^2_\mathcal{F}(k; \mathbb{R}^m)}\big{[}\bar{u}_k^T{\mathbb{O}}_k\bar{u}_k\big{]}\geq 0,
\end{eqnarray*}
which implies (\ref{convex-closed}). Furthermore, based on Lemma \ref{Lemma-matrix-equation} and by reversing the procedure of \emph{i)$\Rightarrow$ii)}, we can assert that $(\Phi, v^{t,x})$ with $v^{t,x}$ given in (\ref{v}) is a mixed equilibrium solution of Problem (LQ)$_{tx}$.  \hfill $\square$


%

\section{Proof of Theorem \ref{Sec3:Theorem-N-S}}\label{appendix:Theorem-N-S}

\emph{i)$\Leftrightarrow$ii).} Note that (\ref{station-strategy}) is equivalent to
$\widetilde{\mathbb{O}}_k \widetilde{\mathbb{O}}_k^\dagger\big{(}\widetilde{\mathbb{L}}_kX^{t,x,*}_k+\widetilde{\theta}_k \big{)}=\widetilde{\mathbb{L}}_kX^{t,x,*}_k+\widetilde{\theta}_k$, $k\in \mathbb{T}_t$.
Letting $k=t$ and taking different $xs$, we have
$\widetilde{\mathbb{O}}_t \widetilde{\mathbb{O}}_t^\dagger\widetilde{\mathbb{L}}_t=\widetilde{\mathbb{L}}_t$, $\widetilde{\mathbb{O}}_t \widetilde{\mathbb{O}}_t^\dagger\widetilde{\theta}_t=\widetilde{\theta}_t$.
Because for any $(t,x)$ with $t\in \mathbb{T}$ and $x\in l^2_\mathcal{F}(t;\mathbb{R}^n)$ Problem (LQ)$_{tx}$ admits a linear feedback equilibrium strategy, we must have the solvability of (\ref{S-all-feedback})-(\ref{pi-all-feedback}). Furthermore, from the solvability of (\ref{S-all-feedback})-(\ref{pi-all-feedback}), it is not hard to confirm the existence of a linear feedback equilibrium strategy.

\emph{ii)$\Rightarrow$iii).} Let
$\Phi=\{-\widetilde{\mathbb{O}}^{\dagger}_{k} \widetilde{\mathbb{L}}_{k},~k\in \mathbb{T}\}$,
$v=\{-\widetilde{\mathbb{O}}^{\dagger}_{k} \widetilde{\theta}_{k},~k\in \mathbb{T}\}$.
Then, for any $(t,x)$ with $t\in \mathbb{T}$ and $x\in l^2(t;\mathbb{R}^n)$, $(\Phi, v)|_{\mathbb{T}_t}$ is a linear feedback equilibrium strategy.

\emph{iii)$\Rightarrow$iv).} Let $\psi=(\Phi,v)$. Then, this $\psi$ satisfies the property of iv).

\emph{iv)$\Rightarrow$ii).}  We adopt a backward procedure to prove ii). First, letting $t=N-1$, then (\ref{defi-closed-loop}) reads as
\begin{eqnarray}\label{Sec4-Necessary-1}
J\big{(}N-1, X_{N-1}^{N-1,x,*}; \psi_{N-1}(X^{N-1,\psi}_{N-1})\big{)}\leq J\big{(}N-1, X_{N-1}^{N-1,x,*}; u_{N-1}\big{)},~~ \forall u_{N-1}\in l^2_\mathcal{F}(N-1; \mathbb{R}^m).
\end{eqnarray}
Noting $X_{N-1}^{N-1,x,*}=X^{N-1,\psi}_{N-1}=x$, it follows that
\begin{eqnarray*}
&&J(N-1,x;u_{N-1})\\[1mm]
&&=x^T\big{[}\mathcal{Q}_{N-1,N-1}+\mathcal{A}^T_{N-1,N-1}\mathcal{G}_{N-1}\mathcal{A}_{N-1,N-1} +2\mathcal{A}^T_{N-1,N-1}F_{N-1}\\[1mm]
&&\hphantom{=}+\sum_{i,j=1}^p({\mathcal{C}^i}_{N-1,N-1})^T{G}_{N-1}{\mathcal{C}^j}_{N-1,N-1}\big{]}x+2\big{[}x^T\big{(}\mathcal{A}_{N-1,N-1}^T\mathcal{G}_{N-1}\mathcal{B}_{N-1,N-1}\\[1mm]
&&\hphantom{=}+\sum_{i,j=1}^p(\mathcal{C}^i_{N-1,N-1})^T {G}_{N-1}\mathcal{D}^j_{N-1,N-1} +F_{N-1}^T\mathcal{B}_{N-1,N-1}\big{)} +f^T_{N-1,N-1}\mathcal{G}_{N-1}\mathcal{B}_{N-1,N-1}\\[1mm]
&&\hphantom{=}+\sum_{i,j=1}^p(d^i_{N-1,N-1})^T{G}_{N-1}\mathcal{D}^j_{N-1,N-1} +\rho^T_{N-1,N-1}+g_{N-1}^T\mathcal{B}_{N-1,N-1}\big{]}u_{N-1} \\[1mm]
&&\hphantom{=}+u_{N-1}^T\big{[}\mathcal{R}_{N-1,N-1}+\mathcal{B}^T_{N-1,N-1}\mathcal{G}_{N-1}\mathcal{B}_{N-1,N-1}+\sum_{i,j=1}^p(\mathcal{D}^i_{N-1,N-1})^T {G}_{N-1}\mathcal{D}^j_{N-1,N-1} \big{]}u_{N-1}\\[1mm]
&&\hphantom{=}+2x^T\big{[}q_{N-1,N-1}+\mathcal{A}^T_{N-1,N-1}\mathcal{G}_{N-1}f_{N-1,N-1}+\mathcal{A}_{N-1,N-1}^Tg_{N-1}\\[1mm]
&&\hphantom{=}+\sum_{i,j=1}^p (\mathcal{C}^i_{N-1,N-1})^T{G}_{N-1}d^j_{N-1,N-1} \big{]}+f_{N-1}^T\mathcal{G}_{N-1}f_{N-1}+2g_{N-1}^Tf_{N-1,N-1}+\sum_{i,j=1}^p(d_{N-1}^i)^T{G}_{N-1}d^j_{N-1}\\[1mm]
&&\geq J\big{(}N-1, x; \psi_{N-1}(x)\big{)}> -\infty.
\end{eqnarray*}
The last inequality occurs because $\psi\in \mathbb{F}_{\mathbb{T}}$. Based on a lemma in \cite{Ni-Zhang-Li}, we have
\begin{eqnarray*}
%
\widetilde{\mathcal{O}}_{N-1}\succeq  0,~~~
\widetilde{\mathcal{O}}_{N-1}\widetilde{\mathcal{O}}_{N-1}^\dagger \widetilde{\mathcal{L}}_{N-1}=\widetilde{\mathcal{L}}_{N-1},~~~
\widetilde{\mathcal{O}}_{N-1}\widetilde{\mathcal{O}}_{N-1}^\dagger \widetilde{\theta}_{N-1}=\widetilde{\theta}_{N-1},
\end{eqnarray*}
and for any $u_{N-1}\in L^2_{\mathcal{F}}(N-1; \mathbb{R}^m)$,
\begin{eqnarray}\label{Sec4-Necessary-2}
&&J(N-1,x;u_{N-1})\nonumber \\[1mm]
&&=(\widetilde{\mathcal{O}}_{N-1}u_{N-1}+\widetilde{\mathcal{L}}_{N-1}x+\widetilde{\theta}_{N-1})^T \widetilde{\mathcal{O}}^\dagger_{N-1}(\widetilde{\mathcal{O}}_{N-1}u_{N-1}+\widetilde{\mathcal{L}}_{N-1}x+\widetilde{\theta}_{N-1})\nonumber \\[1mm]
&&\hphantom{=}+x^T\mathcal{P}_{N-1,N-1}x+2x^T\big{[}-\widetilde{\mathcal{L}}_{N-1}^T\widetilde{\mathcal{O}}_{N-1}\widetilde{\theta}_{N-1} +q_{N-1,N-1}+\mathcal{A}^T_{N-1,N-1}\mathcal{G}_{N-1}f_{N-1,N-1}\nonumber \\[1mm]
&&\hphantom{=}+\sum_{i,j=1}^p(\mathcal{C}^i_{N-1,N-1})^T{G}_{N-1}d^j_{N-1,N-1} +\mathcal{A}_{N-1,N-1}^Tg_{N-1}\big{]}+f_{N-1}^T\mathcal{G}_{N-1}f_{N-1}\nonumber \\[1mm]
&&\hphantom{=}+\sum_{i,j=1}^p(d^i_{N-1})^T{G}_{N-1}d^j_{N-1}+2g_{N-1}^Tf_{N-1,N-1}-\widetilde{\theta}_{N-1}^T\widetilde{\mathcal{O}}_{N-1}\widetilde{\theta}_{N-1}\nonumber \\[1mm]
&&\geq J(N-1,x;\widetilde{u}_{N-1}(x)),
\end{eqnarray}
where
\begin{eqnarray*}
\left\{
\begin{array}{l}
\widetilde{\mathcal{O}}_{N-1}=\mathcal{R}_{N-1,N-1}+\mathcal{B}^T_{N-1,N-1}\mathcal{G}_{N-1}\mathcal{B}_{N-1,N-1}+ \sum_{i,j=1}^p(\mathcal{D}^i_{N-1,N-1})^T{G}_{N-1}\mathcal{D}^j_{N-1,N-1},\\[1mm]
\widetilde{\mathcal{L}}_{N-1}=\mathcal{B}_{N-1,N-1}^T\mathcal{G}_{N-1}\mathcal{A}_{N-1,N-1}+\sum_{i,j=1}^p(\mathcal{D}^i_{N-1,N-1})^T{G}_{N-1}\mathcal{C}^j_{N-1,N-1} +\mathcal{B}_{N-1,N-1}^TF_{N-1},\\[1mm]
\widetilde{\theta}_{N-1}=\mathcal{B}_{N-1,N-1}^T\mathcal{G}_{N-1}f_{N-1,N-1}+ \sum_{i,j=1}^p(\mathcal{D}^i_{N-1,N-1})^T{G}_{N-1}d^j_{N-1,N-1}+\rho_{N-1,N-1}+\mathcal{B}_{N-1,N-1}g_{N-1},
\end{array}
\right.
\end{eqnarray*}
%
$$
\mathcal{P}_{N-1,N-1}=\mathcal{Q}_{N-1,N-1}+\mathcal{A}^T_{N-1,N-1}\mathcal{G}_{N-1}\mathcal{A}_{N-1,N-1}
+\sum_{i,j=1}^p({\mathcal{C}}^i_{N-1,N-1})^T{G}_{N-1}{\mathcal{C}}^j_{N-1,N-1}
-\widetilde{\mathcal{L}}_{N-1}^T\widetilde{\mathcal{O}}^{\dagger}_{N-1}\widetilde{\mathcal{L}}_{N-1},
$$
and
$\widetilde{u}_{N-1}(x)=-\widetilde{\mathcal{O}}_{N-1}^\dagger\widetilde{\mathcal{L}}_{N-1}x
-\widetilde{\mathcal{O}}_{N-1}^\dagger\widetilde{\theta}_{N-1}$. From this and (\ref{Sec4-Necessary-1}), (\ref{Sec4-Necessary-2}), one can select $\psi_{N-1}$ as $(-\widetilde{\mathcal{O}}_{N-1}^\dagger\widetilde{\mathcal{L}}_{N-1}, -\widetilde{\mathcal{O}}_{N-1}^\dagger\widetilde{\theta}_{N-1})$, i.e.,
$\psi_{N-1}(x)=-\widetilde{\mathcal{O}}_{N-1}^\dagger\widetilde{\mathcal{L}}_{N-1}x
-\widetilde{\mathcal{O}}_{N-1}^\dagger\widetilde{\theta}_{N-1}$.

Assume we have obtained $\psi_{\ell}=(\widetilde{\Phi}_\ell, \widetilde{v}_\ell), \ell\in \mathbb{T}_{k+1}$, namely, $\psi(z)=\widetilde{\Phi}_{\ell}z+\widetilde{v}_\ell$, with $(\Phi_{\ell}, \widetilde{v}_\ell)\in \mathbb{R}^{m\times n}\times \mathbb{R}^{m}$. Let us derive the expression of $\psi_{k}$. Now, consider Problem (LQ) for the initial pair $(k,x)$. 
By adding to and subtracting three terms
\begin{eqnarray*}
&&\hspace{-1.5em}\sum_{\ell=k+1}^{N-1}\mathbb{E}_k\Big{[}\big{(}X^{k,u_k,\psi}_{\ell+1}-\mathbb{E}_kX^{k,u_k,\psi}_{\ell+1}\big{)}^TP_{k,\ell+1} \big{(}X^{k,u_k,\psi}_{\ell+1}-\mathbb{E}_kX^{k,u_k,\psi}_{\ell+1} \big{)}\\
&&\hspace{-1.5em}\hphantom{\sum_{\ell=k+1}^{N-1}\mathbb{E}_k\Big{[}}-\big{(}X^{k,u_k,\psi}_{\ell}-\mathbb{E}_kX^{k,u_k,\psi}_{\ell}\big{)}^T P_{k,\ell}\big{(}X^{k,u_k,\psi}_{\ell}-\mathbb{E}_kX^{k,u_k,\psi}_{\ell} \big{)}\Big{]}, \\
&&\hspace{-1.5em}\sum_{\ell=k+1}^{N-1}\Big{[}\big{(}\mathbb{E}_kX^{k,u_k,\psi}_{\ell+1}\big{)}^T\mathcal{P}_{k,\ell+1}\mathbb{E}_kX^{k,u_k,\psi}_{\ell+1}- \big{(}\mathbb{E}_kX^{k,u_k,\psi}_{\ell}\big{)}^T\mathcal{P}_{k,\ell}\mathbb{E}_kX^{k,u_k,\psi}_{}\Big{]}, \\
&&\hspace{-1.5em}2\sum_{\ell=k+1}^{N-1}\mathbb{E}_k\big{[}\sigma_{k,\ell+1}^TX^{k,u_k,\psi}_{k,\ell+1}- \sigma_{k,\ell}^TX^{k,u_k,\psi}_{\ell}\big{]}
\end{eqnarray*}
from $J\big{(}k, x; (u_k,(\psi\cdot X^{k,u_k,\psi})|_{\mathbb{T}_{k+1}})\big{)}$, we have
\begin{eqnarray*}\label{Sec4-Necessary-8}
&&\hspace{-1.5em}J\big{(}k, x; (u_k,(\psi\cdot X^{k,u_k,\psi})|_{\mathbb{T}_{k+1}})\big{)}\\[1mm]
&&\hspace{-1.5em}=\mathbb{E}_k\big{[}x^T\mathcal{Q}_{k,k}x + u_k^T\mathcal{R}_{k,k}u_k+2q_{k,k}^T x+ 2\rho_{k,k}^T u_k\big{]}+\sum_{\ell=k+1}^{N-1}\dbE_k\Big{[}(X^{k,u_k,\psi}_\ell-\mathbb{E}_kX^{k,u_k,\psi}_\ell)^T\Big{(}Q_{k,\ell} +\widetilde{\Phi}_{\ell}^TR_{k,\ell}\widetilde{\Phi}_{\ell}\nonumber \\[1mm]
&&\hspace{-1.5em}\hphantom{=} +(A_{k,\ell}+ B_{k,\ell}\widetilde{\Phi}_{\ell})^TP_{k,\ell+1}(A_{k,\ell}+ B_{k,\ell}\widetilde{\Phi}_{\ell})+\sum_{i,j=1}^p(C^i_{k,\ell}+ D^i_{k,\ell}\widetilde{\Phi}_{\ell})^TP_{k,\ell+1}(C^j_{k,\ell}+ D^j_{k,\ell}\widetilde{\Phi}_{\ell})-P_{k,\ell}\Big{)} \\
&&\hspace{-1.5em}\hphantom{=}\times(X^{k,u_k,\psi}_\ell-\mathbb{E}_kX^{k,u_k,\psi}_\ell)+(\mathbb{E}_kX^{k,u_k,\psi}_\ell)^T\Big{(}
 {\mathcal{Q}}_{k,\ell}+\widetilde{\Phi}_{\ell}^T{\mathcal{R}}_{k,\ell}\widetilde{\Phi}_{\ell}+(\mathcal{A}_{k,\ell}+\mathcal{B}_{k,\ell}\widetilde{\Phi}_{\ell})^T\mathcal{P}_{k,\ell+1} (\mathcal{A}_{k,\ell}+\mathcal{B}_{k,\ell}\widetilde{\Phi}_{\ell})\\[1mm]
&&\hspace{-1.5em}\hphantom{=}+\sum_{i=1}^p(C^i_{k,\ell}+{D}^i_{k,\ell}\widetilde{\Phi}_{\ell})^T{P}_{k,\ell+1} ({C}^j_{k,\ell}+{D}^j_{k,\ell}\widetilde{\Phi}_{\ell})-\mathcal{P}_{k,\ell} \Big{)}\mathbb{E}_kX^{k,u_k,\psi}_\ell\nonumber \\[1mm]
&&\hspace{-1.5em}\hphantom{=}+2\Big{(}q_{k,\ell}+\widetilde{\Phi}_{\ell}^T\mathcal{R}_{k,\ell}\widetilde{v}_{\ell}+\widetilde{\Phi}_{\ell}^T\rho_{k,\ell} +\sum_{i,j=1}^p(C^i_{k,\ell}+D^i_{k,\ell}\widetilde{\Phi}_{\ell})^TP_{k,\ell+1}(D^j_{k,\ell}\widetilde{v}_{\ell}+d^j_{k,\ell})\\[1mm]
&&\hspace{-1.5em}\hphantom{=}+(\mathcal{A}_{k,\ell}+\mathcal{B}_{k,\ell}\widetilde{\Phi}_{\ell})^T\mathcal{P}_{k,\ell+1}(\mathcal{B}_{k,\ell}\widetilde{v}_{k,\ell}+f_{k,\ell})
+(\mathcal{A}_{k,\ell}+\mathcal{B}_{k,\ell}\widetilde{\Phi}_{\ell})^T\sigma_{k,\ell+1}-\sigma_{k,\ell}\Big{)}^T X^{k,u_k,\psi}_\ell\\[1mm]
&&\hspace{-1.5em}\hphantom{=}+2\rho_{k,\ell}^T\widetilde{v}_\ell+\widetilde{v}^T\mathcal{R}_{t,k}\widetilde{v}_\ell+\sum_{i=1}^p(D^i_{t,k}\widetilde{v}_\ell+d^i_{k,\ell})^T P_{k,\ell} (D^j_{k,\ell}\widetilde{v}_\ell+d^j_{k,\ell})\\[1mm]
&&\hspace{-1.5em}\hphantom{=}+(\mathcal{B}_{t,k}\widetilde{v}_\ell+f_{k,\ell})^T\mathcal{P}_{k,\ell+1}(\mathcal{B}_{k,\ell}\widetilde{v}_\ell+f_{k,\ell}) +\sigma_{k,\ell+1}^T(\mathcal{B}_{k,\ell}\widetilde{v}_\ell+f_{k,\ell})\Big{]}\nonumber\\[1mm]
&&\hspace{-1.5em}\hphantom{=}+2x^TF_k^T\mathbb{E}_kX^{k,u_k,\psi}_N+\mathbb{E}_k\Big{[}(X^{k,u_k,\psi}_{k+1}-\mathbb{E}_kX^{k,u_k,\psi}_{k+1})^T P_{k,k+1}(X^{k,u_k,\psi}_{k+1}-\mathbb{E}_kX^{k,u_k,\psi}_{k+1})\Big{]}
\nonumber\\[1mm]
&&\hspace{-1.5em}\hphantom{=}+(\mathbb{E}_kX^{k,u_k,\psi}_{k+1})^T {\mathcal{P}}_{k,k+1}\mathbb{E}_kX^{k,u_k,\psi}_{k+1}+2\sigma_{k,k+1}^T \mathbb{E}_k X^{k,u_k,\psi}_{k+1}.
\end{eqnarray*}
Let
\begin{eqnarray*}
\left\{
\begin{array}{l}
P_{k,\ell}=Q_{k,\ell}+\Phi_\ell^TR_{k,\ell}\Phi_\ell+(A_{k,\ell}+ B_{k,\ell}\Phi_\ell)^TP_{k,\ell+1}(A_{k,\ell}+ B_{k,\ell}\Phi_\ell)\\
\hphantom{P_{k,\ell}=}+\sum_{i,j=1}^p(C^i_{k,\ell}+ D^i_{k,\ell}\Phi_\ell)^TP_{k,\ell+1}(C^j_{k,\ell}+ D^j_{k,\ell}\Phi_\ell),\\[1mm]
\mathcal{P}_{k,\ell}= {\mathcal{Q}}_{k,\ell}+\Phi_\ell^T{\mathcal{R}}_{k,\ell}\Phi_{k,\ell}+(\mathcal{A}_{k,\ell}+\mathcal{B}_{k,\ell}\Phi_{k,\ell})^T\mathcal{P}_{k,\ell+1} (\mathcal{A}_{k,\ell}+\mathcal{B}_{k,\ell}\Phi_\ell)\\[1mm]
\hphantom{\mathcal{P}_{k,\ell}= }+\sum_{i,j=1}^p(C^i_{k,\ell}+{D}^i_{k,\ell}\Phi_\ell)^T{P}_{k,\ell+1} ({C}^j_{k,\ell}+{D}^j_{k,\ell}\Phi_\ell),\\[1mm]
\sigma_{k,\ell}=q_{k,\ell}+\Phi_\ell^T\mathcal{R}_{k,\ell}v_\ell+\Phi_\ell^T\rho_{k,\ell} +\sum_{i=1}^p(C^i_{k,\ell}+D^i_{k,\ell}\Phi_\ell)^TP_{k,\ell+1}(D^j_{k,\ell}v_\ell+d^j_{k,\ell})\\[1mm]
\hphantom{\sigma_{k,\ell}=}+(\mathcal{A}_{k,\ell}+\mathcal{B}_{k,\ell}\Phi_\ell)^T\mathcal{P}_{k,\ell+1}(\mathcal{B}_{k,\ell}\Phi_{k,\ell}+f_{k,\ell})
+(\mathcal{A}_{k,\ell}+\mathcal{B}_{k,\ell}v_\ell)^T\sigma_{k,\ell+1},\\[1mm]
P_{k,N}=G_k,~~\mathcal{P}_{k,N}=\mathcal{G}_{k},~~\sigma_{k,N}=g_k,\quad \ell\in \mathbb{T}_{k+1},
\end{array}
\right.
\end{eqnarray*}
and
\begin{eqnarray*}
&&\hspace{-3em}\gamma_{k,k+1}=\sum_{\ell=k+1}^{N-1}\Big{[}2\rho_{k,\ell}^Tv_\ell+v_\ell^T\mathcal{R}_{k,\ell}v_\ell+\sum_{i,j=1}^p(D^i_{k,\ell}v_\ell +d^i_{k,\ell})^TP_{k,\ell}(D^j_{k,\ell}v_\ell+d^j_{k,\ell})\nonumber\\[1mm]
&&\hspace{-3em}\hphantom{\gamma_{k,k+1}=\sum_{\ell=k+1}^{N-1}}+(\mathcal{B}_{k,\ell}v_\ell+f_{k,\ell})^T\mathcal{P}_{k,\ell+1}(\mathcal{B}_{k,\ell}v_\ell+f_{k,\ell}) +\sigma_{k,\ell+1}^T(\mathcal{B}_{k,\ell}v_\ell+f_{k,\ell})\Big{]}\\
&&\hspace{-3em}\hphantom{\gamma_{k,k+1}=}+F_k(\mathcal{B}_{k,N-1}v_{N-1}+f_{k,N-1})+F_k(\mathcal{A}_{k,N-1}+\mathcal{B}_{k,N-1}\Phi_{N-1})(\mathcal{B}_{k,N-2}v_{N-2}+f_{k,N-2}) \\
&&\hspace{-3em}\hphantom{\gamma_{k,k+1}=}+\sum_{\ell=k}^{N-3}F_k(\mathcal{A}_{k,N-1}+\mathcal{B}_{k,N-1}\Phi_{N-1})\cdots (\mathcal{A}_{k,\ell+1}+\mathcal{B}_{k,\ell+1}\Phi_{\ell+1})(\mathcal{B}_{k,\ell}v_{\ell}+f_{k,\ell})\\
&&\hspace{-3em}\hphantom{\gamma_{k,k+1}=}+F_k(\mathcal{A}_{k,N-1}+\mathcal{B}_{k,N-1}\Phi_{N-1})\cdots (\mathcal{A}_{k,k+1}+\mathcal{B}_{k,k+1}\Phi_{k+1})f_{k,k}.
\end{eqnarray*}
Then, it holds that
\begin{eqnarray*}
&&J\big{(}k, x; (u_k,(\psi\cdot X^{k,u_k,\psi})|_{\mathbb{T}_{k+1}})\big{)}\\[1mm]
&&=x^T\big{[}\mathcal{Q}_{k,k}+\mathcal{A}_{k,k}^T\mathcal{P}_{k,k+1}\mathcal{A}_{k,k}
+\sum_{i,j=1}^p(\mathcal{C}^i_{k,k})^TP_{k,k+1}\mathcal{C}^j_{k,k} +2\widetilde{U}_{k+1}^T\mathcal{A}_{k,k}\big{]}x  \\
&&\hphantom{=} +u_k^T\big{[}\mathcal{R}_{k,k}+\mathcal{B}_{k,k}^T\mathcal{P}_{k,k+1}\mathcal{B}_{k,k} +\sum_{i,j=1}^p(\mathcal{D}^i_{k,k})^TP_{k,k+1}\mathcal{D}^j_{k,k} \big{]}u_k\\
&&\hphantom{=}+2x^T\big{[}\mathcal{A}_{k,k}^T\mathcal{P}_{k,k+1}\mathcal{B}_{k,k}
+\sum_{i,j=1}^p(\mathcal{C}^i_{k,k})^TP_{k,k+1}\mathcal{D}^j_{k,k}
+\widetilde{U}^T_{k+1}\mathcal{B}_{k,k} \big{]}u_k \\
&&\hphantom{=}+2u_k^T\big{[}\rho_{k,k}+\mathcal{B}_{k,k}^T\mathcal{P}_{k,k+1}f_{k,k}+ \sum_{i,j=1}^p(\mathcal{D}^i_{k,k})^TP_{k,k+1}d^j_{k,k}+\mathcal{B}^T_{k,k}\sigma_{k,k+1}\big{]} \\
&&\hphantom{=}+2x^T\big{[}q_{k,k}+\mathcal{A}_{k,k}^T\mathcal{P}_{k,k+1}f_{k,k}
+\sum_{i,j=1}^p(\mathcal{C}^i_{k,k})^TP_{k,k+1}d^j_{k,k} +\mathcal{A}_{k,k}^T\sigma_{k,k+1} \big{]} \\
&&\hphantom{=}+f_{k,k}^T\mathcal{P}_{k,k+1}f_{k,k} +\sum_{i,j=1}^p(d^i_{k,k})^TP_{k,k+1}d^j_{k,k}+2\sigma_{k,k+1}^Tf_{k,k}+\gamma_{k,k+1} \\
&&\geq J\big{(}k, x; (\psi\cdot X^{k,\psi})|_{\mathbb{T}_k}\big{)}> -\infty,
\end{eqnarray*}
where
$\widetilde{U}_{k+1}=(\mathcal{A}_{k,k+1}+\mathcal{B}_{k,k+1}
\Phi_{k+1})\cdots(\mathcal{A}_{k,N-1}+\mathcal{B}_{k,N-1}\Phi_{N-1})F_k$.
Then, we have
$\widetilde{\mathcal{O}}_{k}\succeq  0$,
$\widetilde{\mathcal{O}}_{k}\widetilde{\mathcal{O}}_{k}^\dagger \widetilde{\mathcal{L}}_{k}=\widetilde{\mathcal{L}}_{k}$,
$\widetilde{\mathcal{O}}_{k}\widetilde{\mathcal{O}}_{k}^\dagger \widetilde{\theta}_{k}=\widetilde{\theta}_{k}$,
and
\begin{eqnarray*}\label{Sec4-Necessary-12-2}
J\big{(}k,x;(u_k,(\psi\cdot X^{k,u_k,\psi})|_{\mathbb{T}_{k+1}})\big{)}\geq J\big{(}k,x;(\widetilde{u}_{k}(x),(\psi\cdot X^{k,\widetilde{u}_{k}(x),\psi})|_{\mathbb{T}_{k+1}})\big{)},~~\forall u_{k}\in l^2_{\mathcal{F}}(k; \mathbb{R}^m),
\end{eqnarray*}
where
\begin{eqnarray*}\label{Sec4-Necessary-13}
\left\{
\begin{array}{l}
\widetilde{\mathcal{O}}_{k}=\mathcal{R}_{k,k}+\mathcal{B}^T_{k,k}\mathcal{P}_{k,k+1}\mathcal{B}_{k,k}+\sum_{i,j=1}^p(\mathcal{D}^i_{k,k})^T{P}_{k,k+1}\mathcal{D}^j_{k,k},\\[1mm]
\widetilde{\mathcal{L}}_{k}=\mathcal{B}_{k,k}^T\mathcal{P}_{k,k+1}\mathcal{A}_{k,k}+\sum_{i,j=1}^p(\mathcal{D}^i_{k,k})^T{P}_{k,k+1}\mathcal{C}^j_{k,k},\\[1mm]
\widetilde{\theta}_{k}=\mathcal{B}_{k,k}^T\mathcal{P}_{k,k+1}f_{k,k}+\sum_{i,j=1}^p(\mathcal{D}^i_{k,k})^T{P}_{k,k+1}d^j_{k,k} +\rho_{k,k}+\mathcal{B}_{k,k}\sigma_{k,k+1},
\end{array}
\right.
\end{eqnarray*}
and
$\widetilde{u}_{k}(x)=-\widetilde{\mathcal{O}}_{k}^\dagger\widetilde{\mathcal{L}}_{k}x
-\widetilde{\mathcal{O}}_{k}^\dagger\widetilde{\theta}_{k}$.
Hence, $\psi_{k}$ can be selected as $(-\widetilde{\mathcal{O}}_{k}^\dagger\widetilde{\mathcal{L}}_{k}, -\widetilde{\mathcal{O}}_{k}^\dagger\widetilde{\theta}_{k})$, i.e.,
$
\psi_{k}(x)=-\widetilde{\mathcal{O}}_{k}^\dagger\widetilde{\mathcal{L}}_{k}x-\widetilde{\mathcal{O}}_{k}^\dagger\widetilde{\theta}_{k}.
$
Furthermore, $\widetilde{S}_{k,\ell}=P_{k,\ell}, \widetilde{\mathcal{S}}_{k,\ell}=\mathcal{P}_{k,\ell}, \sigma_{k,\ell}=\pi_{k,\ell}, \ell\in \mathbb{T}_{k+1}$, where $(\widetilde{S}_{k,\ell}, \widetilde{\mathcal{S}}_{k,\ell}, \pi_{k,\ell})$ is given in (\ref{S-all-feedback}) and (\ref{pi-all-feedback}).

By the method of induction, we have the solvability of (\ref{S-all-feedback}) (\ref{pi-all-feedback}), and $\{\psi_k=(-\widetilde{\mathcal{O}}_{k}^\dagger\widetilde{\mathcal{L}}_{k}, -\widetilde{\mathcal{O}}_{k}^\dagger\widetilde{\theta}_{k}), k\in \mathbb{T}\}$ is a feedback equilibrium strategy, that is,
$\psi_k(x)=-\widetilde{\mathcal{O}}_{k}^\dagger\widetilde{\mathcal{L}}_{k}x
-\widetilde{\mathcal{O}}_{k}^\dagger\widetilde{\theta}_{k},~k\in \mathbb{T}$.
This completes the proof. \hfill $\square$

\end{document}